\documentclass[preprint,12pt]{elsarticle}

\usepackage{subcaption}
\usepackage{algorithm}
\usepackage{algpseudocode} 
\usepackage{float}
\usepackage[margin=2.5cm]{geometry}
\usepackage{amsmath,amsfonts,amssymb}
\usepackage{mathptmx,mathrsfs,mathtools}
\usepackage{booktabs}
\usepackage{cleveref}
\usepackage{pifont}
\def\tsc#1{\csdef{#1}{\textsc{\lowercase{#1}}\xspace}}
\tsc{WGM}
\tsc{QE}
\tsc{EP}
\tsc{PMS}
\tsc{BEC}
\tsc{DE}

\biboptions{sort&compress}

\usepackage[bordercolor=white,backgroundcolor=gray!30,linecolor=black,colorinlistoftodos]{todonotes}

\begin{document}

\begin{frontmatter}
\let\WriteBookmarks\relax
\def\floatpagepagefraction{1}
\def\textpagefraction{.001}


\title{$Dmsh$: A Multi-Agent Reinforcement Learning Framework for All-Quad Mesh Generation}

\author[1]{Anirudh Kalyan}
\author[2]{Cosmin Anitescu}
\author[2]{Xiaoying Zhuang}
\author[3]{Timon Rabzcuk}
\author[4]{Somdatta Goswami}
\author[1]{Sundararajan Natarajan\corref{cor1}\fnref{cor2}}
\ead{snatarajan@iitm.ac.in}

\address[1]{Department of Mechanical Engineering, Indian Institute of Technology Madras, Chennai-600036, India}
\address[2]{Institute of Continuum Mechanics, Leibniz Universit\"at Hannover, Hannover, Germany}
\address[3]{Institute of Structural Mechanics, Bauhaus-Universit\"at Weimar, Weimar, Germany}
\address[4]{Department of Civil and Systems Engineering, Johns Hopkins University, Baltimore, Maryland, USA}


\cortext[cor1]{Pandurangan Faculty Fellow}  
\fntext[cor2]{Corresponding author}

\begin{abstract}
Generating high-quality meshes for arbitrary geometries remains a fundamental bottleneck in computational engineering, often demanding heuristic tuning and semi-manual workflows. In this paper, we introduce $Dmsh$, a first fully automated reinforcement learning pipeline that unifies geometric decomposition and quadrilateral mesh generation within a single learning-based framework. $Dmsh$ decomposes the problem through three coordinated agents handling topology simplification, geometric regularization, and mesh generation. The meshing process is formulated as a Markov Decision Process and solved using a parametric Soft Actor-Critic architecture with decoupled critics, enabling efficient exploration of a hybrid discrete–continuous action space. A curriculum learning strategy ensures scalability from simple domains to highly complex geometries, suppressing seed variance. By design, the recursive decomposition enables parallel meshing of subregions, yielding globally conforming all-quadrilateral meshes without post hoc correction. Across a wide range of benchmarks, $Dmsh$ consistently outperforms existing methods in automation, robustness, and mesh quality, establishing a new paradigm for learning-based mesh generation.
\end{abstract}

\begin{graphicalabstract}
\centering
\includegraphics[scale=0.42]{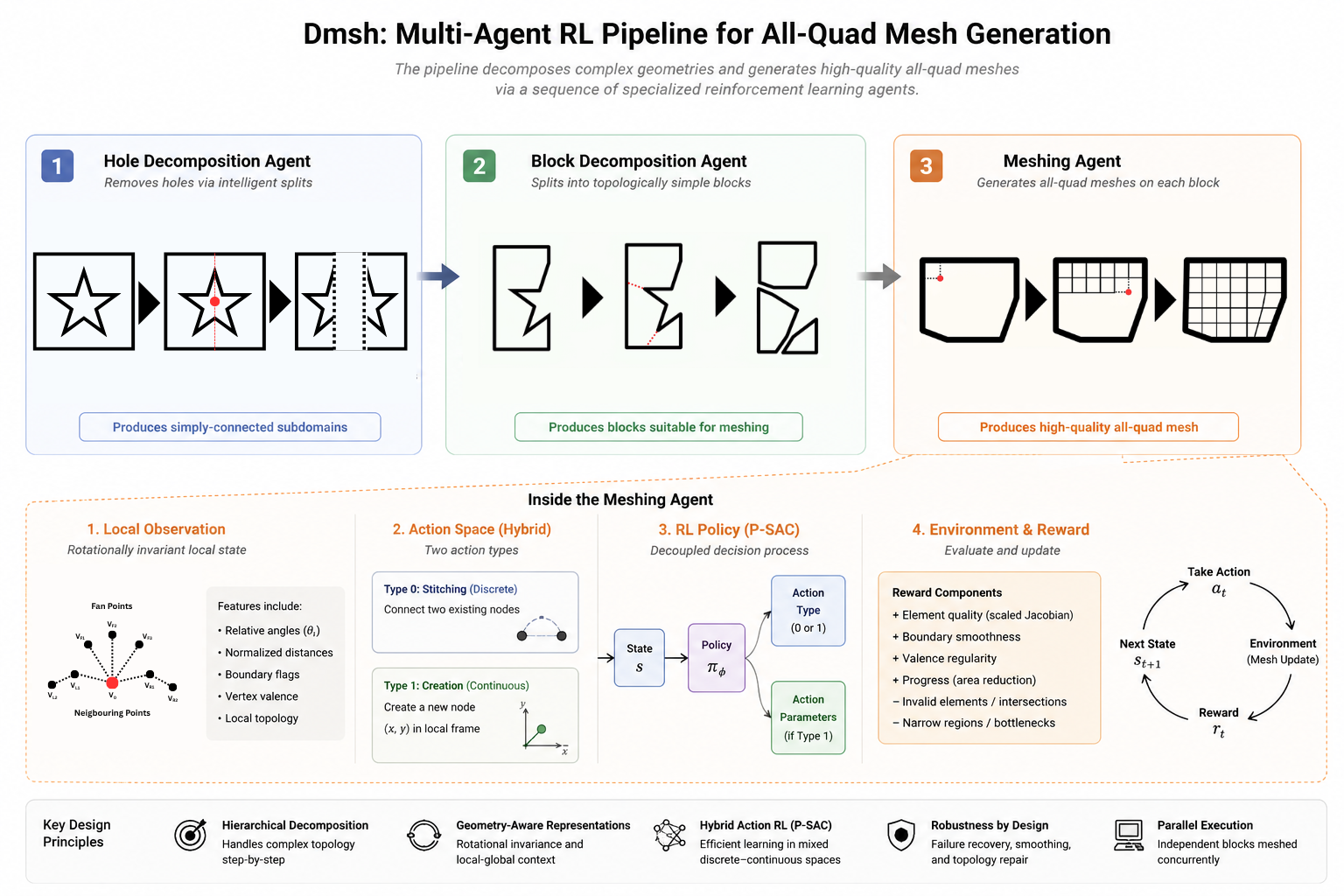}
\end{graphicalabstract}

\begin{highlights}
\item[\ding{113}] Fully automated 3-agent reinforcement learning pipeline for all-quad meshes
\item[\ding{113}]  Unified block/hole decomposition and meshing within a learning framework
\item[\ding{113}]  Computer vision-based state representation for geometry-aware decisions
\item[\ding{113}]  Multi-agent SAC formulation for hybrid discrete–continuous actions
\item[\ding{113}]  Validated on diverse 2D and 3D geometries with improved mesh quality
\end{highlights}


\end{frontmatter}

\section{Introduction}

Modern engineering design is fundamentally constrained by simulation. Every structural analysis, fluid dynamics study, or thermal evaluation that a designer wants to run must wait on a solver, and solvers are slow. The dream of \textit{simulation at the speed of ideation}, sometimes called hardware iteration at the speed of software, has become a defining ambition of computational engineering. Pursuing it, the community has largely converged on two complementary strategies. 

The first, which we term the \textit{top-down} approach, seeks to replace the entire solver stack with neural-network-based surrogates. Neural Operators~\citep{li2021fourier,Lu_2021,goswami2023physics}, which learn mapping between infinite-dimensional functional spaces, promise speedups of several orders of magnitude over classical solvers, by learning the family of partial differential equations. These approaches have achieved notable successes in design optimization, materials modeling, failure prediction, and weather forecasting~\citep{pathak2023fourcastnet}. Yet for many standard engineering applications, these methods still face challenges in robustness, extrapolation, and reliability~\citep{mcgreivy2024weak,roy2025best}, particularly when high accuracy is required across broad operating regimes. 

The second strategy, which we term the \textit{bottom-up} approach, aims to accelerate the existing simulation pipeline that engineers already rely on by targeting its most labor-intensive stages with deep learning and AI \citep{WANG2025118319}, while still being able to obtain solutions with established error guarantees. This is not a concession; it is a recognition that the highest-fidelity simulations in industry will depend on classical solvers for the foreseeable future, and that the greatest leverage lies in removing the human bottlenecks surrounding them. Chief among these is mesh generation.

Meshing is, in a precise sense, the problem that gates everything else. The NASA CFD Vision 2030 study~\citep{slotnick2014cfd} estimated that mesh generation alone can consume the majority of total engineering time in a simulation workflow, more than the solve itself.  While mature software exists for triangular and tetrahedral meshes, all-quadrilateral (2D) and all-hexahedral (3D) meshes, which offer superior convergence, lower numerical diffusion, and better alignment with flow features, still require substantial manual intervention: geometry clean-up, domain decomposition, and careful local adjustment~\citep{owen1998survey, bommes2013quad}. A single complex geometry can take an experienced engineer days to mesh properly. This is the problem we address.

Considerable research has been devoted to automating quad mesh generation through algorithmic approaches. Classical methods include advancing-front techniques such as paving~\citep{blacker1991paving}, which grow the mesh inward from the boundary by sequentially placing elements, and integer-grid map methods~\citep{bommes2013quad}, which formulate the problem as a global optimization over a cross-field. While effective for many geometries, these methods can struggle with complex domains containing sharp features, holes, or highly non-convex regions, often requiring manual decomposition of the domain into simpler blocks before meshing can proceed.

We introduce \textsc{Dmsh}, the first fully automated reinforcement learning pipeline that unifies geometric decomposition and all-quadrilateral mesh generation within a single learning-based framework. Recently, there has been a growing interest in exploring machine learning approaches as a way of automating mesh generation. Zhang et al.~\citep{zhang2020meshingnet} proposed MeshingNet, a deep-learning framework that predicts element sizing fields to guide conventional mesh generators. Comprehensive surveys by Li et al.~\citep{lei2024intelligent} and Owen et al.~\citep{owen2025survey} catalogue the growing body of work applying neural networks, generative models, and reinforcement learning to various stages of the meshing pipeline, from geometry preparation and defeaturing to element placement and quality assessment. Of particular relevance to the present work are two contributions that apply reinforcement learning directly to quad mesh generation via the advancing-front paradigm. Pan et al.~\citep{pan2023reinforcement} formulated quad mesh generation as a Markov Decision Process and trained a Soft Actor-Critic agent to sequentially place quadrilateral elements, demonstrating that reinforcement learning can produce meshes competitive with commercial software without manual intervention. Tong et al.~\citep{tong2023srl} combined supervised pre-training with reinforcement learning fine-tuning to guide an advancing-front method, achieving high accuracy in reproducing commercial-quality meshes for planar domains with sharp features and boundary layers. However, both single-agent formulations are limited in their ability to handle sharp corners robustly, scale poorly to large boundaries, and do not address geometries containing holes.

$Dmsh$ addresses all these limitations by introducing an end-to-end pipeline of three reinforcement learning agents that work in concert to decompose and mesh arbitrary geometries. The meshing agent extends the formulation of Pan et al.\ with a Parametric Soft Actor-Critic architecture tailored to the hybrid discrete-continuous action space of the problem. Two additional vision-based agents operate upstream: a block decomposition agent that systematically partitions the domain to eliminate sharp corners and reduce boundary complexity, and a hole decomposition agent that splits geometries containing internal voids into regular sub-regions. Each sub-region is then meshed independently, and the results are stitched together to produce the final all-quad mesh. This recursive decomposition strategy not only improves element quality but also enables large-scale parallelization, as the lightweight meshing agent can process all sub-regions concurrently. As long as each subdomain contains an even number of points, the result is guaranteed to be an all-quad mesh. To our knowledge, this is the first instance of computer vision being used to tackle the problem of block decomposition. 

The remainder of this paper is organized as follows. Section~2 reviews the foundations of reinforcement learning and the Soft Actor-Critic algorithm, the foundation of our approach. Section~3 describes the meshing agent's environment, the Parametric SAC architecture, and the curriculum learning strategy. Sections~4 and~5 detail the block and hole decomposition agents, respectively. Section~6 presents results in both 2D and 2.5D (surface meshing) and compares them against current state-of-the-art methods.

\section{Overview of Reinforcement Learning Methods}

The meshing pipeline developed in this work formulates the generation of all-quadrilateral meshes as a sequential decision-making problem solved by reinforcement learning (RL) agents. Following the atomic idea established by Zeng et al.~\citep{zeng1993}, the fundamental task is recast as a Markov Decision Process by Pan et al. \citep{pan2023reinforcement} in which, at each step, the agent either connects two existing boundary nodes or places a new node to form a quadrilateral element. The resulting element is evaluated by the environment, a reward is assigned based on its geometric quality, and the updated boundary is recursively passed back to the agent until the domain is fully discretized. Crucially, the agent receives only a partial, rotationally invariant observation of the boundary at each step, which promotes geometry-agnostic behaviour and generalization across diverse domains.

While this single-agent formulation can produce an all-quad mesh on its own, provided the boundary has an even number of nodes, as mentioned in the earlier sections, it suffers from three practical limitations. As discussed earlier, the single agent formulation is not designed to handle holes empirically, and suffers from element quality degredation around convex boundaries. Thirdly, since geometric complexity accumulates and the strict convexity checks imposed by the environment become increasingly difficult to satisfy, the agent's error rate increases, causing a longer runtime. To address these issues, two additional vision-based agents are introduced upstream of the mesher: a block decomposition agent, and a hole decomposition agent.  

All three agents in the pipeline are trained using variants of the Soft Actor-Critic (SAC) algorithm~\citep{haarnoja2018soft}, an off-policy actor-critic method well suited to continuous and hybrid action spaces. The remainder of this section introduces the formal foundations of reinforcement learning and then presents the SAC framework in detail. The later part of the background section then discusses our modifications that better suit the underlying problem at hand. 

\subsection{Fundamentals of Reinforcement Learning}
\label{subsec:rl_fundamentals}

Reinforcement learning problems are typically formalized within the framework of a Markov Decision Process (MDP)~\citep{bellman1957dynamic,sutton2018reinforcement}, defined by the tuple $(S, A, P, R, \gamma)$. Here $S$ denotes the set of states, $s$ the agent may occupy, $A$ is the set of actions, $a$ available to it, $P(s'|s,a)$ specifies the probability of transitioning to state $s'$ after taking action $a$ in state $s$, $R(s,a)$ is the scalar reward received for that transition, and $\gamma \in [0,1)$ is a discount factor that controls the trade-off between immediate and future rewards.

A \textit{policy} $\pi(a|s)$ defines a probability distribution over actions conditioned on the current state, and the central goal of an RL agent is to find a policy that maximizes the expected cumulative discounted reward. The quality of a policy is characterized by two closely related quantities. The state-value function $V^\pi(s) = \mathbb{E}_\pi\!\left[\sum_{t=0}^{\infty} \gamma^t\, r_t \,\middle|\, s_0 = s\right]$ gives the expected return when starting from state $s$ and thereafter following policy $\pi$, while the action-value function $Q^\pi(s,a) = \mathbb{E}_\pi\!\left[\sum_{t=0}^{\infty} \gamma^t\, r_t \,\middle|\, s_0 = s,\, a_0 = a\right]$ gives the expected return when first taking action $a$ in state $s$ and then following $\pi$. Both functions satisfy the Bellman equations~\citep{bellman1957dynamic}, which express a recursive relationship between the value of a state and the values of its successors:
\begin{equation}
\begin{aligned}
    Q^\pi(s,a) &= R(s,a) + \gamma \sum_{s'} P(s'|s,a)\, V^\pi(s'), \\
    V^\pi(s) &= \sum_{a} \pi(a|s)\, Q^\pi(s,a).
\end{aligned}
\label{eq:bellman}
\end{equation}

The empirical goal of these equations is to numerically express all the information about the future in the present. The algorithms developed to solve MDPs can be broadly grouped into three families. \textit{Value-based methods}, such as Q-learning~\citep{watkins1992qlearning}, learn the optimal action-value function $Q^*(s,a)$ directly and derive the policy by acting greedily with respect to it. Deep Q-Networks (DQN)~\citep{mnih2015human} extended this approach by approximating $Q^*$ with a convolutional neural network, demonstrating that RL agents could learn control policies from high-dimensional sensory inputs such as raw pixel observations, something that we take advantage of with $Dmsh$. However, value-based methods are largely restricted to discrete action spaces, which limits their applicability to problems requiring continuous control. 

\textit{Policy gradient methods} address this limitation by parameterizing the policy directly and optimizing it by estimating the gradient of the expected return with respect to the policy parameters. The REINFORCE algorithm~\citep{williams1992simple} is the prototypical example, using Monte Carlo roll-outs to construct an unbiased gradient estimate. Because the policy is represented explicitly, these methods naturally accommodate continuous and stochastic policies, though they can suffer from high variance in the gradient estimates.

\textit{Actor-critic methods}~\citep{konda1999actor} combine the strengths of both approaches: an \textit{actor} network learns the policy while a \textit{critic} network estimates the value function, which serves as a variance-reducing baseline for the policy gradient, thereby the pair learn to estimate the advantage of taking an action. This architecture forms the basis of many modern RL algorithms, including the Soft Actor-Critic method described in the following subsection, which augments the standard actor-critic objective with an entropy regularization term to encourage exploration whilst searching for optimality.

\subsection{Soft Actor-Critic (SAC) Methods}
\label{subsec:sac_methods}

Soft Actor-Critic (SAC) is an off-policy actor-critic deep reinforcement learning algorithm based on the maximum entropy reinforcement learning framework. Unlike standard reinforcement learning, which seeks only to maximize the expected sum of rewards, SAC aims to maximize the expected reward while also maximizing entropy. This entropy regularization encourages exploration and prevents the policy from converging prematurely to a sub-optimal deterministic behavior.

\subsubsection{The Core Objective}
The objective of Maximum Entropy Reinforcement Learning is to learn a policy $\pi$ that maximizes the weighted sum of expected reward and entropy:

\begin{equation}
    J(\pi) = \sum_{t=0}^{T} \mathbb{E}_{(s_t, a_t) \sim \rho_\pi} \left[ r(s_t, a_t) + \alpha \mathcal{H}(\pi(\cdot|s_t)) \right],
    \label{eq:sac_objective}
\end{equation}
where $\mathcal{H}(\pi(\cdot|s_t)) = -\mathbb{E}_{a \sim \pi}[\log \pi(a|s_t)]$ denotes the entropy of the policy at state $s_t$, and $\alpha$ is the \textit{temperature parameter} that controls the stochasticity of the optimal policy.

\subsubsection{Soft Value Functions}
To incorporate the entropy terms, SAC modifies the standard Bellman equations to derive two ``soft" value functions, viz., Soft Q-value and the soft state-value.
The soft Q-value, $Q(s_t, a_t)$, accounts for both the immediate reward and the expected future value of the state, which includes the entropy of future actions, given by:
\begin{equation}
    Q(s_t, a_t) = r(s_t, a_t) + \gamma \mathbb{E}_{s_{t+1} \sim p} \left[ V(s_{t+1}) \right].
    \label{eq:soft_q}
\end{equation}
The value function $V(s_t)$ is the additive sum of expected Q-value and the entropy of the policy at that state, given by:
\begin{equation}
    V(s_t) = \mathbb{E}_{a_t \sim \pi} \left[ Q(s_t, a_t) - \alpha \log \pi(a_t|s_t) \right].
    \label{eq:soft_v}
\end{equation}

\subsubsection{The Learning Process}
The algorithm optimizes the following three components through the stochastic gradient descent: the critic (Q-functions), the actor (policy), and the temperature ($\alpha$).
\paragraph{Soft Q-Function Update}
SAC utilizes \textit{Clipped Double-Q Learning} to mitigate positive bias. Two Q-networks are parameterized by $\theta_1, \theta_2$. The parameters are trained to minimize the soft Bellman residual:
\begin{equation}
    J_Q(\theta) = \mathbb{E}_{(s, a, r, s') \sim \mathcal{D}} \left[ \frac{1}{2} \left( Q_\theta(s, a) - \left(r + \gamma V_{\bar{\theta}}(s')\right) \right)^2 \right].
    \label{eq:critic_loss}
\end{equation}
\noindent where $V_{\bar{\theta}}(s')$ is the target value estimate derived from a target network. It is worth reiterating that only one of the 2 networks are trained via gradient descent and the other is trained as a moving average of the first network to maintain stability. 

\paragraph{Policy Optimization}
The policy parameters $\phi$ are optimized to maximize the expected Q-value plus the entropy term. Using the reparameterization trick $a_t = f_\phi(\epsilon_t; s_t)$ (where $\epsilon_t$ is sampled from a fixed Gaussian), the objective becomes differentiable:
\begin{equation}
    J_\pi(\phi) = \mathbb{E}_{s \sim \mathcal{D}, \epsilon \sim \mathcal{N}} \left[ \alpha \log \pi_\phi(f_\phi(\epsilon; s)|s) - Q_\theta(s, f_\phi(\epsilon; s)) \right].
    \label{eq:actor_loss}
\end{equation}

\paragraph{Temperature Auto-Tuning}
To avoid manual tuning of the exploration trade-off, $\alpha$ is learned by minimizing the following objective with respect to a target entropy $\bar{\mathcal{H}}$:

\begin{equation}
    J(\alpha) = \mathbb{E}_{a \sim \pi} \left[ -\alpha \log \pi(a|s) - \alpha \bar{\mathcal{H}} \right].
    \label{eq:alpha_loss}
\end{equation}

\section{Meshing}

\subsection{Environment Architecture}
In any RL framework, the environment plays a crucial role in supporting convergence, and a delicate balance must be struck to ensure the agent cannot exploit shortcuts while still being conducive to learning the optimal policy. To this end, a custom reinforcement learning environment is written from the ground up to automate the generation of quadrilateral meshes via an advancing front technique~\citep{blacker1991paving}. As discussed in Section~2, the meshing process is modeled as a Markov Decision Process in which a continuous domain $\Omega$ is discretized through the sequential removal of area from the boundary $\partial \Omega$. The environment is built upon the \texttt{Shapely} library for computational geometry and \texttt{NumPy} for vector manipulations. The formulation defines the state space as a partial, rotationally invariant representation of the local geometry. The action space is designed to support both the insertion of new vertices and the establishment of connectivity among existing nodes. In addition, the environment incorporates repair operators that enable recovery from geometrically degenerate configurations.

\subsubsection{Front Selection and Reference Vertex}
At each time step $t$, a specific segment of the boundary, $\partial \Omega_t$ is selected to serve as the local context for the agent. A process is Markov if, and only if, the current state contains all the information needed to predict the distribution of the next state, so truly random selection would break the Markov assumption. Hence, a deterministic priority-based selection mechanism is employed.

A priority score $S(v_i)$ is computed for every vertex $v_i \in \partial \Omega_t$. This score is a weighted linear combination of the interior angle $\theta_i$ and the deviation of the adjacent edge lengths from a target element size $L_{target}$:
\begin{equation}
    S(v_i) = \lambda \cdot \theta_i + (1 - \lambda) \cdot \frac{|\bar{l}_i - L_{target}|}{L_{target}},
\end{equation}
where $\bar{l}_i$ denotes the mean length of the two edges incident to $v_i$, and $\lambda$ is a weighting factor empirically set to the inverse golden ratio ($\approx 0.618$) to balance angular validity with dimensional uniformity. The vertex $v_0$ minimizing $S(v_i)$ is selected as the reference vertex. To prevent the agent from cycling through unsolvable configurations and entering an infinite loop recursively, a ``blacklist'' mechanism is utilized; vertices associated with repeated invalid actions are temporarily excluded from selection. Moreover, to ensure that the best quality elements are placed near the original boundary, nodes representing the original boundary are further multiplied by a weighing factor. 

\subsubsection{State Representation}\label{sssec:MeshStateRep}
The pursuit for creating a rotationally and translationally invariant representation rules out passing global coordinates. Consequently, a \textbf{local, rotationally invariant polar coordinate system} is established for each observation, centered at $v_0$, the reference vertex.

\paragraph{Local Coordinate Frame Definition}
The local frame is constructed by fixing the origin at the reference vertex $v_0$ and computing a bisector vector $\vec{b}$ from the unit vectors pointing toward the left neighbor $v_{L1}$ and the right neighbor $v_{R1}$. A candidate $x$-axis is aligned with $\vec{b}$; to ensure that this axis points into the interior of the domain $\Omega$, a test point $p_{test} = v_0 + \epsilon \vec{b}$ is generated, and if $p_{test} \notin \Omega$ (indicating a convex corner) the axis orientation is inverted. The $y$-axis is then derived via a $90^\circ$ counter-clockwise rotation of the corrected $x$-axis. This construction ensures that the agent empirically learns that an angle of zero corresponds to a point along the bisector, which speeds up convergence and keeps the action space symmetric with respect to the bisector.

\paragraph{Feature Vector Construction}
The state vector $s_t$, with dimension $D_{state} = 2(2N_{neigh} + N_{fan}) + 1$, is composed of features projected into this local frame and normalized by a scaling factor $L_{scale} = \beta \cdot L_{target}$ (where $\beta=6.0$). The first component consists of the polar coordinates $(r, \phi)$ of $N_{neigh}=3$ neighboring vertices to the left and right of $v_0$. To perceive non-adjacent boundary walls or obstacles deeper within the domain, $N_{fan}=3$ fan-sensing rays are cast within the angular wedge defined by $v_{L1}$ and $v_{R1}$; the intersection points of these rays with the global boundary are computed, transformed to the local frame, and appended to the state. Finally, a scalar domain ratio $\rho_t = A_t / A_{initial}$ is included to indicate the remaining area relative to the initial domain size. This domain ratio is the only sense of global context the agent receives, and proves to be sufficient for it to develop an effective global picture of the meshing progress.

\subsubsection{Action Space Dynamics}
A hybrid discrete-continuous action space is utilized, denoted as $(a_{type}, a_{params})$, with validity masks applied to filter action types based on the vertex count of the current boundary. The first action type (\textit{Type~0: Stitching}) is discrete: a quadrilateral is formed by connecting $v_0$ to an existing non-adjacent vertex, with candidate quads constructed from index sets $\{(i-1, i, i+1, i+2)\}$ and $\{(i-2, i-1, i, i+1)\}$. This stitching action is critical for closing the mesh loop and is solely responsible for reducing the global number of boundary nodes, since the alternative action adds and removes a node simultaneously. The second action type (\textit{Type~1: Creation}) is continuous: a new vertex $v_{new}$ is generated by mapping the continuous parameters $a_{params} \in [-1, 1]^2$ to local polar coordinates $(r, \theta)$, where the radius is scaled to $r \in [0, \alpha \cdot L_{target}]$ (with $\alpha=2.0$) and the angle is mapped relative to the current wedge aperture. The global coordinates of $v_{new}$ are then reconstructed via an inverse transformation of the local frame matrix.

It is worth noting that for the stitching action, the specific choice of which quadrilateral to close is performed by the environment rather than the agent, relieving the agent of some combinatorial complexity. Upon action proposal, strict geometric validity checks are enforced: the resulting quadrilateral must be convex, strictly contained within the prepared geometry of $\partial \Omega_t$, and free of intersections with any existing edges.

\subsubsection{Reward Formulation}
The reward function $R(s, a)$ is designed to drive the policy toward generating regular, isotropic elements while strictly penalizing topological defects and ensuring efficient termination. We define the total reward as a composite of geometric quality, boundary fidelity, and topological constraints:
\begin{equation}
    R(s, a) = (\eta_e + \eta_b - 1) + P_{speed} + P_{valence}.
\end{equation}
We introduce several modifications to the reward framework introduced by  Pan et al. \citep{pan2023reinforcement} to enhance convergence and mesh regularity. First, rather than relying on simple aspect ratios, the element quality term $\eta_e$ utilizes the Scaled Jacobian \citep{knupp2001algebraic} of the generated quadrilateral, normalized to $[0, 1]$, to rigorously detect and penalize inverted elements; we frame our entire optimization around this value. Second, we redefine the boundary quality term $\eta_b$ as a multiplicative aggregate of angular deviation, boundary narrowness, and transition smoothness, ensuring that the agent prioritizes the long-term "health" of the evolving boundary over greedy element placement. To enforce size constraints, a penalty $P_{speed}$ is applied when the element area deviates from the optimal range $[A_{min}, A_{max}]$, actively discouraging the formation of degenerate slivers as well as abnormally large elements. Finally, we introduce a discrete valence penalty $P_{valence}$ that is triggered whenever a vertex in the new element exceeds a valence of four; this explicitly promotes a regular grid topology and minimizes singularities.

\noindent To proactively prevent the advancing front from collapsing into bottlenecks or creating self-intersections, we introduce a specific narrowness penalty $\phi_{narrow}$ into the reward structure. This term evaluates the proximity of the active reference vertex $v_{ref}$ to non-adjacent segments of the boundary, guarding against the formation of singularities that are geometrically impossible to mesh. We first define a narrowness ratio $\rho$, calculated as the ratio of the minimum clearance distance $d_{min}$ to the local element scale:
\begin{equation}
    \rho = \frac{\min_{v \in \mathcal{B} \setminus \mathcal{N}(v_{ref})} \| v - v_{ref} \|}{2 \cdot \bar{l}_{local}},
\end{equation}
where $\mathcal{B}$ denotes the set of boundary vertices, $\mathcal{N}(v_{ref})$ represents the immediate topological neighbors of the reference vertex (excluded to avoid trivial zero distances), and $\bar{l}_{local}$ is the average length of the edges connected to $v_{ref}$. The penalty is applied via a piecewise function that imposes a soft barrier as the front narrows ($\rho < 1.0$) and a severe penalty when collision is imminent ($\rho < 0.5$):
\begin{equation}
    \phi_{narrow} = 
    \begin{cases} 
        1.0 & \text{if } \rho \ge 1.0 \\
        \rho & \text{if } 0.5 \le \rho < 1.0 \\
        0.3\rho & \text{if } \rho < 0.5.
    \end{cases}
\end{equation}
Furthermore, to ensure the updated boundary maintains both angular flatness and metric uniformity, critical factors for the stability of subsequent element generation. We implement a smoothness penalty $\phi_{smooth}$. This factor is computed as a weighted aggregate of an angular consistency term $Q_{ang}$ and a metric consistency term $Q_{met}$, balancing geometric regularity with size constraints:
\begin{equation}
    \phi_{smooth} = 0.6 \cdot Q_{ang} + 0.4 \cdot Q_{met}.
\end{equation}

\noindent The angular term $Q_{ang}$ penalizes sharp kinks in the boundary by measuring the maximum deviation of local interior angles $\theta_i$ from a straight line ($\pi$ radians) within the immediate neighborhood of the new element. Simultaneously, the metric term $Q_{met}$ evaluates the uniformity of edge lengths, penalizing high variance ($\sigma_l / \mu_l$) and deviation from the target characteristic length $l_{target}$:
\begin{equation}
    Q_{ang} = 1 - \min\left(1, \frac{\max_{i} | \theta_i - \pi |}{\pi/2} \right), \quad
    Q_{met} = \exp\left( - \left( \frac{\sigma_l}{\mu_l} + 0.5 \cdot \frac{|\mu_l - l_{target}|}{l_{target}} \right) \right).
\end{equation}
Collectively, these modifications empower the agent to transcend simple element quality maximization, strategizing node placement to maintain a boundary that is free of narrow bottlenecks or sharp irregularities detrimental to convergence.

\subsubsection{Failure Recovery and Mesh Repair}
Given the fragile nature of the computational geometry, a robust recovery loop is implemented to resolve states where no valid action can be sampled. If the consecutive-failure counter surpasses a threshold, a multi-stage repair sequence is initiated. The first stage performs a heuristic geometric expansion in which vertices exhibiting reflex angles ($>180^\circ$) or acute kinks are iteratively shifted toward ideal positions defined by their neighbors (targeting $90^\circ$ or $135^\circ$ interior angles). If this heuristic pass fails to restore validity, the second stage applies a constrained Laplacian smoothing in which all non-fixed nodes are updated according to:
\begin{equation}
    v_i^{k+1} = v_i^k + \gamma ( \bar{v}_{neighbors} - v_i^k ).
\end{equation}
\noindent Strict domain containment is enforced, where nodes are permitted to slide along the boundary or traverse the interior but are strictly prevented from crossing domain boundaries. Following smoothing, a global integrity check is executed. If self-intersections persist, the environment state is reverted to the pre-recovery snapshot and the episode is terminated with a heavy penalty.

\begin{algorithm}
\caption{Failure Recovery Loop for Geometric Deadlocks}
\label{alg:recovery}
\begin{algorithmic}[1]
\State \textbf{Input:} Current boundary $\partial \Omega_t$, target element size $L_{target}$, threshold $T_{fail}$
\State \textbf{Initialize:} $c_{fail} \gets 0$
\If{$action\_valid == \text{False}$}
    \State $c_{fail} \gets c_{fail} + 1$
\EndIf

\If{$c_{fail} \ge T_{fail}$}
    \State \textbf{Stage 1: Geometric Expansion (Heuristic)}
    \For{each vertex $v_i \in \partial \Omega_t$}
        \If{$\theta_i > 180^\circ$ (Reflex Angle)}
            \State Calculate ideal position $v_{ideal}$ based on $90^\circ/135^\circ$ targets
            \State $v_i \gets v_i + \lambda(v_{ideal} - v_i)$ \Comment{Shift toward local convexity}
        \EndIf
    \EndFor
    
    \State \textbf{Stage 2: Constrained Laplacian Smoothing}
    \While{not converged \textbf{and} $k < k_{max}$}
        \State $v_{new} = v_i^k + \gamma \left( \frac{1}{|\mathcal{N}_i|} \sum_{j \in \mathcal{N}_i} v_j^k - v_i^k \right)$
        \If{$\text{Prepared}(\partial \Omega_{init}).\text{contains}(v_{new})$}
            \State $v_i^{k+1} \gets v_{new}$
        \Else
            \State Project $v_{new}$ to nearest valid boundary segment
        \EndIf
        \State $k \gets k + 1$
    \EndWhile
    
    \State \textbf{Stage 3: Global Integrity Check}
    \If{$\partial \Omega_{t}$ is self-intersecting \textbf{or} $\text{Area}(\partial \Omega_t) < \epsilon$}
        \State Revert to snapshot $\partial \Omega_{t-1}$
        \State \textbf{Terminate:} Return Failure Penalty
    \Else
        \State $c_{fail} \gets 0$
        \State \textbf{Return:} Updated $\partial \Omega_t$
    \EndIf
\EndIf
\end{algorithmic}
\end{algorithm}

The smoothing contributes the most to helping repair the intermediate boundary to facilitate further meshing, while the geometric expansion adds some amount of variance to the smoothed nodes to enable distinct vertex selection and higher quality elements. 

The strict checks, combined with the occasional use of the recovery loops, mean that the environment takes up a major chunk of computational overhead ($\approx 72\%$) during meshing, with the actual inferencing of the agent taking a very small chunk of the time ($\approx 18\%$).

\subsubsection{Termination Criteria}
An episode terminates under one of three conditions. The first, and desired, outcome is successful discretization of the domain, which occurs when the number of boundary vertices reaches $N \le 4$; in this case, a final high-fidelity Laplacian smoothing is applied, and a finishing reward $R_{finish}=10$ is granted. The episode also terminates if the recovery mechanism described above fails to resolve a geometric deadlock, or if a maximum step count or maximum consecutive-failure limit is reached. Consequently, relevant penalties/rewards are given depending on the criteria/degree of offense. 

\subsection{Parametric Soft Actor-Critic (P-SAC) Agent}

We now turn our attention to $Dmsh$'s first agent, the meshing agent. Empirically, it is important that the agent architecture aligns with the structure of the action space to support convergence to the optimal policy. While the current SOTA implementation employs a standard Soft Actor-Critic method, the hybrid discrete-continuous nature of choosing an action, then choosing the parameters of that action in meshing, motivates the adoption of a Parametric Soft Actor-Critic (P-SAC) framework inspired by Fan et al., Lin et al.~\citep{fan2019hybrid,lin2024discretionarylanechangedecisioncontrol,Discretecontactionrep}. In this formulation, the action space $\mathcal{A}$ is defined as a tuple $(k, z)$, where $k \in \{0, 1\}$ represents a discrete action type and $z \in \mathbb{R}^n$ represents the continuous parameters associated with that type. 

Unlike standard SAC, which assumes a homogeneous action space (either purely discrete or purely continuous), this implementation utilizes a unified actor with split heads and a modular critic architecture that decouples the value estimation of distinct action topologies. In practice, this allows the actor to maintain two separate output spaces: one to decide which action type to execute, and another to determine the continuous parameters if type~1 is chosen.

We maintain 2 classes of critics, one for type~0 discrete action and one for type~1 continuous action. Within each class of critic, there are 4 neural networks, 2 of which are online and are updated through gradient descent (twin networks to prevent overestimation bias), and 2 networks that are updated using polyak averaging to maintain stability. The averaged networks act as target networks for the Bellman update.

\subsubsection{Unified Parametric Actor}

\begin{figure}[htpb]
     \centering
     \includegraphics[width=0.50\textwidth]{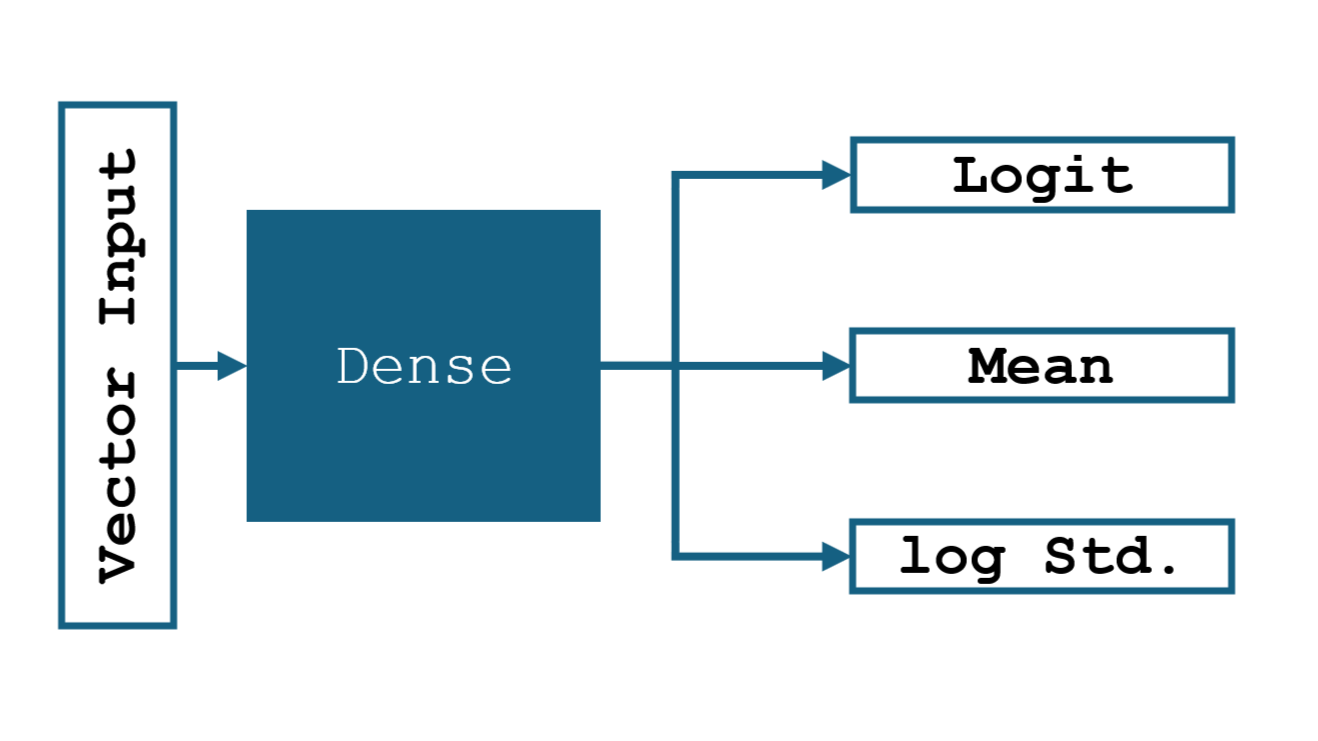}
     \caption{Schematic of the meshing agent architecture}
\end{figure}

The policy network, parameterized by $\theta$, utilizes a shared feature extractor followed by two specialized heads. The discrete head outputs logits $l \in \mathbb{R}^2$ to model the categorical distribution $\pi_{disc}(k|s)$ over action types, while the continuous head outputs the mean $\mu(s)$ and log standard deviation $\log \sigma(s)$ to model the conditional Gaussian distribution $\pi_{cont}(z|s, k=1)$. For action type $k=0$, no continuous parameters are generated since the actual execution is handled by the environment, and is deterministic in nature, given the current state. The full policy density is given by the chain rule:
\begin{equation}
    \pi(a|s) = \pi_{disc}(k|s) \cdot \mathbb{I}_{[k=1]} \pi_{cont}(z|s) + \pi_{disc}(k|s) \cdot \mathbb{I}_{[k=0]}.
\end{equation}

\subsubsection{Asymmetric Target Entropy Tuning}

To appropriately balance exploration and exploitation in this hybrid space, the P-SAC formulation utilizes a dual temperature tuning mechanism. While the continuous action space utilizes a standard target entropy heuristic, the discrete action space requires a specialized approach due to the asymmetric nature of the optimal policy. 

Standard discrete SAC typically sets a target entropy that encourages a uniform distribution across actions. However, our environment exhibits a highly skewed optimal action distribution. Since we place more nodes inside the boundary than connect existing nodes (type~0 action is called as many times as there are boundary nodes), the continuous action is called almost $4\times$ as often as the discrete action. A uniform target entropy would artificially penalize the agent for converging to this optimal, low-entropy behavioral prior. 

To prevent this, we calculate the discrete target entropy $\mathcal{H}_{target}^{disc}$ using the Shannon entropy of our expected asymmetric binomial distribution: 
\begin{equation}
    \mathcal{H}_{target}^{disc} = - (0.18 \log 0.18 + 0.82 \log 0.82).
\end{equation}
By explicitly defining the target entropy to match this 18/82 split, the temperature parameter $\alpha_{disc}$ correctly scales down when the agent discovers the optimal skewed policy, preventing the algorithm from forcing sub-optimal exploration. The 18-82 split is arrived at after empirical experimentation. 

\subsubsection{Decoupled Critic System}
To handle the variable dimensionality of the input space, the Q-function is decomposed into separate neural networks: $Q_{\phi_0}(s)$ estimates the value of the non-parametric stitching action, and $Q_{\phi_1}(s, z)$ estimates the value of the parametric creation action. Twin Delayed Q-learning is employed for each branch, resulting in four critic networks total ($Q_{\phi_{0,1}}, Q_{\phi_{0,2}}, Q_{\phi_{1,1}}, Q_{\phi_{1,2}}$) to mitigate overestimation bias. This decoupled structure is a natural fit for the meshing problem, as each value function can be individually evaluated based on the action type that was taken, thereby de-correlating their gradient updates and improving learning stability.  

\subsubsection{Critic Update (Masked Loss)}
Standard Bellman updates are modified to strictly update the critic corresponding to the sampled action type. Given a batch of transitions, the loss for the critics is computed using a binary mask $M_k = \mathbb{I}_{[action\_type = k]}$:
\begin{equation}
    J(\phi_k) = \mathbb{E}_{\tau \sim \mathcal{D}} \left[ M_k \cdot \left( Q_{\phi_k}(s, z_k) - y_{target} \right)^2 \right],
\end{equation}
where the target $y_{target}$ aggregates values from both action types weighted by their predicted probabilities:
\begin{equation}
    y_{target} = r + \gamma (1-d) \sum_{k \in \{0,1\}} \pi(k|s') \left( \min_{j=1,2} Q_{\text{target}, j}(s', \tilde{z}') - \alpha \log \pi(\cdot) \right).
\end{equation}

\subsubsection{Dual-Entropy Tuning}
The implementation addresses the scale difference between discrete logits and continuous densities by utilizing two distinct temperature parameters, $\alpha_{disc}$ and $\alpha_{cont}$, which are optimized independently:
\begin{equation}
\begin{aligned}
    J(\alpha_{disc}) &= \mathbb{E} [ -\alpha_{disc} (\mathcal{H}(\pi_{disc}) - \bar{\mathcal{H}}_{disc}) ] \\
    J(\alpha_{cont}) &= \mathbb{E} [ -\alpha_{cont} (\log \pi_{cont}(z|s) + \bar{\mathcal{H}}_{cont}) ].
\end{aligned}
\end{equation}
This ensures that exploration is maintained effectively in both the high-level strategy (action type selection) and low-level control (parameter selection).

Overall, this agent setup allows us to scale to more actions in a much more straightforward way that enables a future expansion to generate more polygonal elements, such as hexagonal elements, and even hexahedral elements for 3D meshing. We simply would have to increase the number of logits and add a critic for each type of new action created. For the purposes of this paper, we stick to all-quad meshing.

\subsection{Curriculum Learning}

A major issue with RL algorithms is that convergence is usually highly seed-dependent, as most optimal policies contain too many degrees of freedom to optimize in parallel, and hence a good starting point matters much more. In the case of quad meshing, the agent has to ensure completion, create high-quality elements, maintain a smooth, meshable boundary, as well as reduce singularities, all while staying relatively lightweight at around 100k parameters. To address this, a two-phase training procedure is introduced that utilizes curriculum learning~\citep{bengio2009curriculum} to progressively increase the problem complexity and achieve reliable convergence.

\subsubsection{Phase 1: Pre-training}
Before elements can be adjusted, the agent needs to understand what a good element contributes to. Two out of the four terms in the loss function are dedicated to creating good-quality elements; hence, it is important that the agent reaches this checkpoint before it proceeds with a more contextual response to the state. 

This is achieved by simply limiting the outer boundary to a perfect square, as it is the simplest geometry that perfectly splits into smaller squares, thereby creating a form of local minima for the agent to reach quickly, regardless of initialization. 

\begin{figure}[htpb]
     \centering
     \begin{subfigure}[b]{0.3\textwidth}
         \centering
         \includegraphics[width=\textwidth]{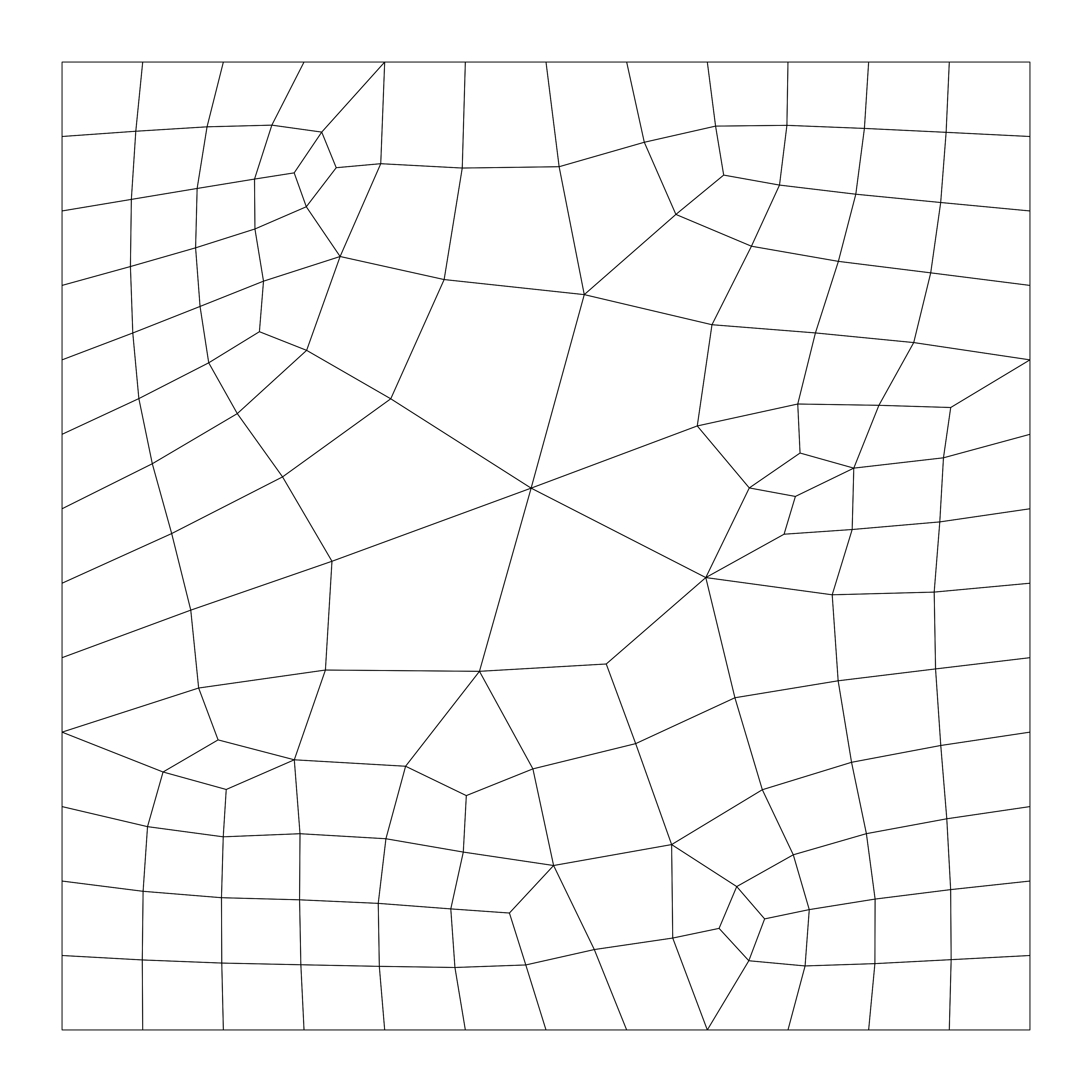}
         \caption{5k Steps}
     \end{subfigure}
     \hfill
     \begin{subfigure}[b]{0.3\textwidth}
         \centering
         \includegraphics[width=\textwidth]{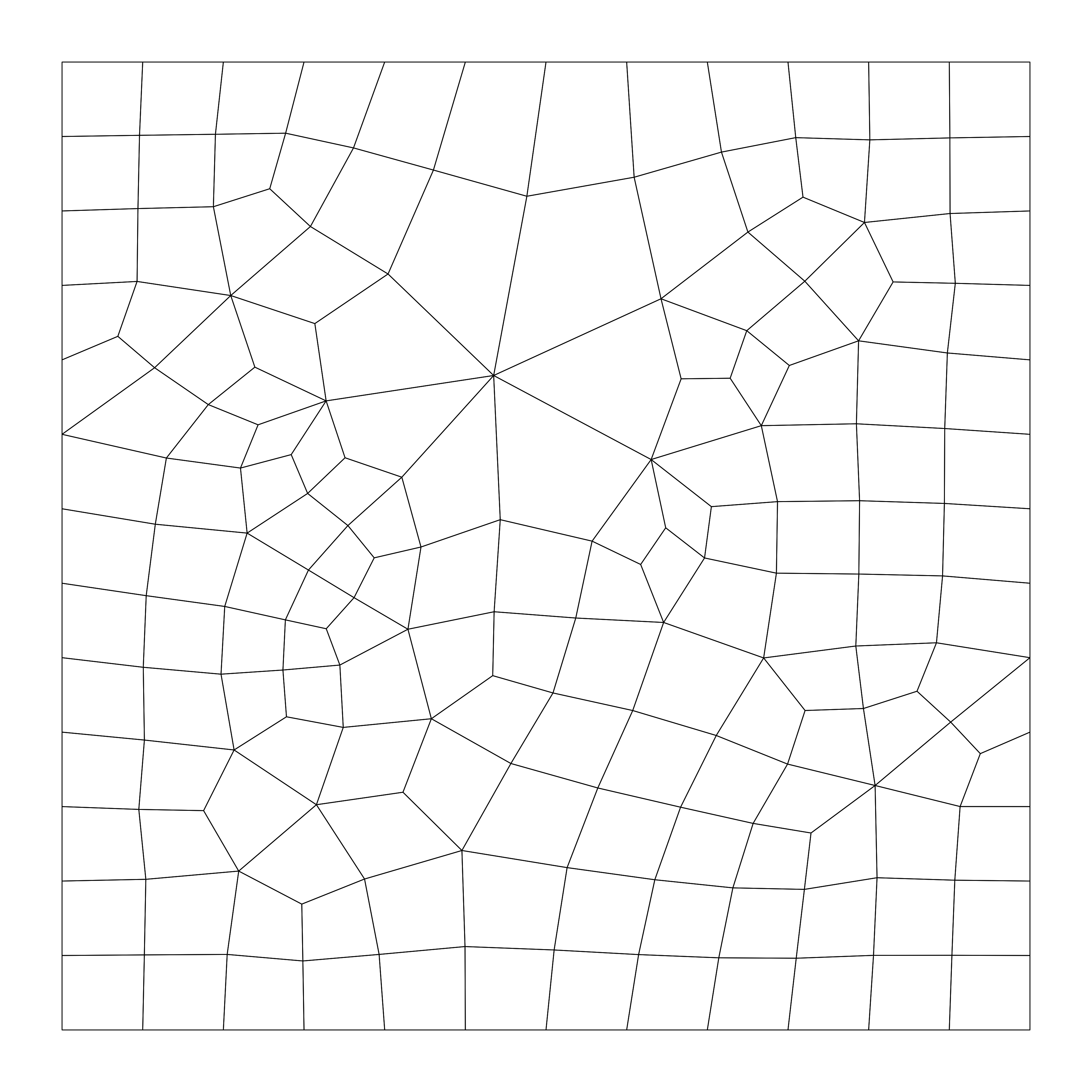}
         \caption{7k Steps}
     \end{subfigure}
     \hfill
     \begin{subfigure}[b]{0.3\textwidth}
         \centering
         \includegraphics[width=\textwidth]{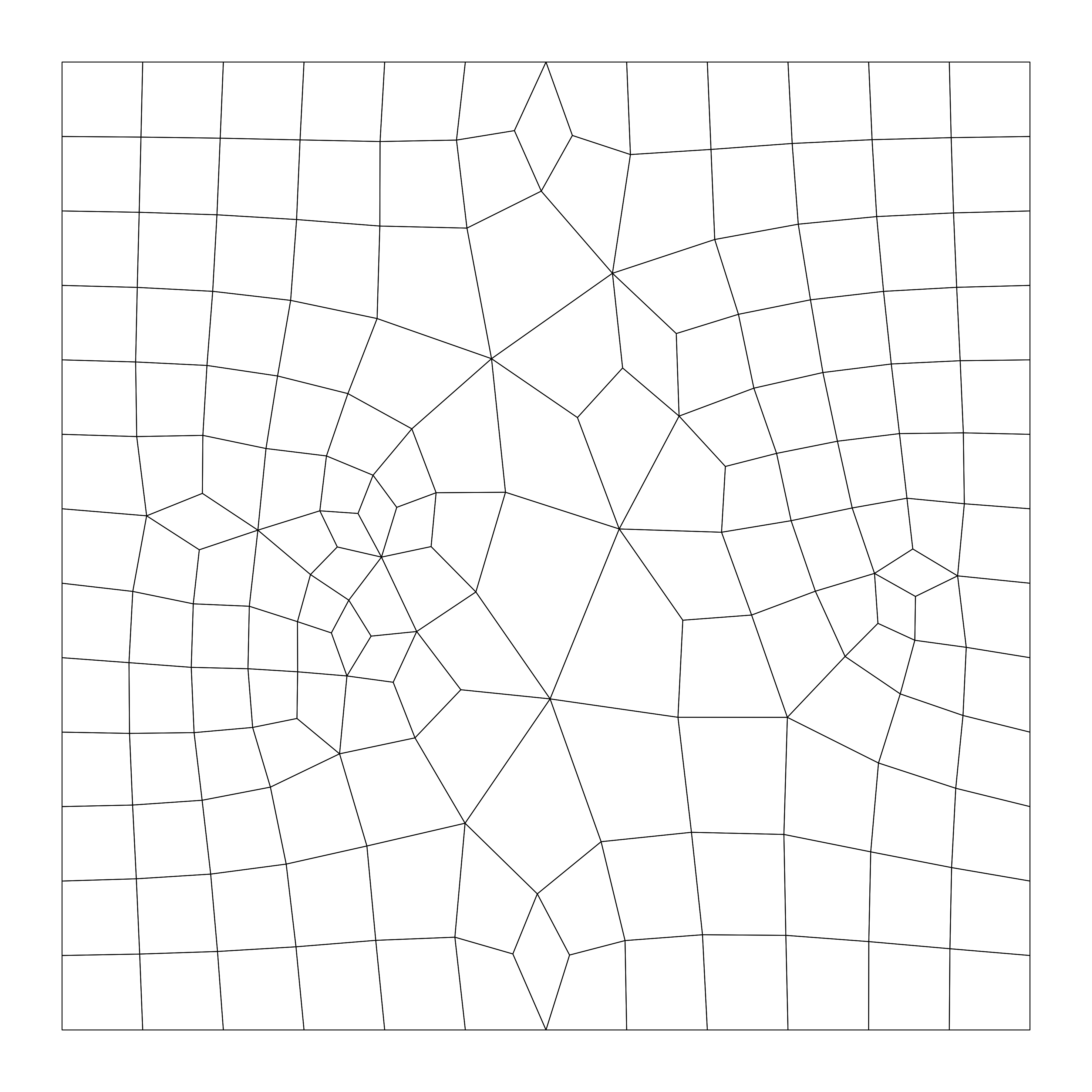}
         \caption{10k Steps}
     \end{subfigure}

     \vspace{1em} 

     \begin{subfigure}[b]{0.3\textwidth}
         \centering
         \includegraphics[width=\textwidth]{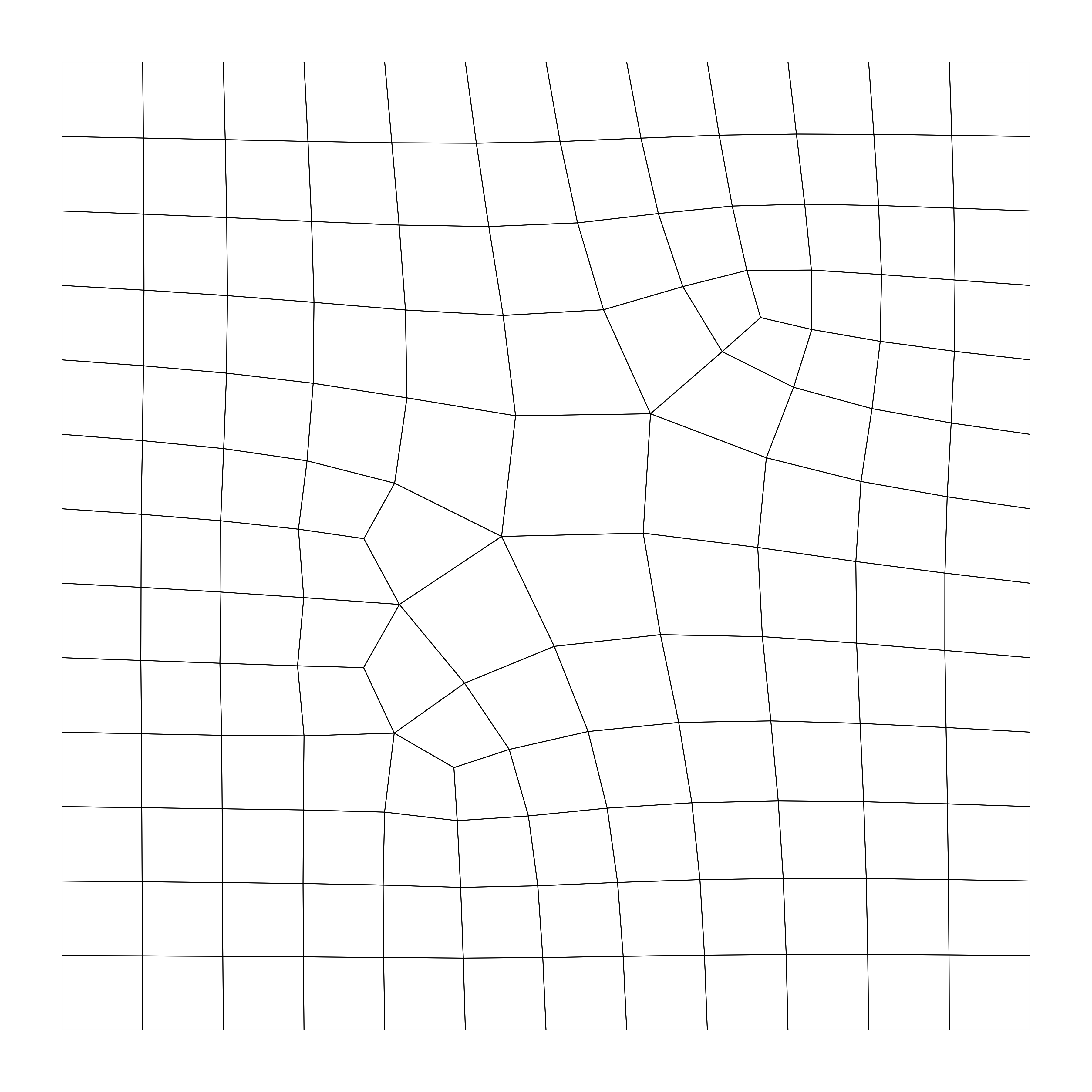}
         \caption{20k Steps}
     \end{subfigure}
     \hfill
     \begin{subfigure}[b]{0.3\textwidth}
         \centering
         \includegraphics[width=\textwidth]{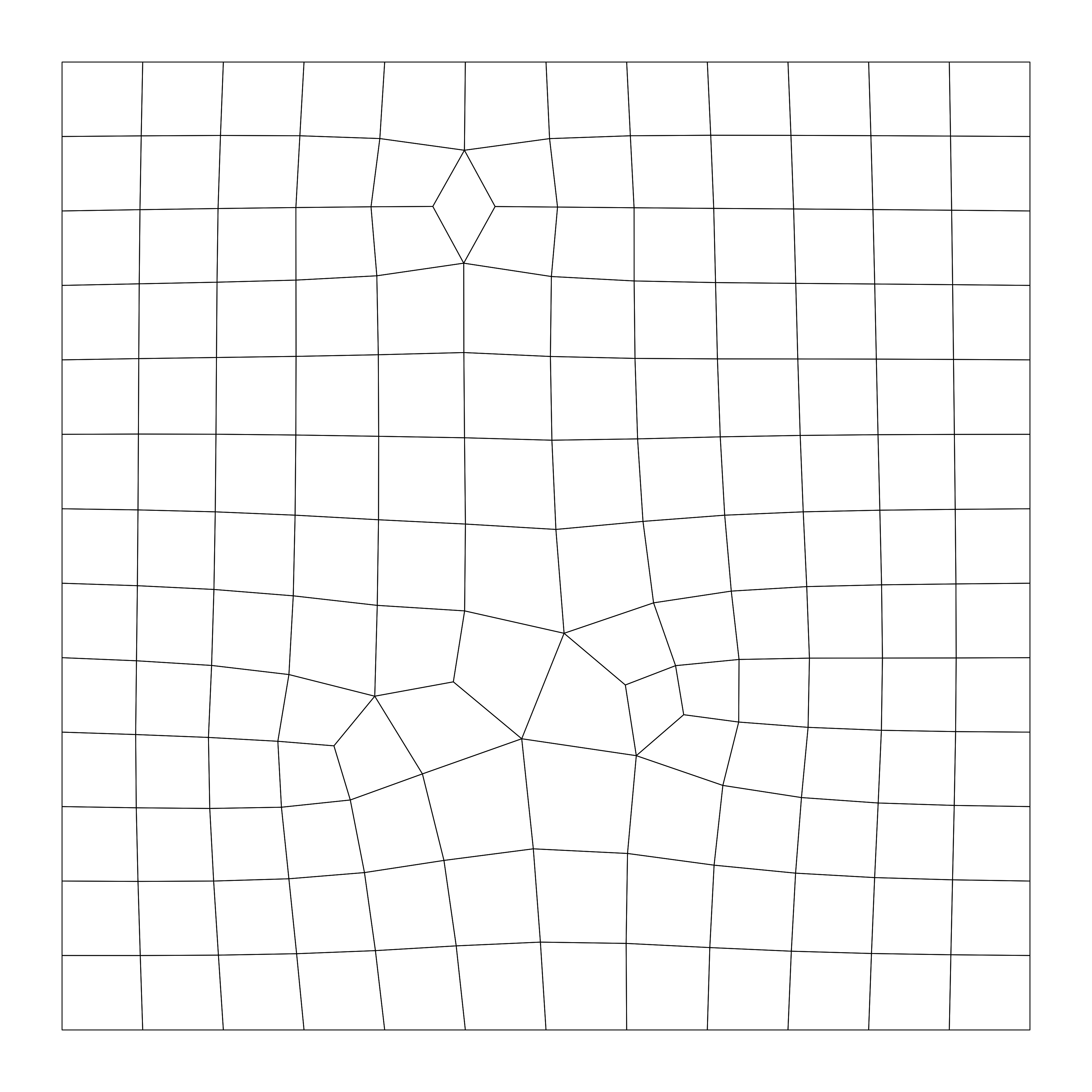}
         \caption{50k Steps}
     \end{subfigure}
     \hfill
     \begin{subfigure}[b]{0.3\textwidth}
         \centering
         \includegraphics[width=\textwidth]{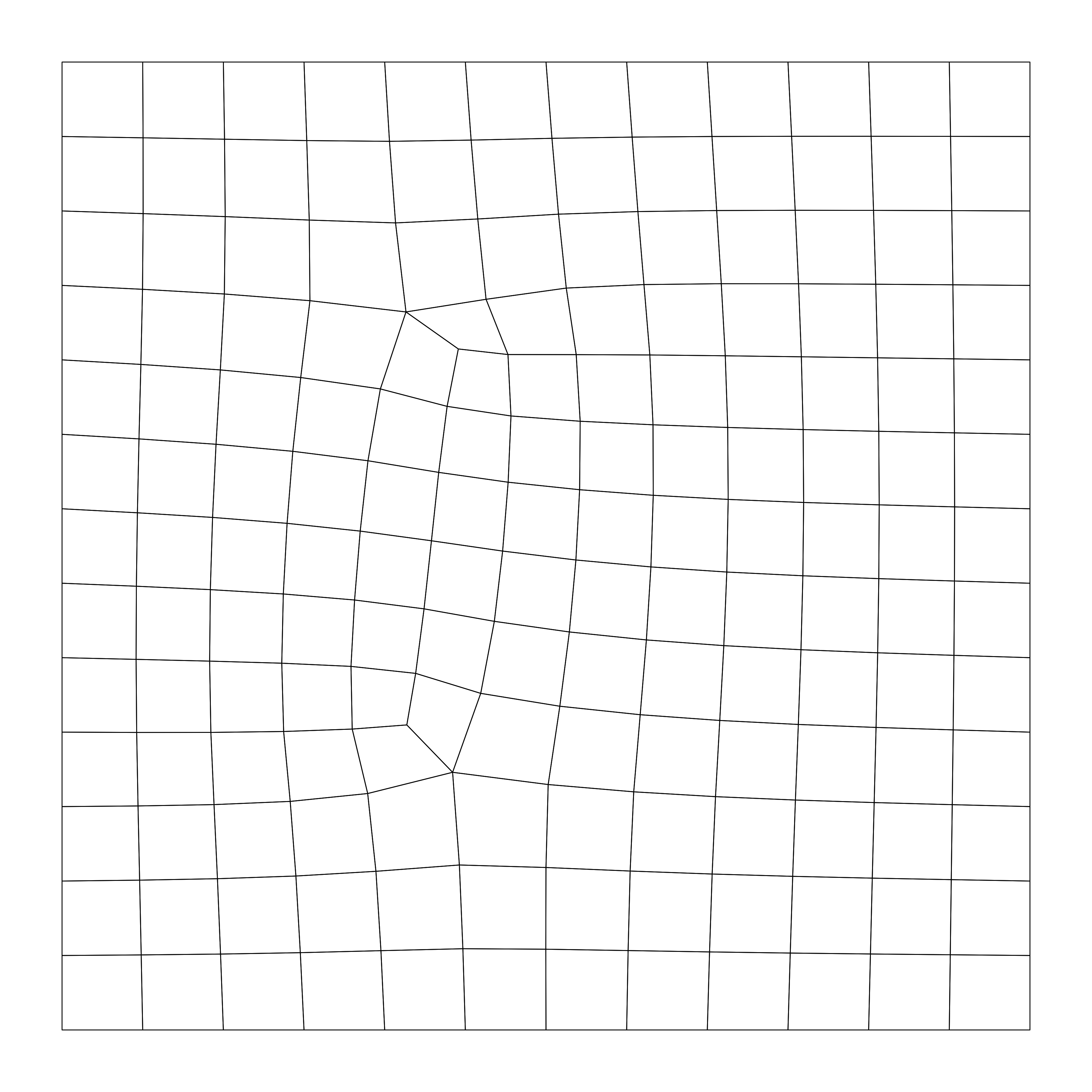}
         \caption{100k Steps}
     \end{subfigure}
     
     \caption{Evolution of the meshing agent during the curriculum pre-training phase on a square domain. Early in training (5k–10k steps), the agent produces irregular and poorly aligned quadrilateral elements. As training progresses (20k–100k steps), the agent learns stable advancing-front behavior, generating smoother boundaries, improved element orthogonality, and increasingly uniform all-quadrilateral meshes.}
     \label{fig:agent_progression}
\end{figure}

Empirically, this is set to be around 100,000 steps before the second phase of training is initiated. Figure~\ref{fig:agent_progression} illustrates the progression of the agent during this phase.

\subsubsection{Phase 2: Ensemble Training}

\begin{figure}[htpb]
     \centering
     \begin{subfigure}[b]{0.3\textwidth}
         \centering
         \includegraphics[width=\textwidth]{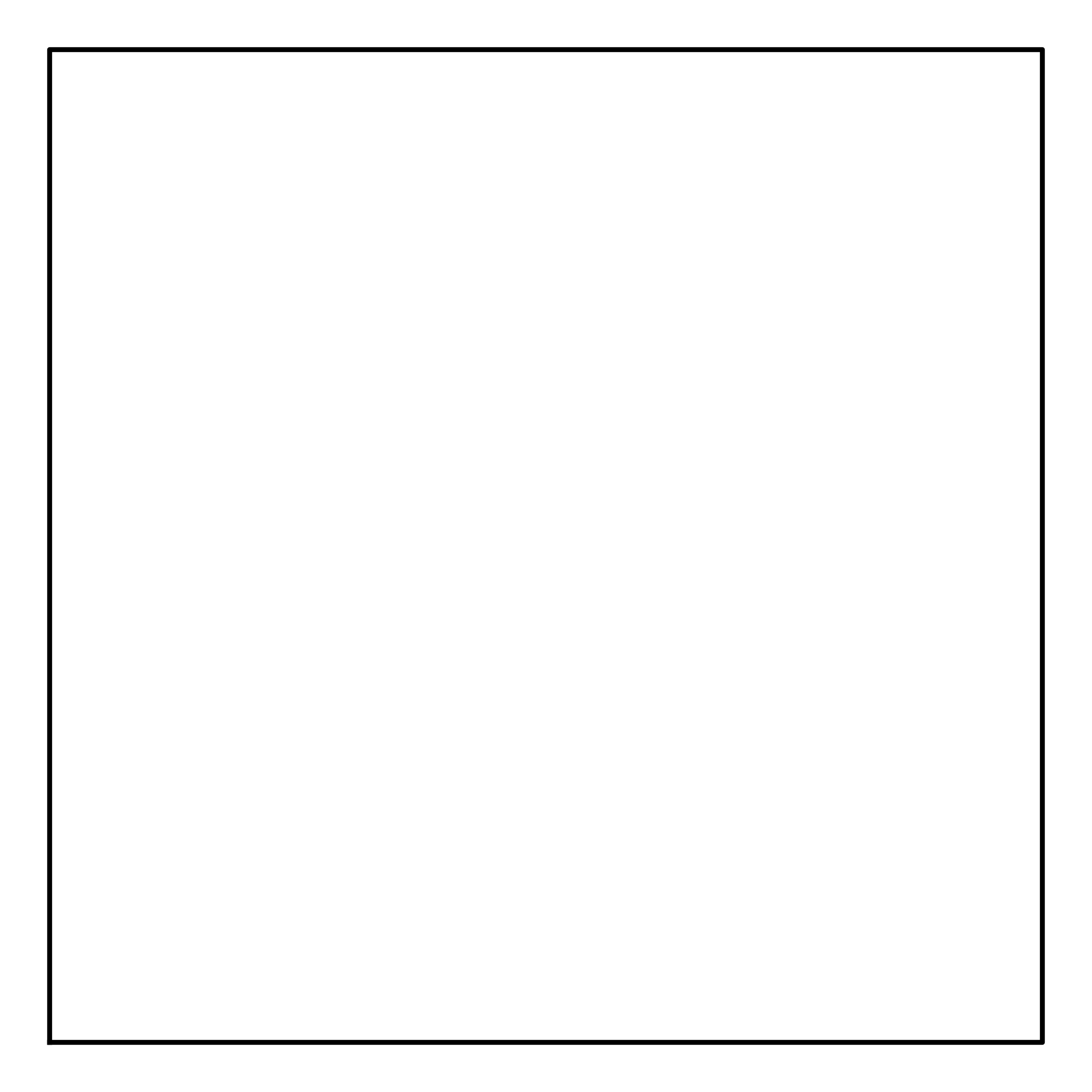}
     \end{subfigure}
     \hfill
     \begin{subfigure}[b]{0.3\textwidth}
         \centering
         \includegraphics[width=\textwidth]{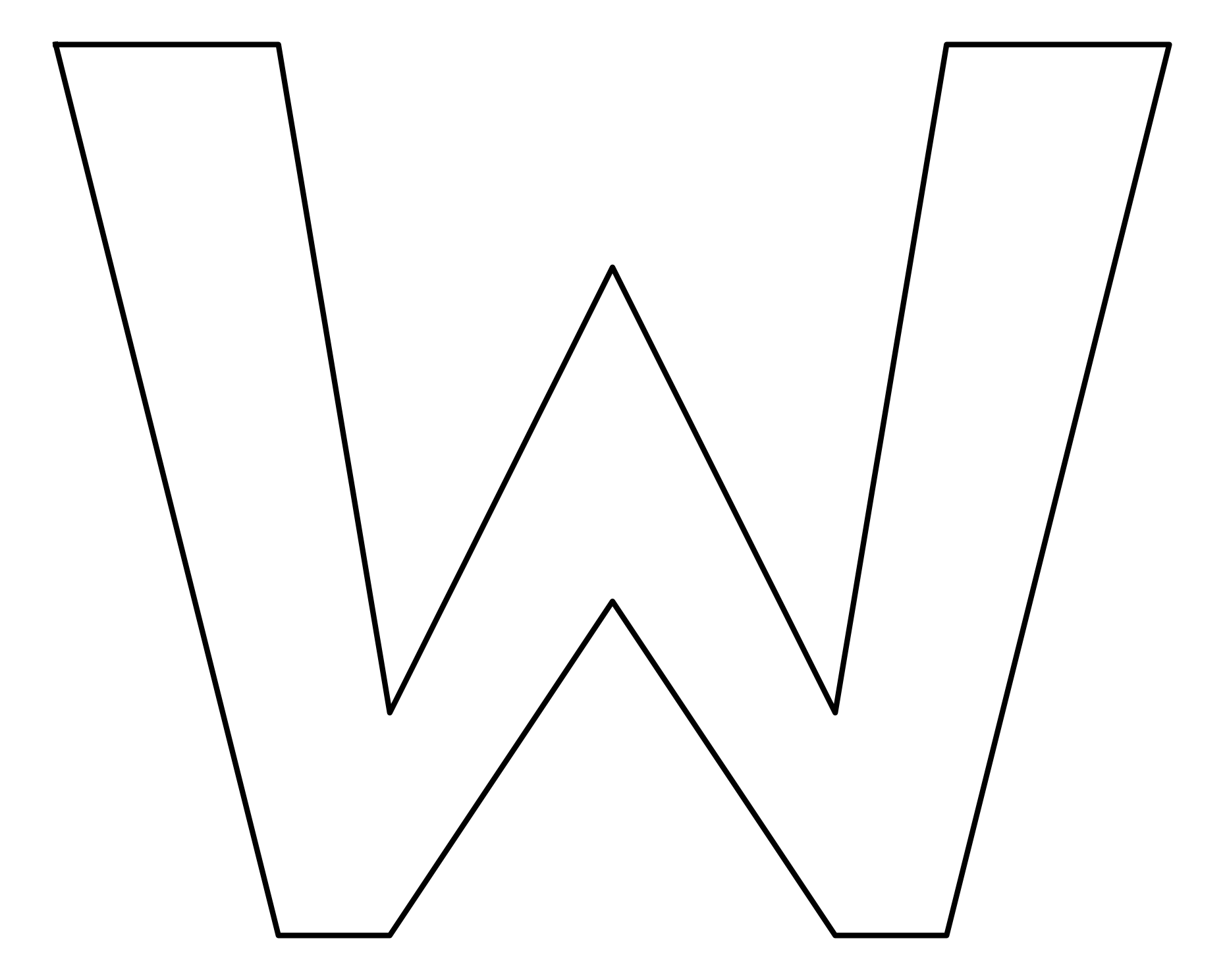}
     \end{subfigure}
     \hfill
     \begin{subfigure}[b]{0.3\textwidth}
         \centering
         \includegraphics[width=\textwidth]{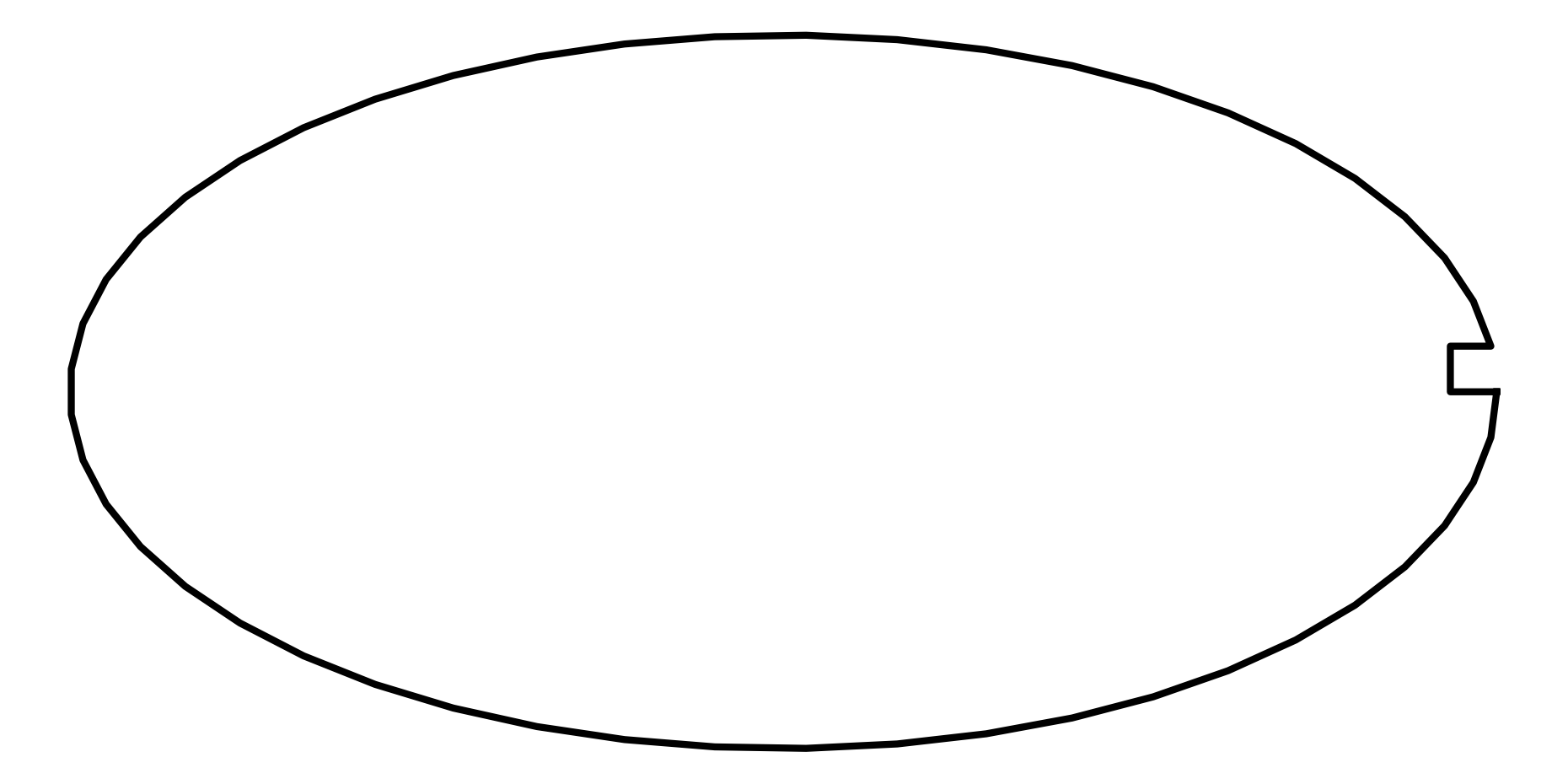}
     \end{subfigure}

     \vspace{1em} 

     \begin{subfigure}[b]{0.3\textwidth}
         \centering
         \includegraphics[width=\textwidth]{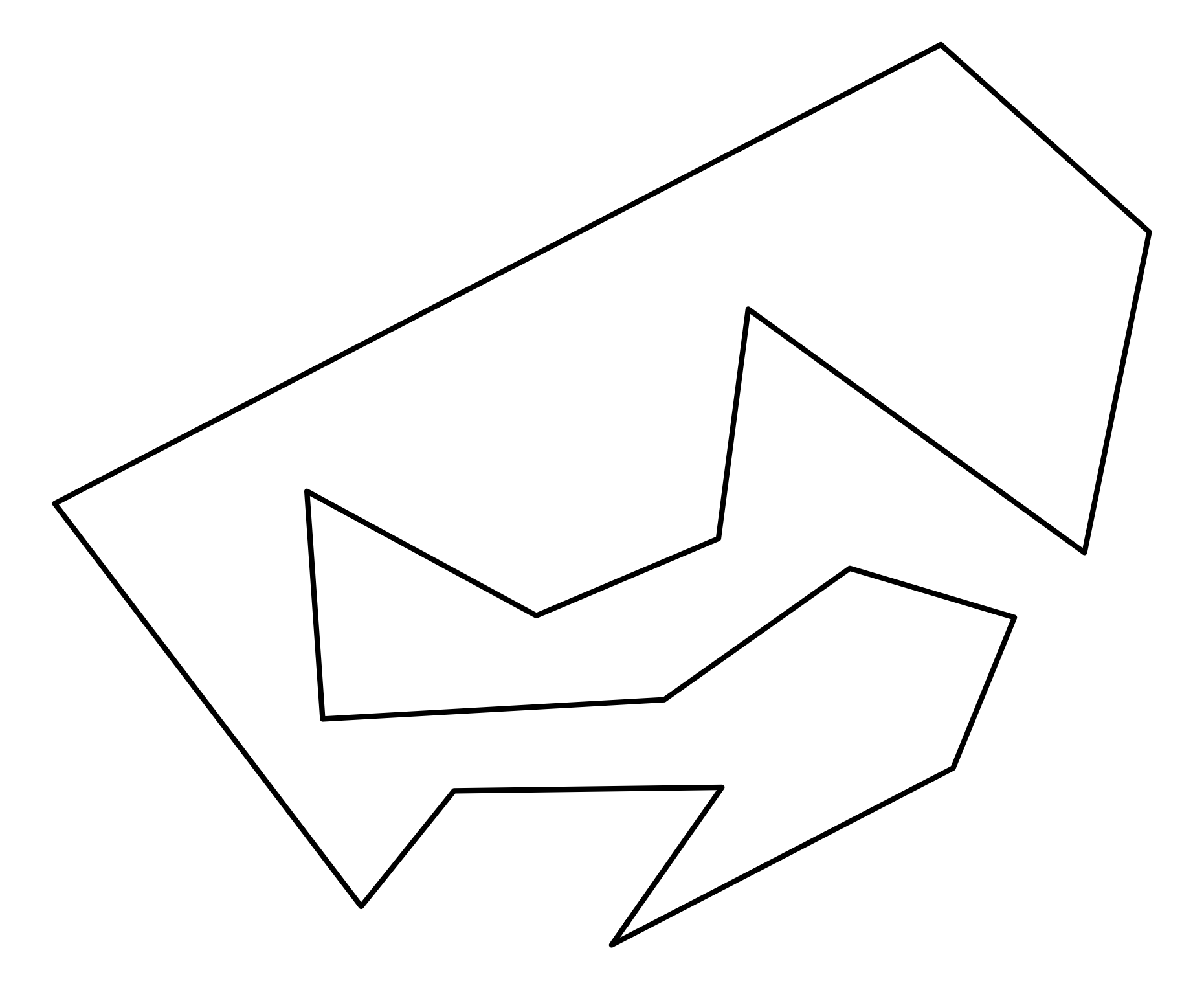}
     \end{subfigure}
     \hfill
     \begin{subfigure}[b]{0.3\textwidth}
         \centering
         \includegraphics[width=\textwidth]{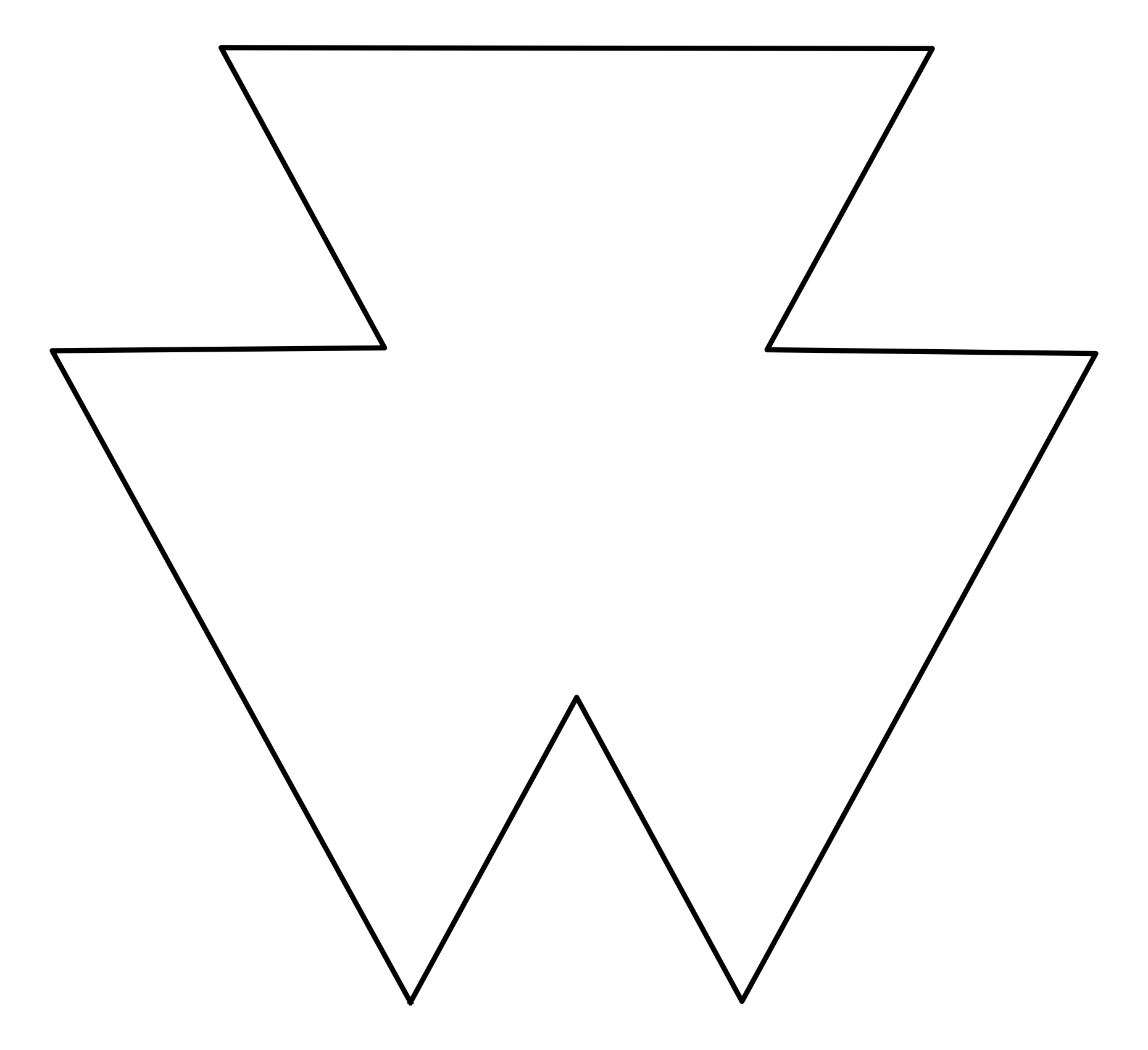}
     \end{subfigure}
     \hfill
     \begin{subfigure}[b]{0.3\textwidth}
         \centering
         \includegraphics[width=\textwidth]{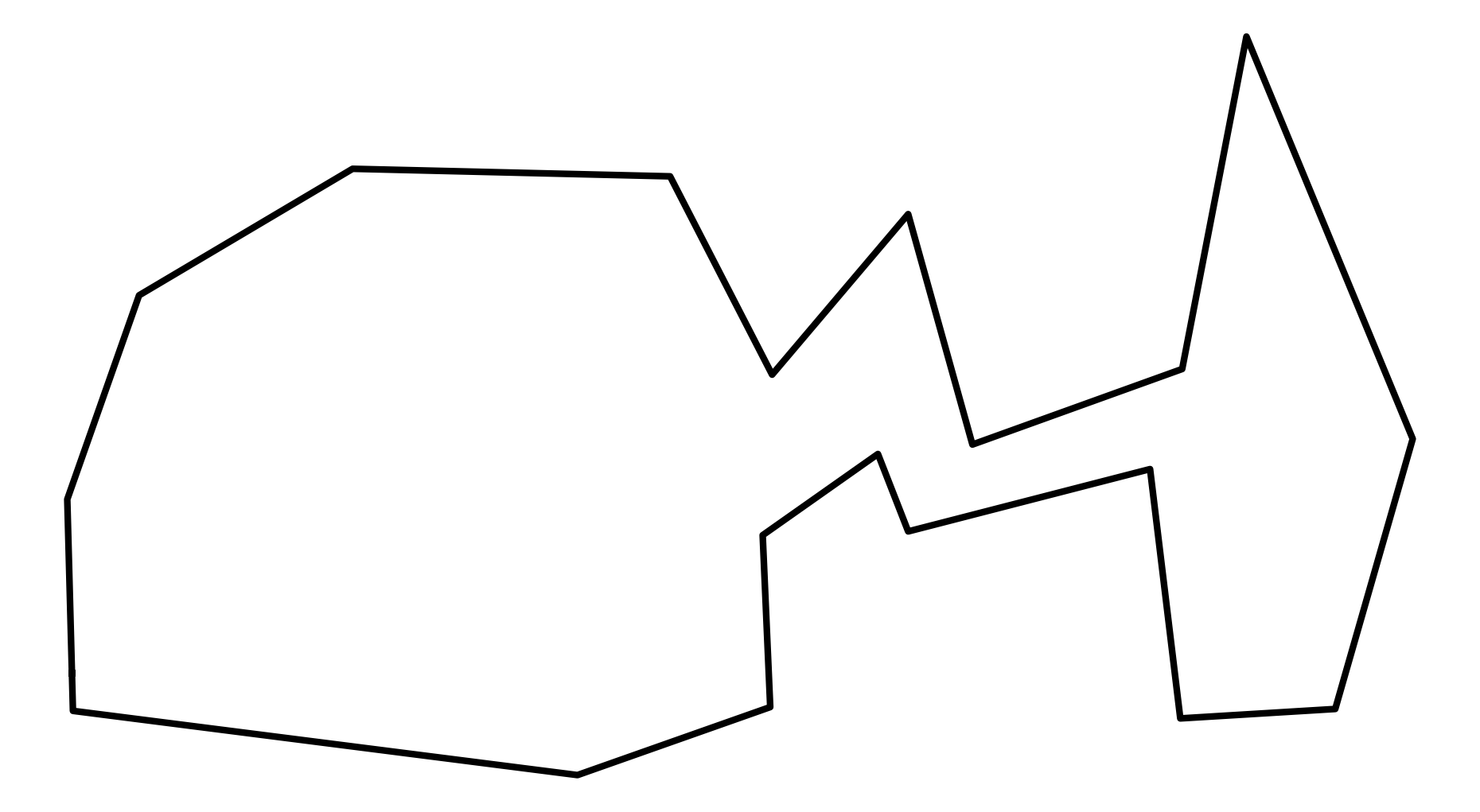}
     \end{subfigure}
     
\caption{Representative set of geometries used during the ensemble-training phase. The agent is exposed to domains with varying curvature, concavity, aspect ratios, and local topological complexity, enabling the policy to generalize beyond the simple square geometry used during curriculum pre-training.}
\label{fig:represenatative_geometries}
\end{figure}

We now introduce the most significant of the two phases. In this phase, the agent learns temporal relations between the geometry and the planning of elements to accommodate the ultimate geometry. An interesting advantage of this recursive formulation is that although the input geometries are only in the single digits, since each step creates a unique boundary and in turn a unique local state representation the agent sees a wide variety of different states, leading to a dense implicit dataset generated automatically. 

This effect is further exemplified using dataset \textbf{bagging}, where a set of geometries is associated with a score (net averaged reward), and a domain is chosen every $n$ training episodes on a softmax basis for the agent to learn from. The boundaries are updated every $n$ episodes.  

This prevents any form of overfitting and generates a dense variety of local data. This phase is significantly longer than the pre-training phase, lasting 1.5--2 million time-steps. The converged agent has the ability to work around most geometries and guarantees an all-quad mesh, unless an odd number of points are given, in which case the last element would remain a triangle. We report the results of the benefits of using curriculum learning in the results section.

\section{Block Decomposition}

Although the new RL-based meshing agent exhibits more robustness to oblique geometries, its fundamental paving method limits it from wrapping around highly convex points, causing either catastrophic breakdowns or low-quality elements near the sharp ends. Moreover, due to the increasing complexity of an incomplete boundary, reduced freedom, and strict convexity checks by the environment, the error rate increases with the number of elements made by a single agent, increasing meshing time in a non-linear fashion. Hence, another agent is developed that sits upstream of the meshing agent that splits the existing boundary into new blocks in a recursive manner. This decomposition also enables parallelisation, as now, each blocks are meshed independently by a separate agent, with no need for mutual communication. The decomposition discretizes the nodes in such a way that mesh conformity is maintained at the end. 

Note that our definition of block decomposition is more relaxed than the traditional sense, which expects perfect rectangular blocks to be generated~\citep{diprete2023reinforcementlearningblockdecomposition}. Since the downstream mesher demonstrates impressive robustness to oblique shapes, the primary goal of this agent is to remove reflex angles without creating any slivers/bad elements, rather than creating rectangles. 

In fact, this relaxation creates more complexity as the agent can now cast rays in any direction to generate blocks, rather than only at right angles. Inspired by the pioneering work done by Silver et al. in using RL to play Atari \citep{mnih2013playingatarideepreinforcement}, this complexity is handled by introducing a hybrid, multimodal, vector-vision-based input that allows for human-like, nuanced decision-making by the agent. Unlike the meshing agent, which only observes a local view, this agent observes a local and a global view of the boundary with the help of a \textbf{siamese convolutional neural network}~\citep{bromley1993signature}, which establishes separate encoders for both views, enabling rich feature learning at both scales. 

\subsection{Environment Formulation}
We formulate the decomposition process as a Markov Decision Process (MDP) defined by the tuple $(\mathcal{S},~\mathcal{A},~\mathcal{P},~\mathcal{R},~\gamma)$, where the agent interacts with a dynamic geometric boundary $\partial \Omega$. The environment operates via a recursive decomposition strategy: at each step $t$, the agent selects a target sub-region (block) and a point of interest $v_0$ to execute a split, partitioning the domain into two simpler sub-domains. The points of interest are generated based on the internal angle being greater than a set threshold $\alpha_t$. 

\subsubsection{State Representation}

\begin{figure}[htpb]

     \centering
     \begin{minipage}{0.92\textwidth}
     \begin{subfigure}[b]{0.11\textwidth}
         \centering
         \includegraphics[width=\textwidth]{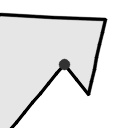}
     \end{subfigure}
     \hfill
     \begin{subfigure}[b]{0.11\textwidth}
         \centering
         \includegraphics[width=\textwidth]{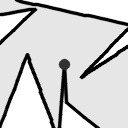}
     \end{subfigure}
     \hfill
     \begin{subfigure}[b]{0.11\textwidth}
         \centering
         \includegraphics[width=\textwidth]{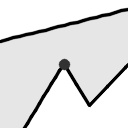}
     \end{subfigure}
     \hfill
     \begin{subfigure}[b]{0.11\textwidth}
         \centering
         \includegraphics[width=\textwidth]{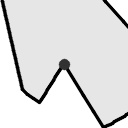}
     \end{subfigure}
     \hfill
     \begin{subfigure}[b]{0.11\textwidth}
         \centering
         \includegraphics[width=\textwidth]{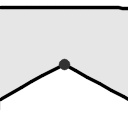}
     \end{subfigure}
     \hfill
     \begin{subfigure}[b]{0.11\textwidth}
         \centering
         \includegraphics[width=\textwidth]{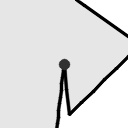}
     \end{subfigure}
     \hfill
     \begin{subfigure}[b]{0.11\textwidth}
         \centering
         \includegraphics[width=\textwidth]{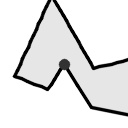}
     \end{subfigure}
     \hfill
     \begin{subfigure}[b]{0.11\textwidth}
         \centering
         \includegraphics[width=\textwidth]{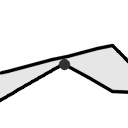}
     \end{subfigure}

     \vspace{10pt} 

     \begin{subfigure}[b]{0.11\textwidth}
         \centering
         \includegraphics[width=\textwidth]{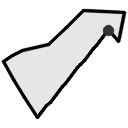}
     \end{subfigure}
     \hfill
     \begin{subfigure}[b]{0.11\textwidth}
         \centering
         \includegraphics[width=\textwidth]{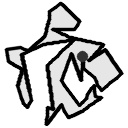}
     \end{subfigure}
     \hfill
     \begin{subfigure}[b]{0.11\textwidth}
         \centering
         \includegraphics[width=\textwidth]{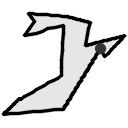}
     \end{subfigure}
     \hfill
     \begin{subfigure}[b]{0.11\textwidth}
         \centering
         \includegraphics[width=\textwidth]{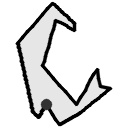}
     \end{subfigure}
     \hfill
     \begin{subfigure}[b]{0.11\textwidth}
         \centering
         \includegraphics[width=\textwidth]{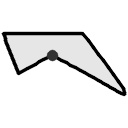}
     \end{subfigure}
     \hfill
     \begin{subfigure}[b]{0.11\textwidth}
         \centering
         \includegraphics[width=\textwidth]{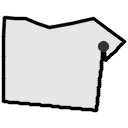}
     \end{subfigure}
     \hfill
     \begin{subfigure}[b]{0.11\textwidth}
         \centering
         \includegraphics[width=\textwidth]{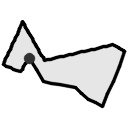}
     \end{subfigure}
     \hfill
     \begin{subfigure}[b]{0.11\textwidth}
         \centering
         \includegraphics[width=\textwidth]{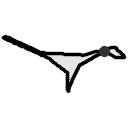}
     \end{subfigure}
     \end{minipage}
     \caption{Representative set of geometries used during the ensemble-training phase. The agent is exposed to domains with varying curvature, concavity, aspect ratios, and local topological complexity, enabling the policy to generalize beyond the simple square geometry used during curriculum pre-training.}
     \label{fig:block_decomposition_grid}
\end{figure}

To capture both local geometric fidelity and global topological context, the block decomposition agent employs a hybrid grayscale state representation $\mathbf{s}_t = (\mathbf{x}_{vec}, \mathbf{x}_{img})$ consisting of a rotationally invariant vector component and a visual component. We source a dataset of geometries with highly convex angles from Tong et al.~\citep{tong2023srl}. Being a visual CNN, this network requires far more data to achieve a good policy. The vector state $\mathbf{x}_{vec}$ encodes the immediate neighborhood of the focus vertex $v_0$, and is the same representation as we discussed in \ref{sssec:MeshStateRep}. 

While a highly partial view suffices for the meshing agent, the block decomposition agent requires greater fidelity to determine the precise angle at which a splitting ray should be cast. Fan-sensing rays alone tend to miss crucial geometric features due to their limited angular coverage, motivating the introduction of a visual state $\mathbf{x}_{img}$ consisting of a dual-channel rasterized image generated via OpenCV. The first channel is a high-resolution local crop centered at $v_0$ and oriented along the bisector, capturing detailed vertex configurations in the immediate vicinity. The second channel is a downsampled global representation of the entire domain $\Omega$, providing context for the overall shape. Two separate 128$\times$128 images are used rather than a single high-resolution image to prevent over-parameterization of the CNN feature extractor while still conveying both local detail and global context.

The local images are rendered such that the agent is always in the center (denoted by the black dot), while the global images are rendered to give positional context for the agent. Both images maintain the same orientation. These consistencies are maintained to eliminate variation before passing into the pseudo-siamese encoder.

\subsubsection{Action Space and Steering}
The action space is continuous, $\mathbf{a}_t \in [-1, 1]$. The primary component, $a_0$, determines the splitting direction relative to the local bisector, once again reflecting the coordinate system of the meshing agent. The steering angle $\phi$ is derived as:
\begin{equation}
    \phi = a_0 \cdot \left( \frac{\alpha_{v_0}}{2} \right) \cdot \eta
\end{equation}
where $\alpha_{v_0}$ is the interior angle and $\eta <= 1$ is a safety scaling factor. Essentially, the agent can cast a ray across the entire interior angle range of the point of interest. A ray is cast from $v_0$ in the direction $\phi$. If the ray intersects the boundary at point $p_{int}$, the environment attempts to snap $p_{int}$ to the nearest existing vertex $v_{near}$ to avoid creating degenerate sliver elements, maintain conformity, and not change the input element count. 

\subsubsection{Topological Constraints and Robust Parity}
A critical requirement for all-quadrilateral meshing is that the number of vertices in any boundary loop must be even. We enforce this via a \textit{Robust Parity Check}. Let the domain be split into two loops, $\partial \Omega_1$ and $\partial \Omega_2$, by a cut $\mathcal{C}$ consisting of $m$ internal vertices. For a valid quad-meshable topology, the vertex count $|\partial \Omega_i|$ must satisfy:
\begin{equation}
    |\partial \Omega_i| \equiv 0 \pmod 2
\end{equation}
If a split results in an odd vertex count, the environment invokes a correction subroutine. Unlike naive approaches that might split the newly created cut $\mathcal{C}$, our method strictly searches for the longest edge $e \in \partial \Omega_{initial}$ (edges belonging to the original domain boundary) and bisects it. This preserves the geometric quality of the internal cuts while satisfying mesh conformity. 

\subsubsection{Reward Function}
The reward function $R(s, a)$ balances geometric quality with topological balance:
\begin{equation}
    R = w_q \cdot Q_{ortho} + w_b \cdot B_{split} + R_{progress} - P_{penalty}
\end{equation}
\noindent The orthogonality term $Q_{ortho}$ measures the deviation of the split angles from $\pi/2$ at both endpoints and carries the largest weight, since maintaining angles close to $90^\circ$ or $180^\circ$ is essential for generating healthy blocks that the downstream meshing agent can process effectively. The balance term $B_{split}$ is the ratio of the perimeter lengths of the two resulting sub-loops, encouraging roughly equal partitions. A sparse progress reward $R_{progress}$ is triggered whenever a reflex vertex (a ``point of interest'') is eliminated by the split. The environment monitors and rewards the new angles created by each cut, providing the agent with a dense learning signal.

\subsection{Block Decomposition Agent}

\begin{figure}[htpb]
     \centering
     \includegraphics[width=0.75\textwidth]{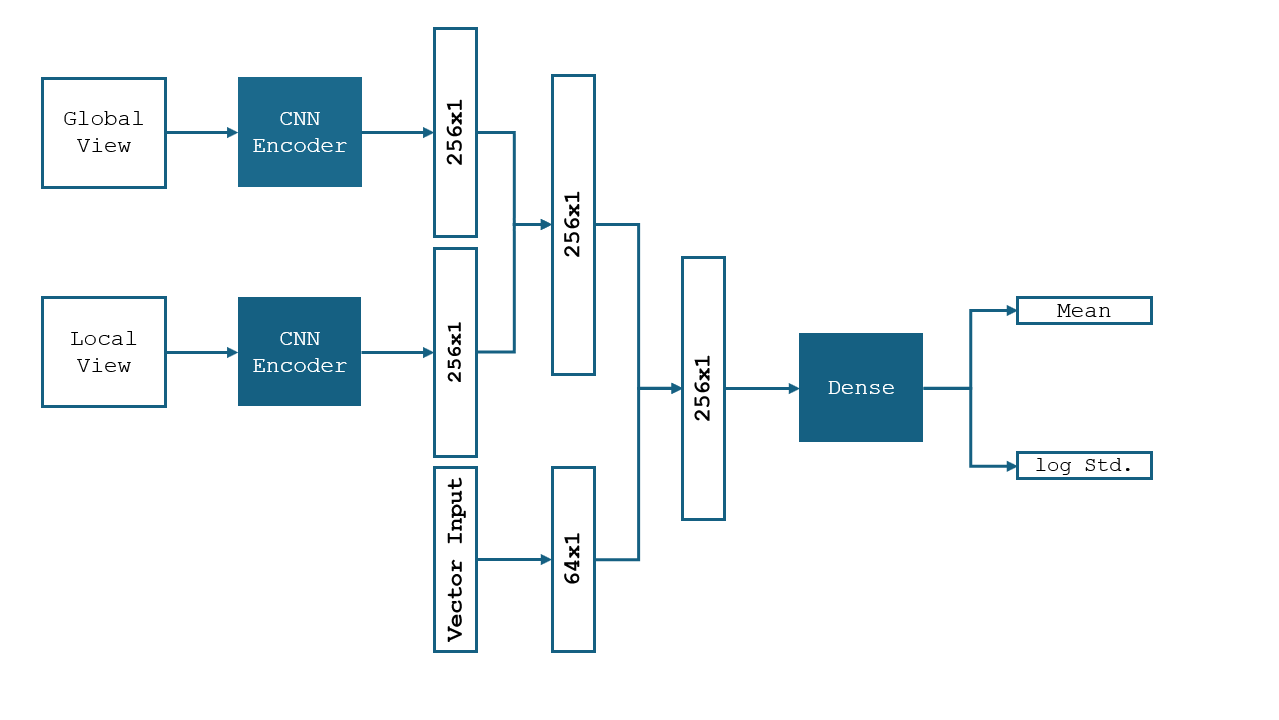}
     \caption{Architecture of the block decomposition reinforcement learning agent. Local and global rasterized geometry representations are processed by separate convolutional encoders and fused with geometric vector features to form a hybrid latent representation. The resulting embedding is used by the Soft Actor-Critic policy and critic networks to predict splitting directions that simplify the geometry into meshable subdomains.}
\end{figure}

The block decomposition agent utilizes a Soft Actor-Critic (SAC) framework~\citep{haarnoja2018soft} adapted for hybrid state spaces, processing high-dimensional visual inputs alongside low-dimensional vector coordinates. The state at time step $t$ is defined as a tuple $s_t = (I_t, \mathbf{v}_t)$, where $I_t \in \mathbb{R}^{2 \times H \times W}$ represents the visual input comprising the local and global single-channel views ($I_t^{loc}$, $I_t^{glob}$), and $\mathbf{v}_t \in \mathbb{R}^{d_v}$ represents the coordinate vector state. The architecture consists of a shared feature-extraction strategy employed across decoupled actor and critic networks. Due to the high fidelity, multimodal input, the block decomposition agent is the biggest among the 3 in the pipeline, with the actor requiring around 1.5M parameters end-to-end. 

\subsubsection{Feature Extraction and Policy Network}

The visual input $I_t$ is processed by a pseudo-siamese dual encoder that, unlike standard siamese networks with shared weights, employs two distinct Convolutional Neural Networks (CNNs) to capture features specific to the local and global perspectives. Let $f_{\phi_{loc}}$ and $f_{\phi_{glob}}$ denote the local and global CNNs, respectively. Each network consists of three convolutional layers with $3 \times 3$ kernels and a stride of 2, followed by a flattening operation, a linear projection, and Layer Normalization. The visual embedding $\mathbf{z}_{img}$ is obtained by concatenating and fusing the two branch outputs:
\begin{equation}
\begin{aligned}
    &\mathbf{h}_{loc} = \text{LayerNorm}(\text{ReLU}(f_{\phi_{loc}}(I_t^{loc}))) \\
    &\mathbf{h}_{glob} = \text{LayerNorm}(\text{ReLU}(f_{\phi_{glob}}(I_t^{glob}))) \\
    &\mathbf{z}_{img} = \text{ReLU}(\mathbf{W}_{fusion}[\mathbf{h}_{loc} \oplus \mathbf{h}_{glob}]),
\end{aligned}
\end{equation}
where $\oplus$ denotes concatenation and $\mathbf{W}_{fusion}$ is a learned linear transformation. This produces a high-dimensional feature representation of the visual observation.

The policy network $\pi_\theta(a_t | s_t)$ maps the hybrid state to a probability distribution over actions. The vector state $\mathbf{v}_t$ is encoded via a separate Multi-Layer Perceptron (MLP), $g_\theta$, yielding $\mathbf{z}_{vec}$. A joint state representation is formed by fusing the visual and vector embeddings, $\mathbf{z}_{joint} = [\mathbf{z}_{img} \oplus \mathbf{z}_{vec}]$, which is then passed through a trunk MLP that bifurcates into two heads predicting the mean $\mu_\theta$ and the log-standard deviation $\log \sigma_\theta$ of a Gaussian distribution. The action $a_t$ is sampled using the reparameterization trick:
\begin{equation}
    a_t = \tanh(\mu_\theta(s_t) + \sigma_\theta(s_t) \cdot \epsilon), \quad \epsilon \sim \mathcal{N}(0, I).
\end{equation}

\subsubsection{Critic Network and Training}

Twin critics $Q_{\psi_1}$ and $Q_{\psi_2}$ are employed to mitigate value overestimation, with the target network updated via Polyak averaging. Crucially, the critics maintain their own separate instances of the image and vector encoders, ensuring that gradients from the Bellman updates do not interfere with the policy's feature learning. Each critic approximates the Q-value by fusing the state features with the action:
\begin{equation}
    Q_{\psi_i}(s_t, a_t) = \text{MLP}_{\psi_i}([\mathbf{z}_{img} \oplus \mathbf{z}_{vec} \oplus a_t])
\end{equation}
for $i \in \{1, 2\}$. The critics are trained to minimize the mean squared soft Bellman residual:
\begin{equation}
    J_Q(\psi_i) = \mathbb{E}_{(s_t, a_t) \sim \mathcal{D}} \left[ (Q_{\psi_i}(s_t, a_t) - y_t)^2 \right]
\end{equation}
where the target $y_t = r_t + \gamma (1 - d_t) (\min_{j=1,2} Q_{\bar{\psi}_j}(s_{t+1}, a_{t+1}) - \alpha \log \pi_\theta(a_{t+1}|s_{t+1}))$. The actor parameters are optimized to maximize the expected return while encouraging entropy:
\begin{equation}
    J_\pi(\theta) = \mathbb{E}_{s_t \sim \mathcal{D}} \left[ \alpha \log \pi_\theta(a_t|s_t) - \min_{j=1,2} Q_{\psi_j}(s_t, a_t) \right]
\end{equation}
Finally, the temperature parameter $\alpha$ is automatically tuned to satisfy a target entropy constraint $\bar{\mathcal{H}}$ by minimizing:
\begin{equation}
    J(\alpha) = \mathbb{E}_{a_t \sim \pi_t} [-\alpha(\log \pi_t(a_t|s_t) + \bar{\mathcal{H}})]
\end{equation}

\section{Hole Decomposition}

\begin{figure}[htpb]

     \centering
     \begin{minipage}{0.92\textwidth}
     \begin{subfigure}[b]{0.11\textwidth}
         \centering
         \includegraphics[width=\textwidth]{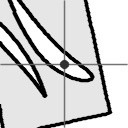}
     \end{subfigure}
     \hfill
     \begin{subfigure}[b]{0.11\textwidth}
         \centering
         \includegraphics[width=\textwidth]{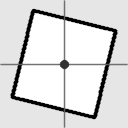}
     \end{subfigure}
     \hfill
     \begin{subfigure}[b]{0.11\textwidth}
         \centering
         \includegraphics[width=\textwidth]{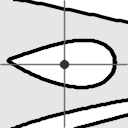}
     \end{subfigure}
     \hfill
     \begin{subfigure}[b]{0.11\textwidth}
         \centering
         \includegraphics[width=\textwidth]{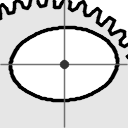}
     \end{subfigure}
     \hfill
     \begin{subfigure}[b]{0.11\textwidth}
         \centering
         \includegraphics[width=\textwidth]{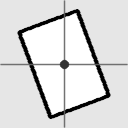}
     \end{subfigure}
     \hfill
     \begin{subfigure}[b]{0.11\textwidth}
         \centering
         \includegraphics[width=\textwidth]{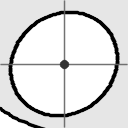}
     \end{subfigure}
     \hfill
     \begin{subfigure}[b]{0.11\textwidth}
         \centering
         \includegraphics[width=\textwidth]{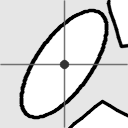}
     \end{subfigure}
     \hfill
     \begin{subfigure}[b]{0.11\textwidth}
         \centering
         \includegraphics[width=\textwidth]{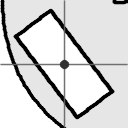}
     \end{subfigure}

     \vspace{10pt} 

     \begin{subfigure}[b]{0.11\textwidth}
         \centering
         \includegraphics[width=\textwidth]{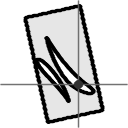}
     \end{subfigure}
     \hfill
     \begin{subfigure}[b]{0.11\textwidth}
         \centering
         \includegraphics[width=\textwidth]{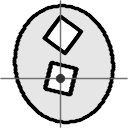}
     \end{subfigure}
     \hfill
     \begin{subfigure}[b]{0.11\textwidth}
         \centering
         \includegraphics[width=\textwidth]{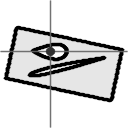}
     \end{subfigure}
     \hfill
     \begin{subfigure}[b]{0.11\textwidth}
         \centering
         \includegraphics[width=\textwidth]{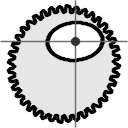}
     \end{subfigure}
     \hfill
     \begin{subfigure}[b]{0.11\textwidth}
         \centering
         \includegraphics[width=\textwidth]{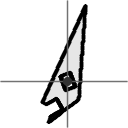}
     \end{subfigure}
     \hfill
     \begin{subfigure}[b]{0.11\textwidth}
         \centering
         \includegraphics[width=\textwidth]{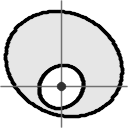}
     \end{subfigure}
     \hfill
     \begin{subfigure}[b]{0.11\textwidth}
         \centering
         \includegraphics[width=\textwidth]{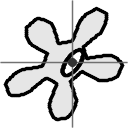}
     \end{subfigure}
     \hfill
     \begin{subfigure}[b]{0.11\textwidth}
         \centering
         \includegraphics[width=\textwidth]{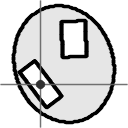}
     \end{subfigure}
     \end{minipage}
     \caption{Dual-scale visual observation used by the hole decomposition agent. The global image encodes the overall geometry and hole configuration, while the local image provides a magnified view centered at the pole of inaccessibility inside the target hole. This representation enables the agent to infer geometrically valid cutting directions independent of the underlying vertex count.}
     \label{fig:local_global_vision}
\end{figure}

To make this agent truly scalable and create an opportunity for end-to-end parallelization, it emerged that another framework would be required to split geometries with holes into blocks. Hence, another agent is trained specifically to split hole-geometries should they be present. This problem warrants another agent to perfectly cut a hole; the agent is best placed inside the hole instead of on the boundary. We restrict ourselves purely to the visual state, rather than the hybrid state as viewed by the block decomposition agent. 

\subsection{Environment Formulation}
\label{sec:env_formulation}

The proposed environment models the hole decomposition task as a sequential decision-making process. The objective is to decompose a complex polygon $\mathcal{P}_{init}$ containing holes $\mathcal{H} = \{h_1, \dots, h_n\}$ into a set of convex, rectangular-like blocks $\mathcal{B}$ suitable for high-quality mesh generation. Crucially, if the agent creates sharp angles of incidence, they will severely affect the quality of elements generated by the meshing module downstream. We model our reward function around this crucial condition.

\subsubsection{Observation Space}
The observation space is image-based to allow the agent to learn geometric features agnostic to vertex count, thereby enabling scale invariance. At each time step $t$, the state $S_t$ consists of two stacked binary masks representing local and global contexts:
\begin{equation}
    S_t = \{I_{global}, I_{local}\} \in \{0, 1\}^{H \times W \times 2}
\end{equation}
The agent is positioned at a point $v_0$ located inside the current target hole. To ensure robust positioning distinct from the polygon centroid (which may fall outside concave geometries), the \textit{Pole of Inaccessibility} is used---defined as the center of the largest inscribed circle within the hole $h_i$. This choice also avoids sensitivity to non-uniform boundary discretizations that can shift the geometric centroid. The global channel $I_{global}$ captures the entire bounding box of the remaining geometry, normalized to the image size, while the local channel $I_{local}$ is a cropped view centered on $v_0$ with a zoom factor of $5\times$ (covering 20\% of the global span) to resolve fine details near the cutting site. We add a crosshair to further promote orientation, since this agent does not have coordinate information that the hole decomposition agent has.

\subsubsection{Action Space, Dual-Ray Cutting, and Parity Enforcement}
The action space is continuous and defined as $\mathcal{A} \in [-1, 1]^2$. The agent selects two cutting angles relative to the geometric center $v_0$:
\begin{equation}
    a_t = [\theta_1, \theta_2] \quad \text{where} \quad \theta_i \in [-\pi, \pi]
\end{equation}
These angles define two ray directions $\vec{d}_1, \vec{d}_2$. The environment computes the intersection of each ray with the polygon boundary $\partial \mathcal{P}$:
\begin{equation}
    p_i = v_0 + t_i \vec{d}_i, \quad t_i = \min \{t > 0 \mid (v_0 + t\vec{d}_i) \in \partial \mathcal{P}\}
\end{equation}
A topological cut is then performed by constructing a cutter string via \texttt{Shapely} operations, partitioning the domain into two blocks per step. These blocks are evaluated by the environment, and an appropriate reward is assigned.

As with the block decomposition environment, an even-vertex parity constraint is enforced to facilitate downstream quad meshing. For every generated polygon $P$ whose vertex count $|V_P|$ is odd, the environment performs an implicit edge split. The splitting logic prioritizes geometric fidelity: letting $E_P$ denote the edges of $P$ and $E_{init}$ the edges of the initial geometry $\mathcal{P}_{init}$, the split candidate $e^*$ is selected as
\begin{equation}
    e^* = \operatorname*{argmax}_{e \in E_P} \text{length}(e) \quad \text{subject to} \quad \text{dist}(e_{mid}, E_{init}) < \epsilon
\end{equation}
This ensures that auxiliary vertices are added only to original boundaries (outer walls or original holes) and never to the internal cut interfaces created by the agent, thereby preserving mesh conformity.

\subsubsection{Reward Function}
The reward function $R(s, a)$ is a hierarchical sum of immediate action quality and terminal shape quality:
\begin{equation}
    R_{total} = R_{action} + \sum_{b \in \mathcal{B}_{new}} R_{outcome}(b)
\end{equation}
\noindent The action reward $R_{action}$ evaluates the geometric validity of the cut before topological processing. It applies a Gaussian penalty when the angle between the two rays $|\theta_1 - \theta_2|$ falls below $15^\circ$, discouraging sliver cuts, and provides a positive signal based on the angle $\phi$ between each ray and the impacted wall normal, reaching its maximum at $\phi = 0$ (orthogonal impact).

The outcome reward $R_{outcome}$ evaluates the quality of each finalized block $b$ by combining three metrics:
\begin{equation}
    R_{outcome}(b) = w_1 \cdot \mathcal{S}(b) + w_2 \cdot \mathcal{A}_{\text{score}}(b) + w_3 \cdot \mathcal{R}_{asp}(b)
\end{equation}
\noindent Here the convexity term $\mathcal{S}(b) = \text{Area}(b) / \text{Area}(\text{ConvexHull}(b))$ measures how close the block is to its convex hull, the angle score $\mathcal{A}_{\text{score}}$ penalizes internal angles below $30^\circ$ while rewarding angles near $90^\circ$, and the aspect ratio $\mathcal{R}_{asp}$ is defined as the ratio of the shortest to the longest edge. This reward formulation extends the methodology established for the block decomposition agent; however, the transition from a single-ray to a dual-ray action space necessitates a recalibration of the weighting terms. In particular, greater emphasis is placed on the convexity coefficient $w_1$, since simultaneous dual-ray casting partitions the geometry into more complex sub-regions than single cuts, and the amplified convexity constraint is needed to ensure geometric stability and prevent the generation of non-convex blocks.

\subsection{Hole Decomposition Agent}
\label{subsec:hole_agent}

\begin{figure}[htpb]
     \centering
     \includegraphics[width=0.75\textwidth]{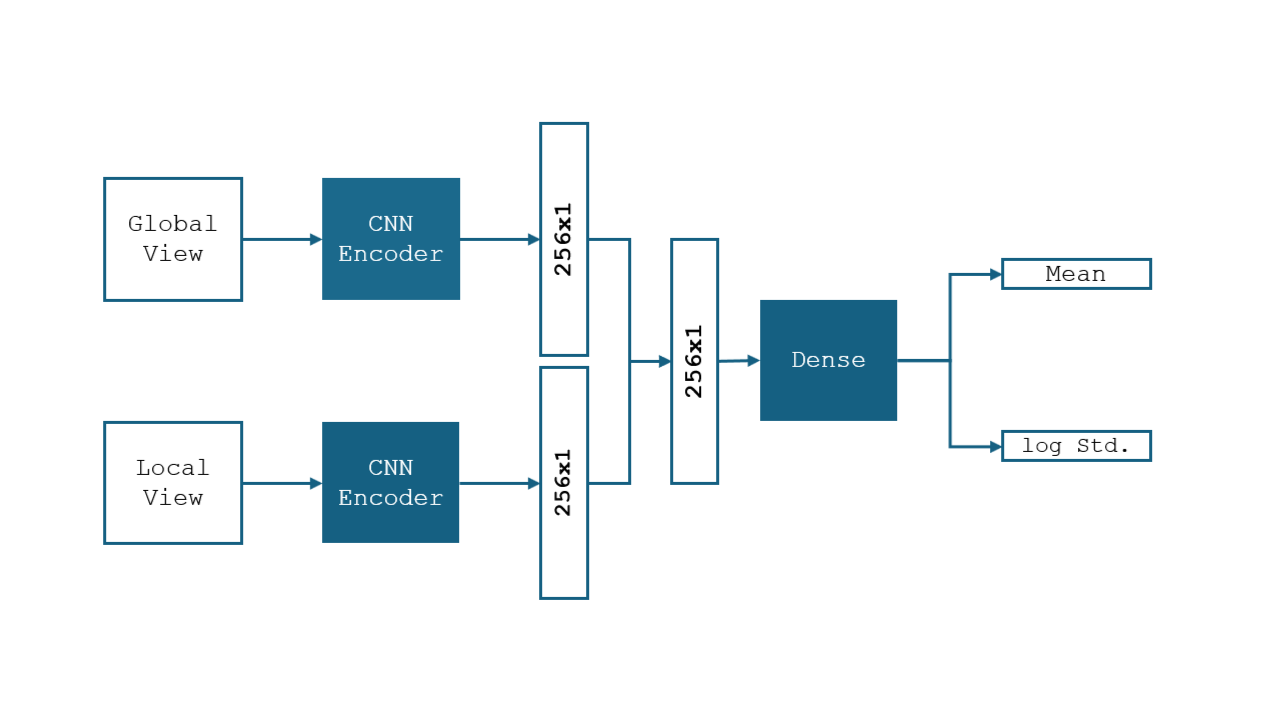}
     \caption{Vision-based architecture of the hole decomposition agent. Local and global binary masks are independently processed by dual convolutional encoders and fused into a shared latent representation. The Soft Actor-Critic policy predicts two cutting-ray directions used to partition hole-containing geometries into simpler meshable regions.}
     \label{fig:vision_based_arch}
\end{figure}

To tackle the specific topology of hole-containing geometries, a modified Soft Actor-Critic agent~\citep{haarnoja2018soft} is employed that relies exclusively on visual perception. Unlike the block decomposition agent, which benefits from hybrid vector-visual states, the hole decomposition agent operates from the Pole of Inaccessibility, a point strictly internal to the void (see Figure~\ref{fig:vision_based_arch}). In this context, relative coordinate vectors to boundary nodes are less informative than the holistic shape of the surrounding void, so the vector state encoder is discarded in favour of a pure vision-based architecture. This also marginally lightens the architecture, coming in at approx 1.3M parameters for the actor. 

The agent observes the state $s_t$ as a dual-channel tensor $I_t \in \mathbb{R}^{2 \times H \times W}$, stacking the local and global binary masks described in Section~\ref{sec:env_formulation}. These channels are processed by the same dual-stream CNN encoder used in the block decomposition agent:
\begin{equation}
    \mathbf{z}_{local} = f_{\theta_{loc}}(I_{t}^{local}), \quad \mathbf{z}_{global} = f_{\theta_{glob}}(I_{t}^{global}),
\end{equation}
where each $f_{\theta}$ is a 3-layer convolutional neural network following the architecture of Mnih et al.~\citep{mnih2015human} (64 filters, stride 2), followed by Layer Normalization. The resulting embeddings are concatenated and fused via a dense projection layer to form the state embedding:
\begin{equation}
    \mathbf{h}_t = \text{ReLU}(\mathbf{W}_{fusion}[\mathbf{z}_{local} \oplus \mathbf{z}_{global}]).
\end{equation}

This visual embedding $\mathbf{h}_t$ serves as the sole input to the actor and critic networks. The policy $\pi_\phi(a_t | s_t)$ outputs the mean and log-standard deviation for the dual-ray cutting angles $a_t \in \mathbb{R}^2$, and the learning process follows the standard SAC formulation (Eq.~\ref{eq:sac_objective}--\ref{eq:alpha_loss}), utilizing automatic entropy tuning to prevent premature convergence to suboptimal cutting strategies. By relying solely on high-dimensional visual inputs, the agent learns to recognize geometric primitives, such as L-shapes or U-shapes, and orient the cut lines to maximize convexity, regardless of the specific vertex count or scale of the hole.

\section{Results}

We begin by comparing the raw results of our meshing solution against the current state of the art, after which we demonstrate the output of each agent in the pipeline, as well as the combined results of the pipeline, and provide relevant numerical results. 

\subsection{Curriculum Learning}
\begin{figure}[htpb]
     \centering
     \begin{subfigure}[b]{0.45\textwidth}
         \centering
         \includegraphics[width=\textwidth]{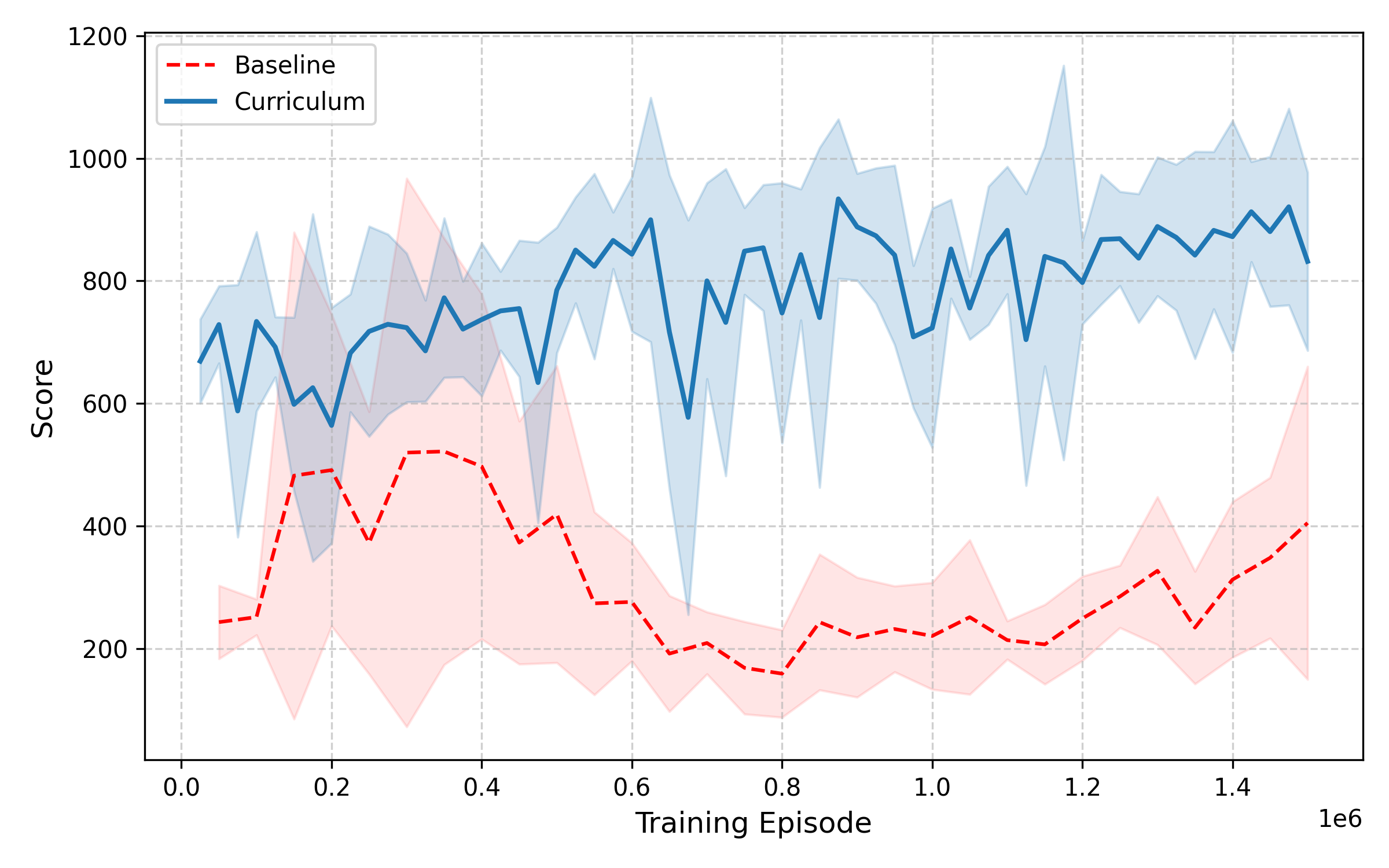}
         \caption{Scaled Score (higher is better)}
     \end{subfigure}
     \hfill
     \begin{subfigure}[b]{0.45\textwidth}
         \centering
         \includegraphics[width=\textwidth]{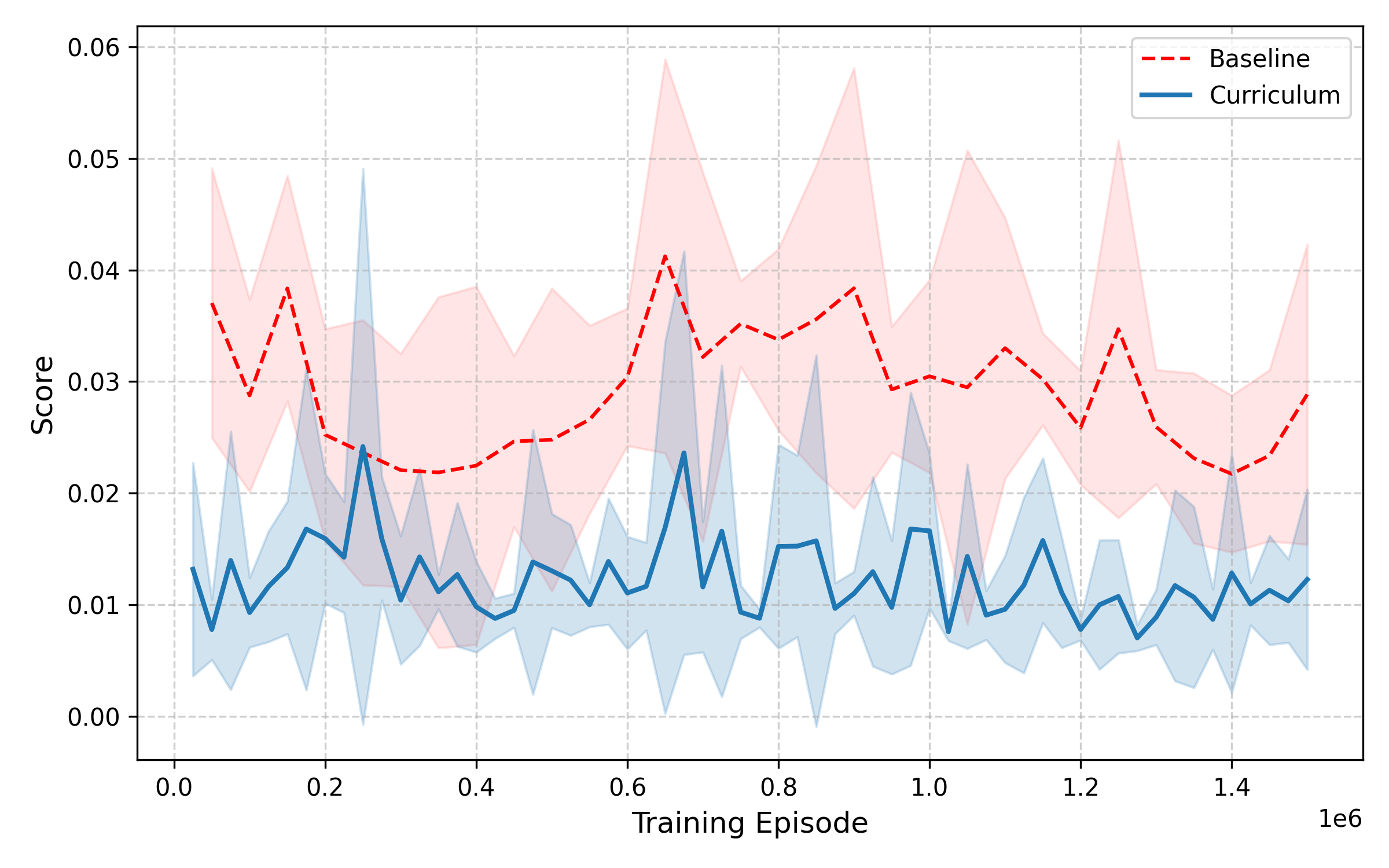}
         \caption{Element variance (lower is better)}
     \end{subfigure}
     \caption{Vision-based architecture of the hole decomposition agent. Local and global binary masks are independently processed by dual convolutional encoders and fused into a shared latent representation. The Soft Actor-Critic policy predicts two cutting-ray directions used to partition hole-containing geometries into simpler meshable regions.}
     \label{fig:seedavgruns}
\end{figure}

To demonstrate the importance of curriculum learning in achieving higher quality elements and faster convergence, we show a scaled score, calculated by multiplying average element quality by the number of elements generated by the meshing process, as well as element variance. The results displayed are averaged across 5 seeds each (see Figure~\ref{fig:seedavgruns}). 

\begin{figure}[htpb]
     \centering
     \begin{subfigure}[b]{0.45\textwidth}
         \centering
         \includegraphics[width=\textwidth]{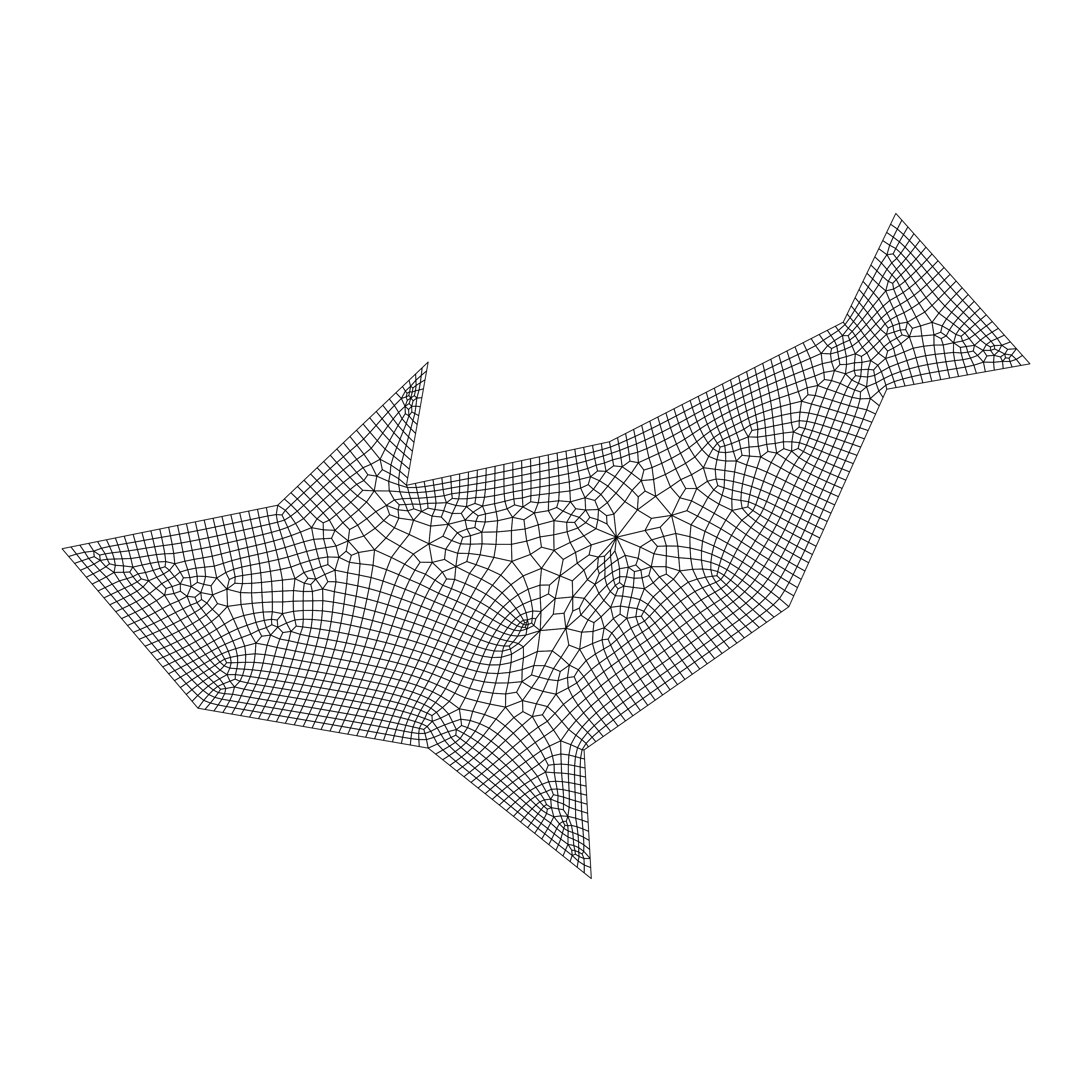}
         \caption{Baseline Training}
         \label{fig:left}
     \end{subfigure}
     \hfill
     \begin{subfigure}[b]{0.45\textwidth}
         \centering
         \includegraphics[width=\textwidth]{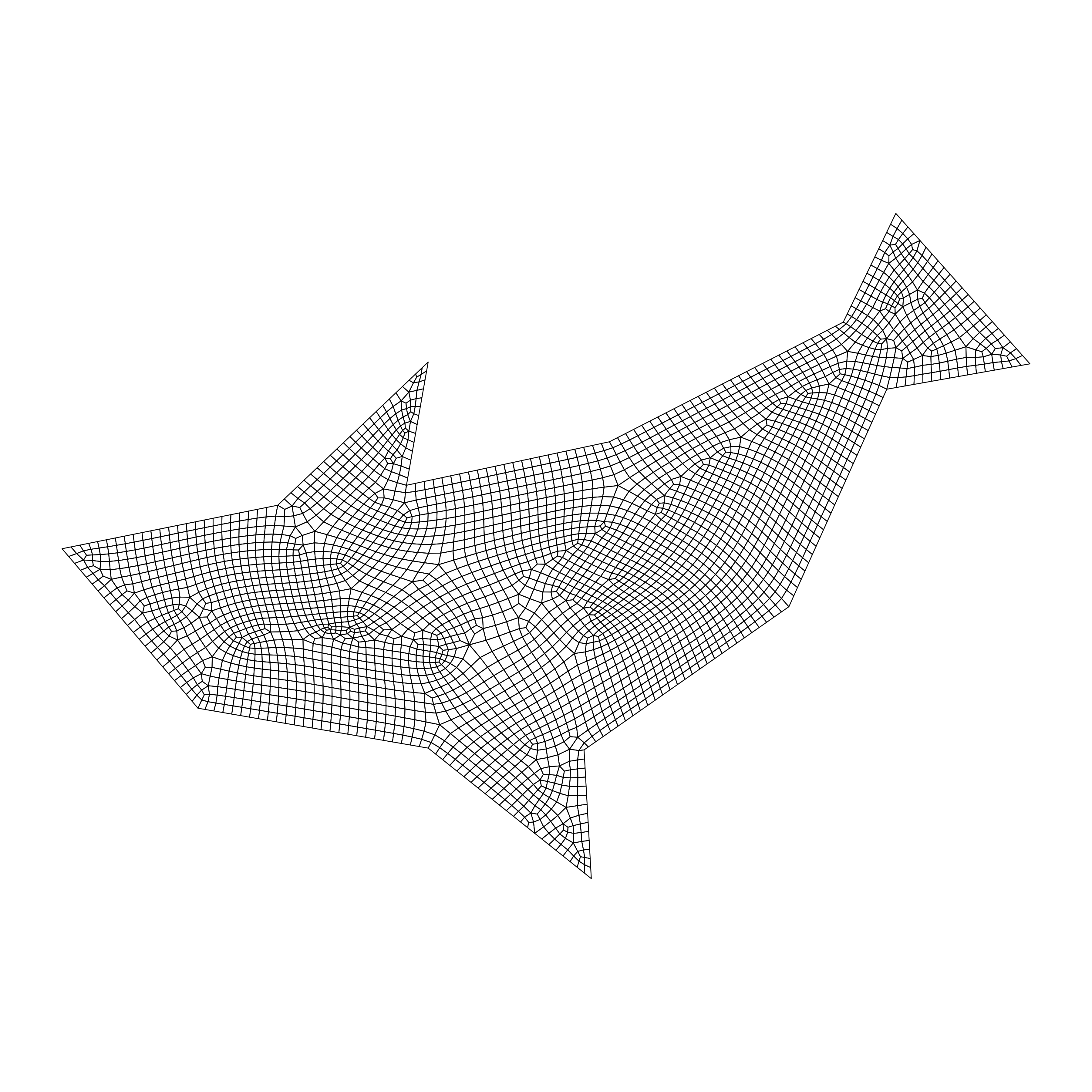}
         \caption{Curriculum Training}
         \label{fig:right}
     \end{subfigure}
        \caption{Qualitative comparison between a baseline agent and a curriculum-trained agent after 600k training steps. The curriculum-trained policy produces smoother advancing fronts, more uniform element sizing, and fewer irregular singular regions, demonstrating improved global planning and mesh completion behavior.}
        \label{fig:overall}
\end{figure}
As we can see, the consistent initial exposure to a simple square geometry allows for not only lower variance, thereby leading to a more consistent behavior, but also improves long-range planning, which allows for completion of the mesh almost 100\% of the time, which is not the case for the baseline method, as indicated by the scaled score.   

\begin{table}[htpb]
    \centering
    \caption{Mesh quality comparison between a baseline meshing agent (trained without curriculum learning) and the proposed curriculum-trained meshing agent, evaluated on identical test geometries and averaged across five independent random seeds. Quality is measured using the Scaled Jacobian metric, normalized to [0,1], where a value of 1.0 denotes a perfect square element. Higher mean and global minimum values, and lower variance, indicate superior mesh quality. The global minimum captures the worst single element produced across all test cases and is a critical indicator of mesh reliability.}
    \label{tab:my_table}
    \begin{tabular}{lccc} 
        \toprule
        \textbf{Method} & \textbf{Mean} & \textbf{Variance} & \textbf{Global Minimum} \\ 
        \midrule
        Baseline        & 0.9432        & 0.007930          & 0.2644                 \\
        Curriculum      & \textbf{0.9638}        &\textbf{ 0.003192}          & \textbf{0.4755}                 \\
        \bottomrule
    \end{tabular}
\end{table}
The consistency in element generation and minimization of singular nodes is apparent, as is the consistent density of elements all across the inner region (see Figure~\ref{fig:overall}). It is worth noting that despite the sub-par elements created by the agent in the baseline model, the environment still manages to recover a complete mesh, demonstrating the robustness of the recovery mechanisms built in.

\subsection{Block Decomposition}

\begin{figure}[htpb]
     \centering
     \begin{subfigure}[b]{0.3\textwidth}
         \centering
         \includegraphics[width=\textwidth]{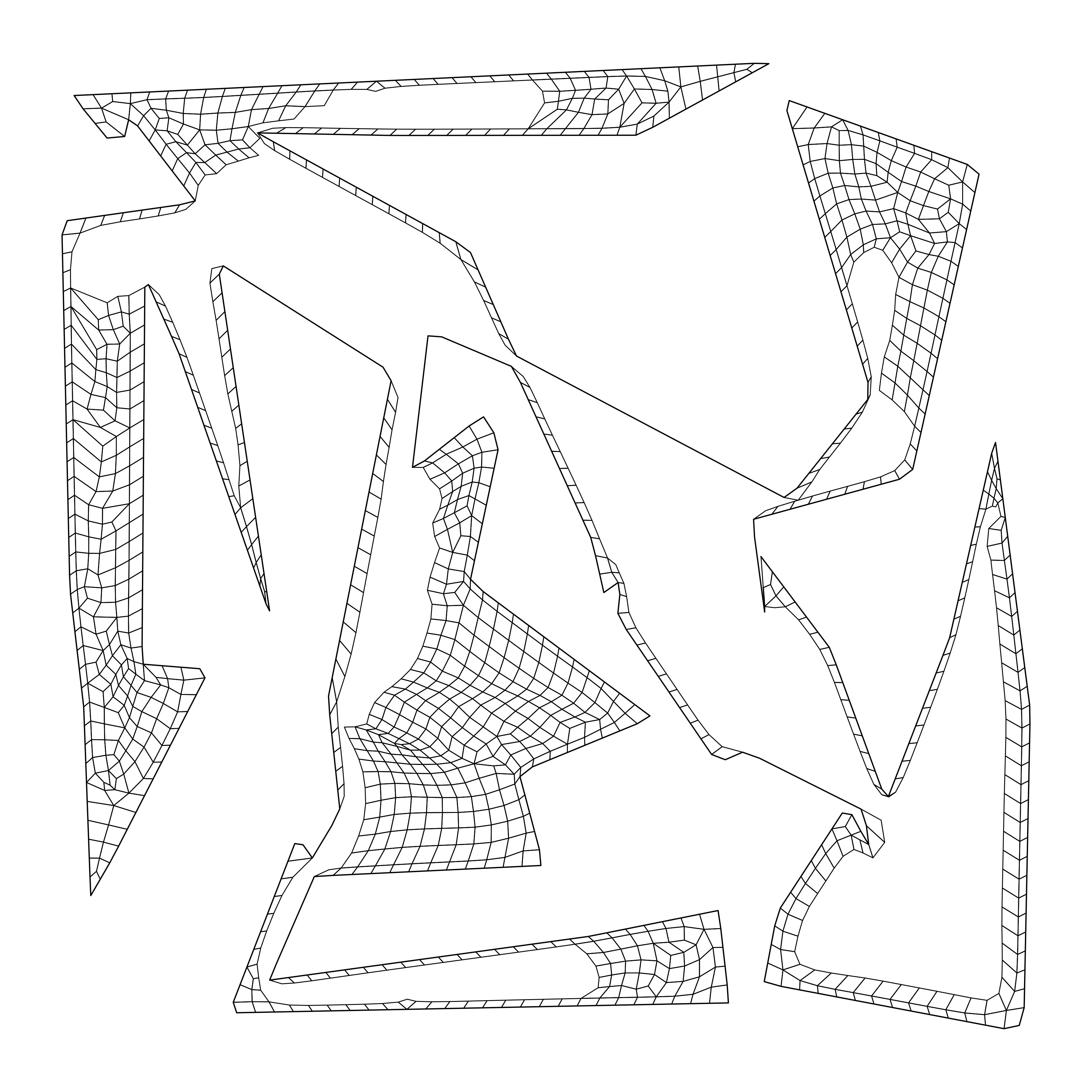}
         \caption{Incomplete Mesh}
     \end{subfigure}
     \hfill
     \begin{subfigure}[b]{0.3\textwidth}
         \centering
         \includegraphics[width=\textwidth]{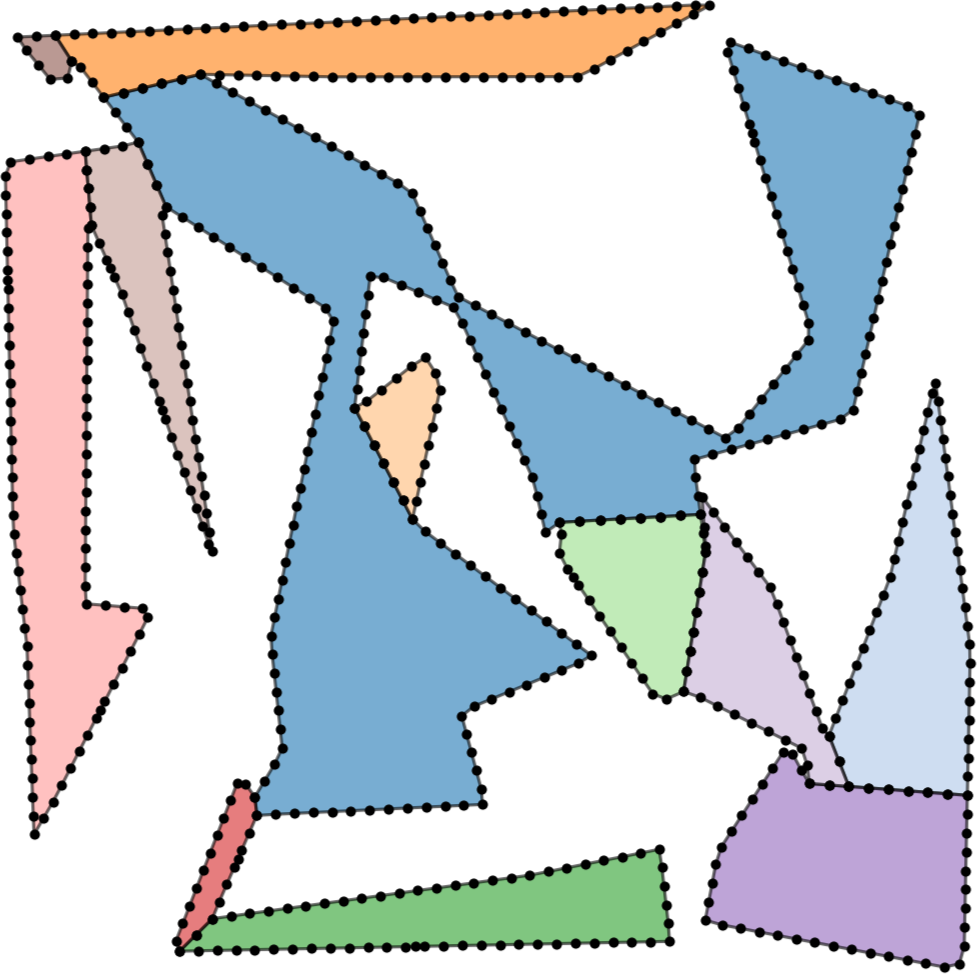}
         \caption{Block Decomposition}
     \end{subfigure}
     \hfill
     \begin{subfigure}[b]{0.3\textwidth}
         \centering
         \includegraphics[width=\textwidth]{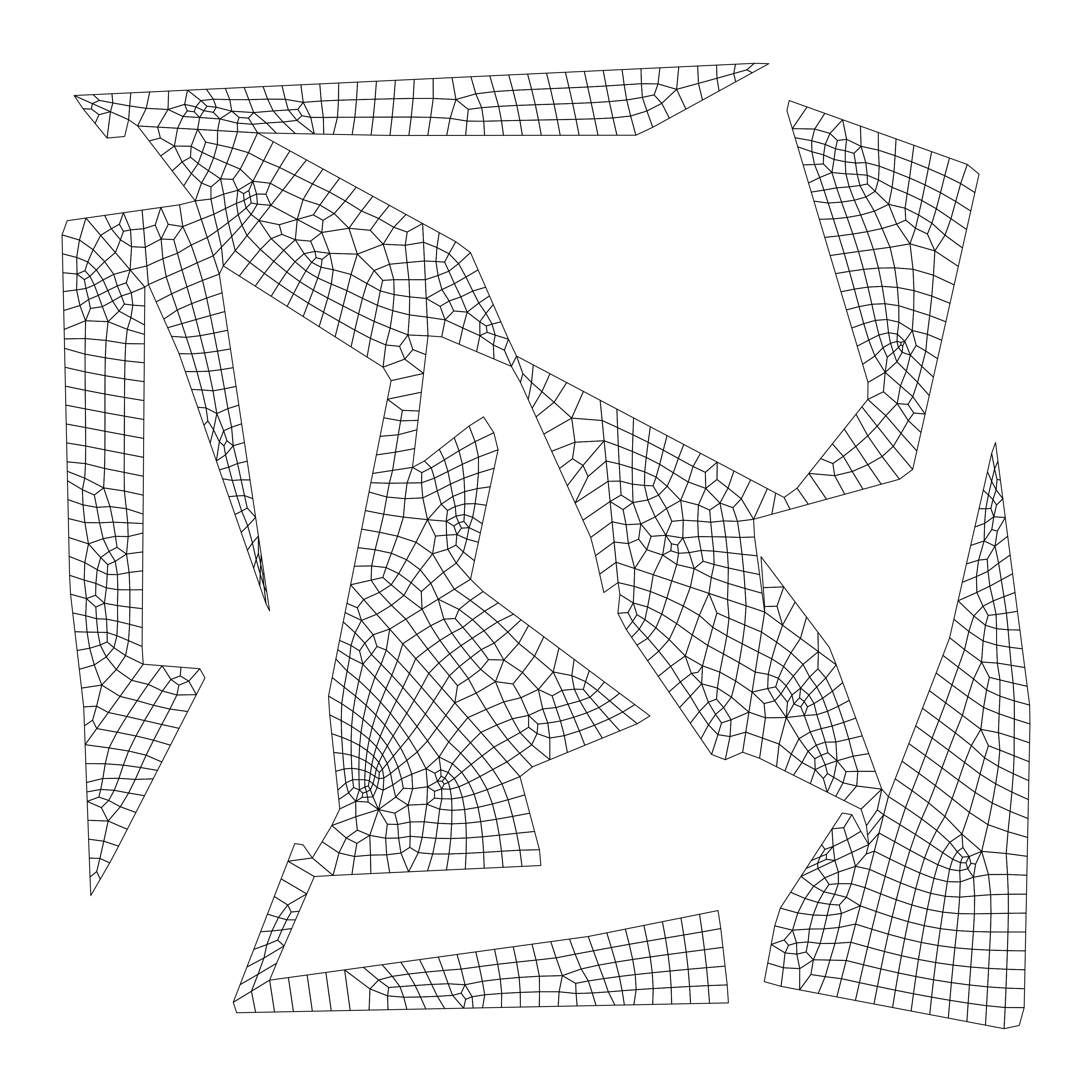}
         \caption{Complete Mesh}
     \end{subfigure}
      \caption{Effect of block decomposition on meshing robustness. Without decomposition (left), the advancing-front meshing agent fails to fully resolve the geometry due to the presence of strong concavities. The block decomposition agent (center) partitions the domain into simpler subregions, enabling successful generation of a complete all-quadrilateral mesh (right).}
     \label{fig:block_ability}
\end{figure}

The block decomposition agent allows the agent to mesh geometries it could not by itself, and allows for parallel processing of the individual elements, thereby speeding up inference time and allowing for multiple advancing fronts. We showcase the result of both aspects here (see Figure~\ref{fig:block_ability}). 

Note that the goal of the block decomposition agent is to eliminate large reflex angles that disrupt the meshing process while minimizing interference. This is because the increase in the number of blocks causes higher variance and a marginal decrease in element quality that the meshing agent otherwise smooths over.  
\begin{figure}[htpb]
     \centering
     \begin{subfigure}[b]{0.3\textwidth}
         \centering
         \includegraphics[width=\textwidth]{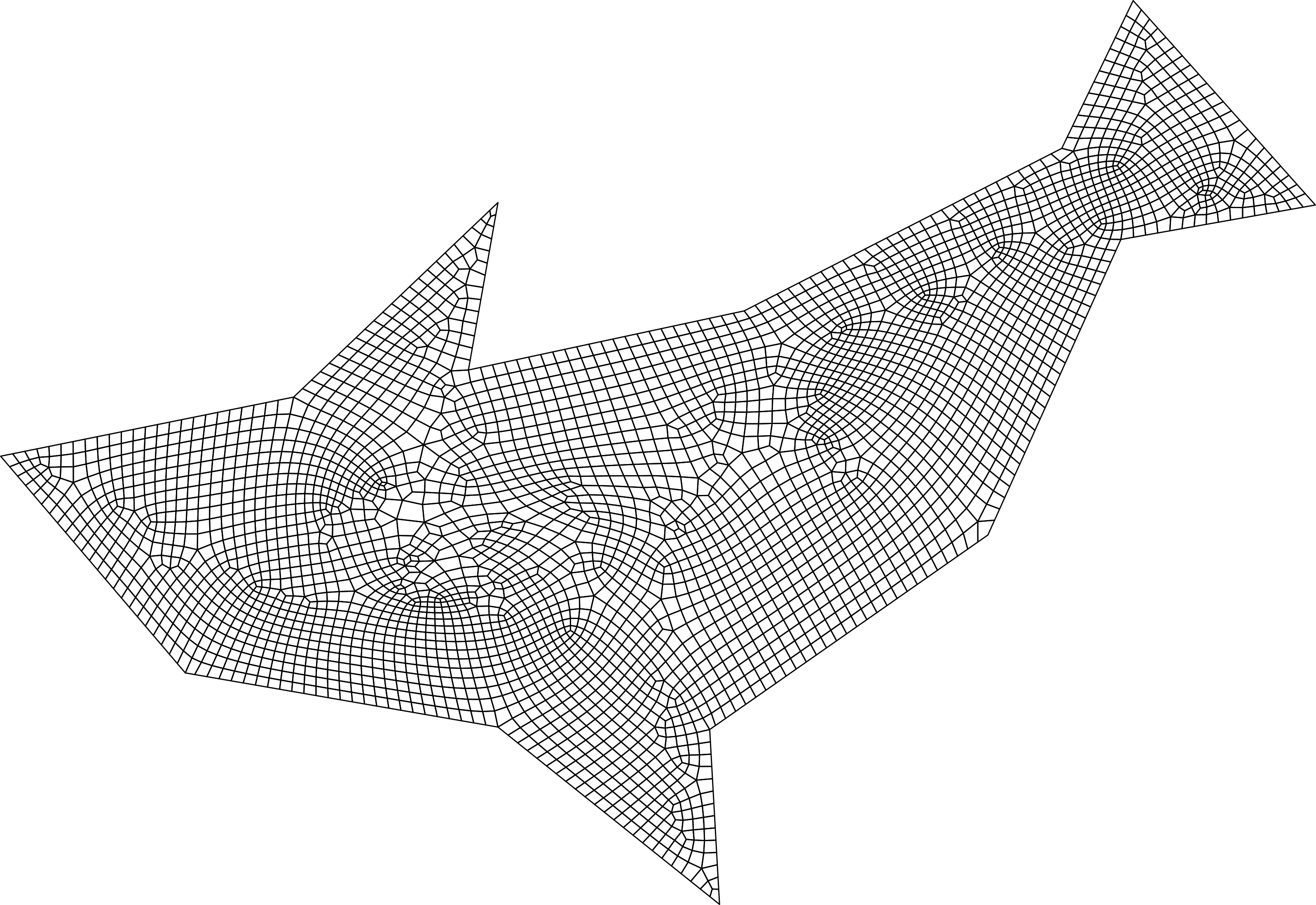}
         \caption{No Block Decomposition}
     \end{subfigure}
     \hfill
     \begin{subfigure}[b]{0.3\textwidth}
         \centering
         \includegraphics[width=\textwidth]{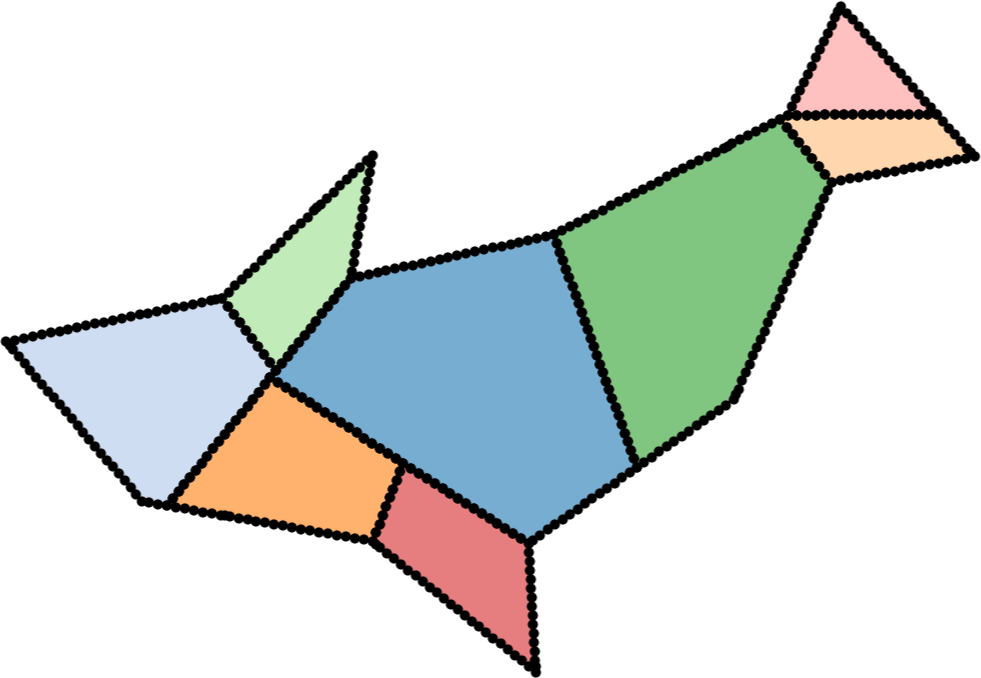}
         \caption{Block Decomposition}
     \end{subfigure}
     \hfill
     \begin{subfigure}[b]{0.3\textwidth}
         \centering
         \includegraphics[width=\textwidth]{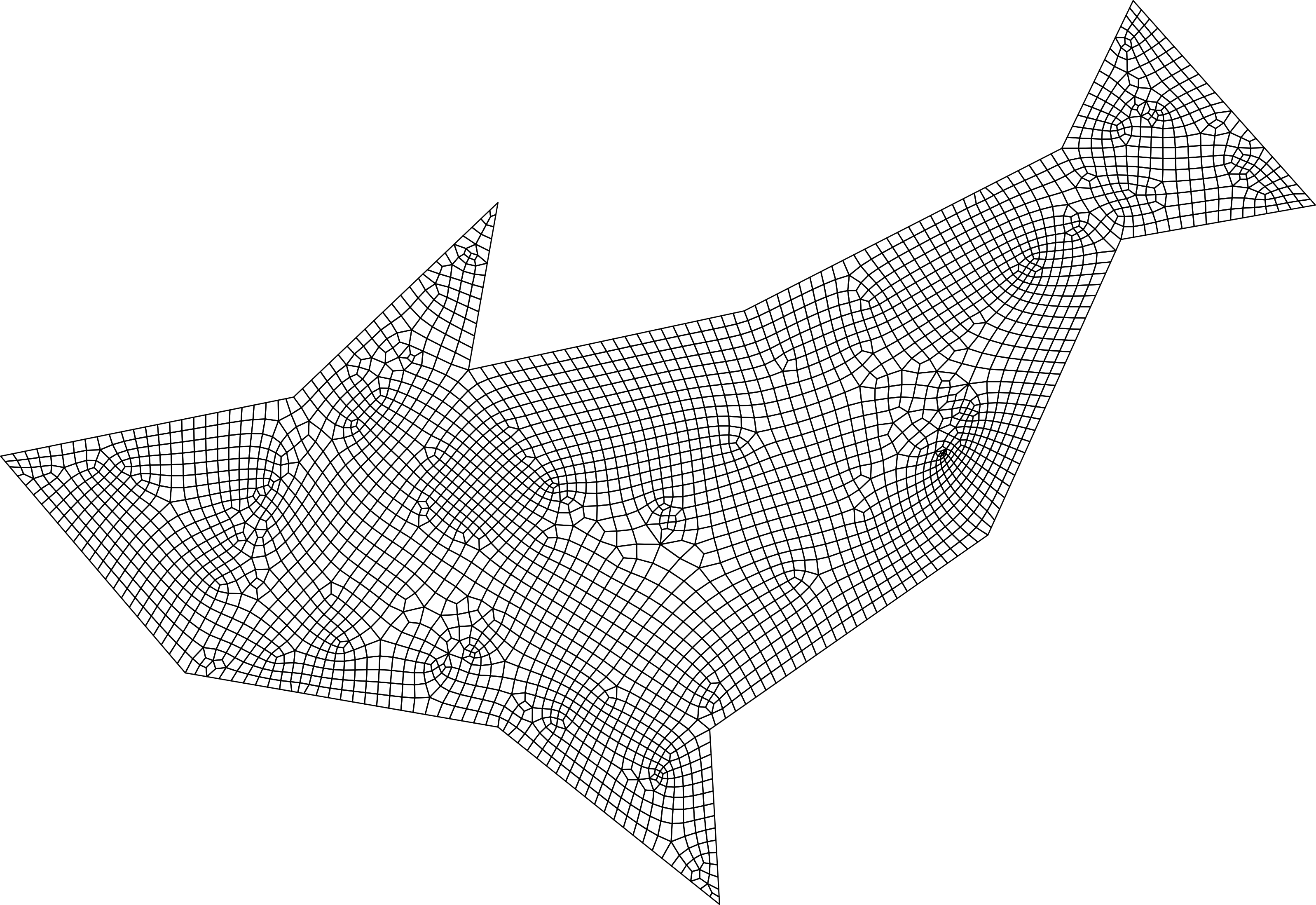}
         \caption{Parallel Mesh}
     \end{subfigure}
      \caption{Parallel meshing enabled through block decomposition. The original geometry (left) is first partitioned into independent subdomains (center), allowing multiple meshing agents to operate concurrently. The resulting mesh (right) demonstrates comparable element quality while significantly reducing overall meshing time.}
     \label{fig:parallel_ability}
\end{figure}

\begin{table}[H]
    \centering
    \caption{Quantitative comparison of meshing performance on a complex concave geometry processed either directly by a single meshing agent (No Block Decomposition) or first partitioned by the block decomposition agent and then meshed in parallel (Block Decomposition). Meshing time is reported in wall-clock seconds on identical hardware. Element quality is measured using the Scaled Jacobian metric; higher is better for mean and global minimum, lower is better for variance. The results illustrate the trade-off between computational speed-up through parallelism and marginal increases in element quality variance introduced by internal block boundaries.}
    \label{tab:comparison_w_and_w/o_blockdecomposition}
    \begin{tabular}{lcc} 
        \toprule
        \textbf{Entity} & \textbf{No Block Decomposition} & \textbf{Block Decomposition} \\ 
        \midrule
        Meshing time             & 61.87s              & 37.12s          \\ 
        Average Element Quality        & 0.9577               & 0.9545              \\ 
        Element Quality Variance             & 0.0040         & 0.0055              \\ 
        Global Minimum             & 0.39         & 0.30           \\ 
        \bottomrule
    \end{tabular}
\end{table}

As we can see, there is a tradeoff that is to be made between parallelism and the consistency of elements generated. It is worth mentioning that the marginal increase in variance can be explained by a handful of sharp elements generated due to the new inner boundaries, and hence, the overall quality of the mesh still remains largely the same. The computational overhead of establishing parallel threads is superseded by its benefits, especially as the number of elements increases, during which a single agent begins to reach its limitations. Figure~\ref{fig:parallel_ability} illustrates the parallel meshing workflow that Table~\ref{tab:comparison_w_and_w/o_blockdecomposition} quantifies.

\subsection{Hole Decomposition \& 2D Meshing}
Here we demonstrate the hole decomposition agent, which is a purely visual agent trained on a dataset of $\approx$ 100,000 geometries with varying hole geometries. This has enabled a robust visual agent that splits a hole geometry into multiple filled geometries, which are then passed downstream into the block decomposition and finally, the meshing agent, which then stitches the individual blocks to create the final mesh. We now show the agent's effectiveness on varying geometries, as well as compare the results against the Soft Actor Critic approach alone. 

\begin{figure}[htpb]

     \centering
     \begin{minipage}{0.99\textwidth}
     \begin{subfigure}[b]{0.36\textwidth}
         \centering
         \includegraphics[width=\textwidth]{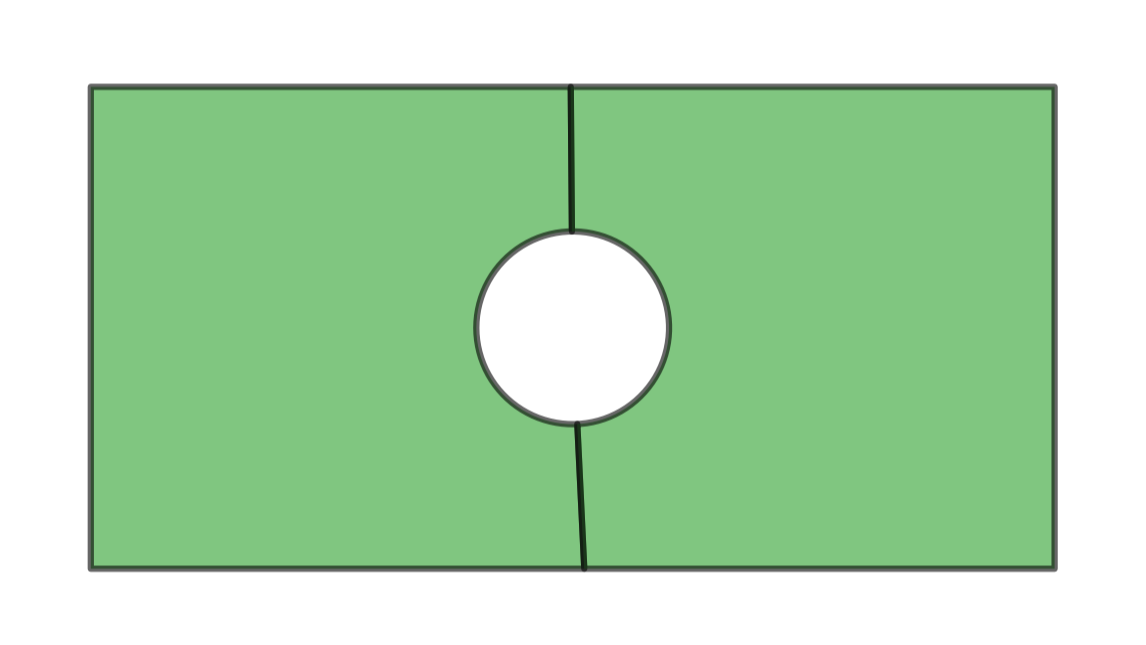}
     \end{subfigure}
     \hfill
     \begin{subfigure}[b]{0.36\textwidth}
         \centering
         \includegraphics[width=\textwidth]{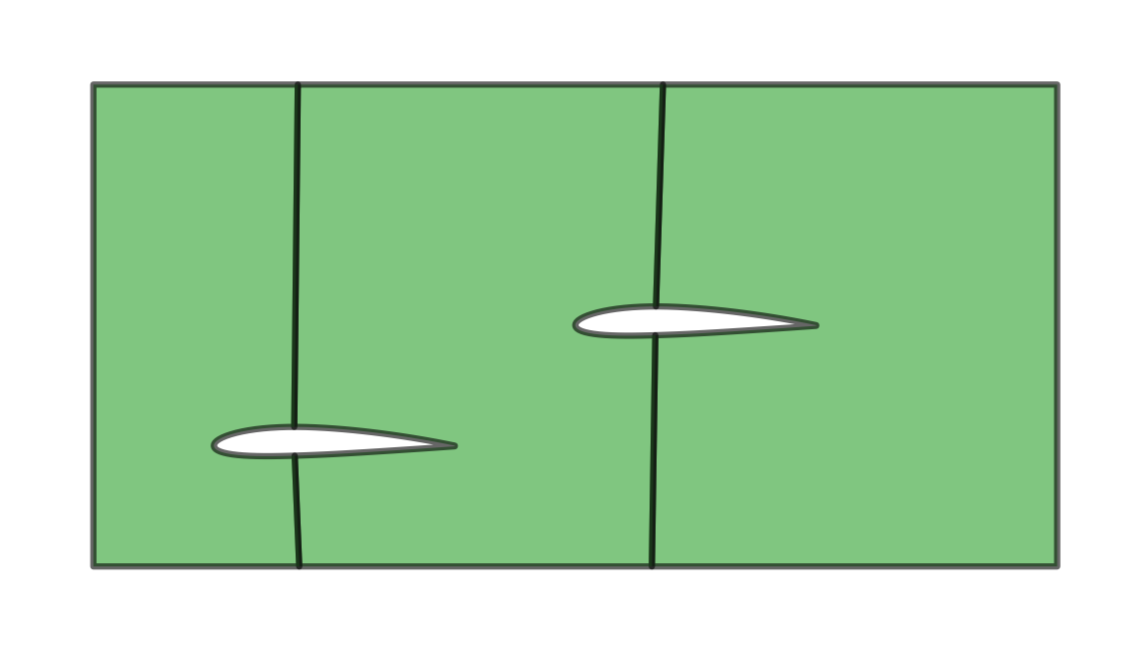}
     \end{subfigure}
     \hfill
     \begin{subfigure}[b]{0.21\textwidth}
         \centering
         \includegraphics[width=\textwidth]{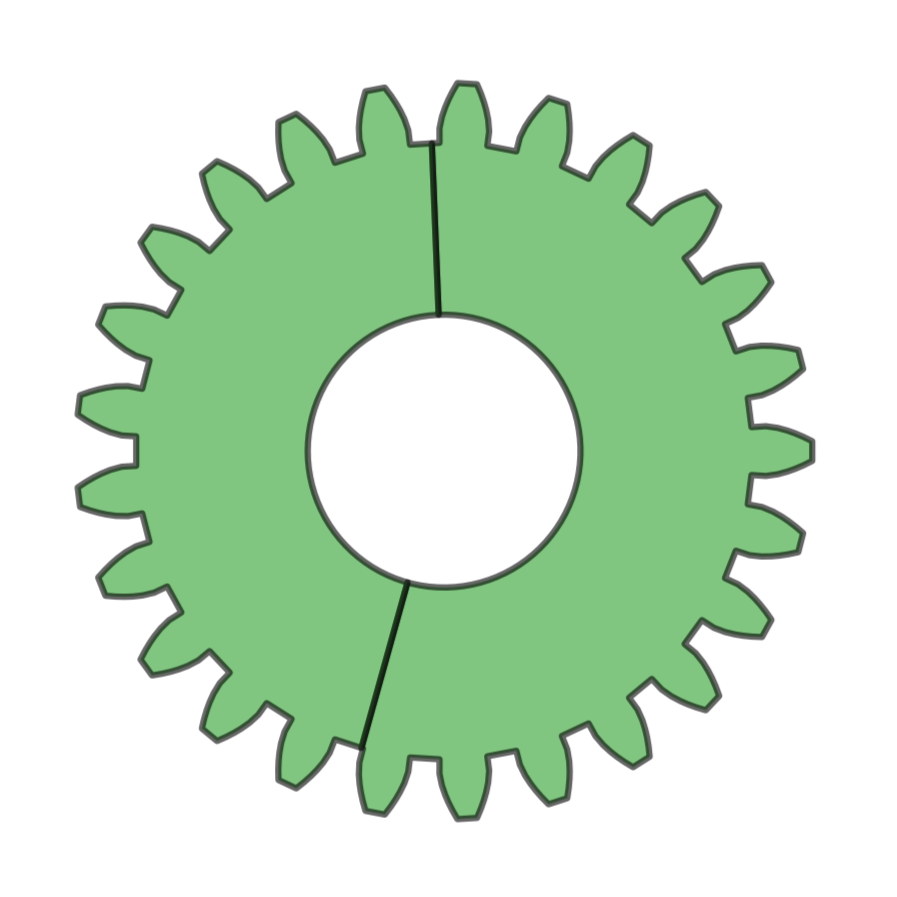}
     \end{subfigure}
     \hfill

     \vspace{10pt} 

     \begin{subfigure}[b]{0.36\textwidth}
         \centering
         \includegraphics[width=\textwidth]{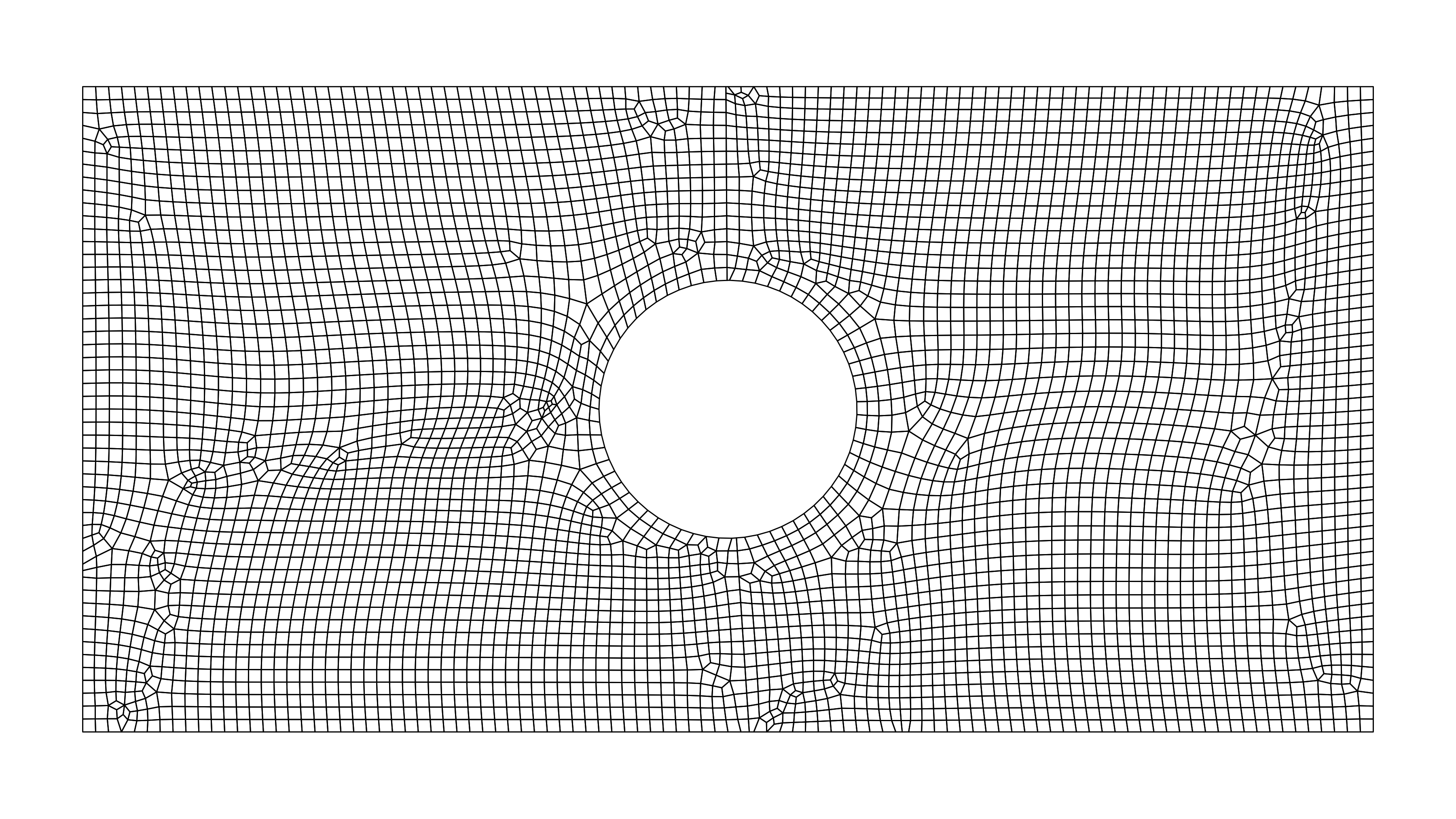}
     \end{subfigure}
     \hfill
     \begin{subfigure}[b]{0.36\textwidth}
         \centering
         \includegraphics[width=\textwidth]{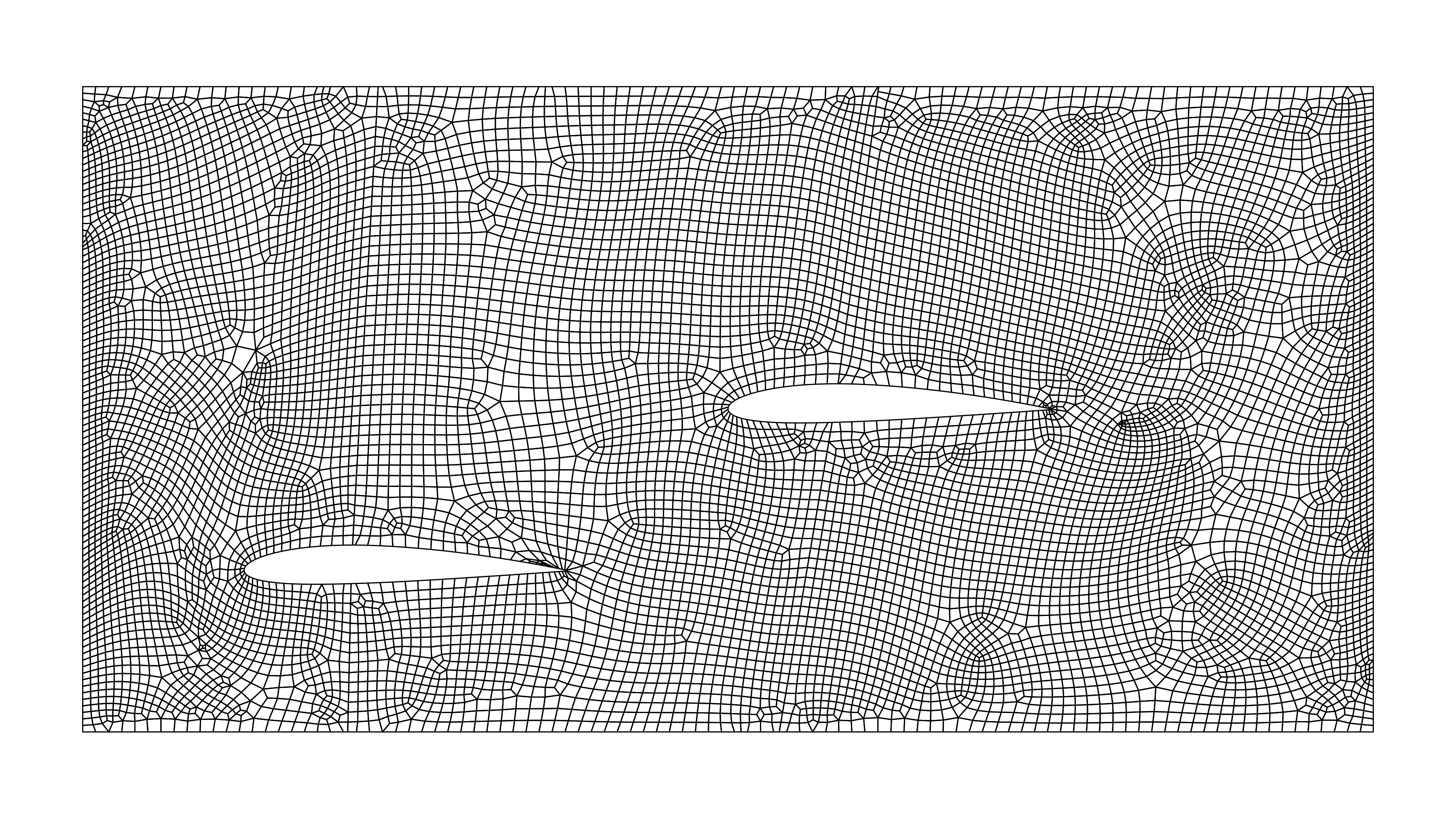}
     \end{subfigure}
     \hfill
     \begin{subfigure}[b]{0.21\textwidth}
         \centering
         \includegraphics[width=\textwidth]{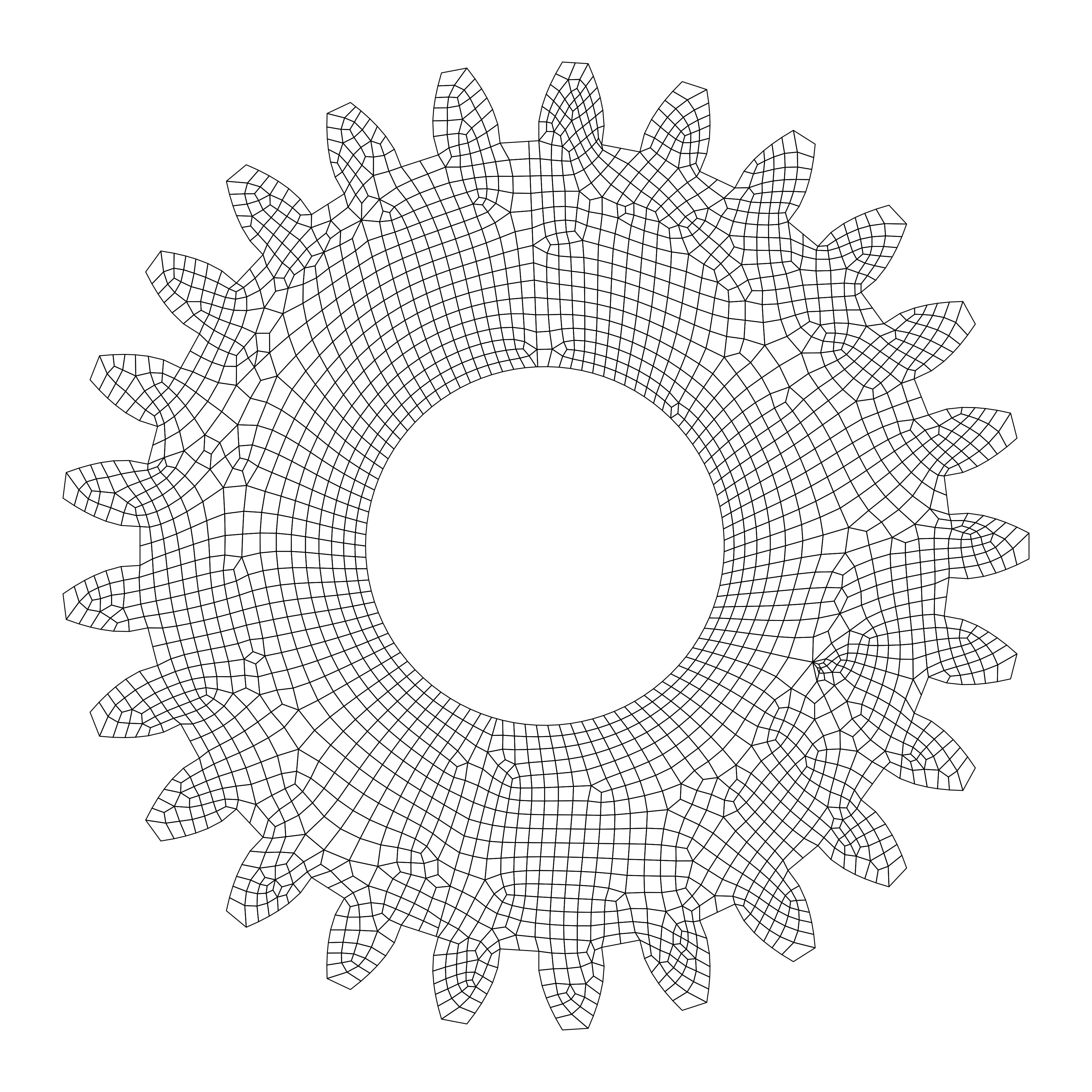}
     \end{subfigure}
     \end{minipage}
     \caption{End-to-end demonstration of the proposed multi-agent meshing pipeline on geometries containing holes. The hole decomposition agent first partitions the geometry into simply connected regions, followed by block decomposition and downstream quadrilateral mesh generation. The final meshes preserve conformity while maintaining high element quality near complex topological features.}
     \label{fig:holeDecompositionResult}
\end{figure}
As we can see with these examples, the hole decomposer has developed a good understanding of lines to cast to create a good split that doesn't cause any irregular elements. This allows for a good, clean split and ultimately a good final mesh. As you can see, the agent manages to cast good quality rays on 3 types of geometries that pertain to 3 different kinds of geometries. 

We now turn our attention to some quantitative comparisons of our pipeline $Dmsh$ against the current state of the art method for RL-based meshing: FreeMesh-RL. 

\begin{figure}[htpb]

     \centering
     \begin{minipage}{0.99\textwidth}
     \begin{subfigure}[b]{0.31\textwidth}
         \centering
         \includegraphics[width=\textwidth]{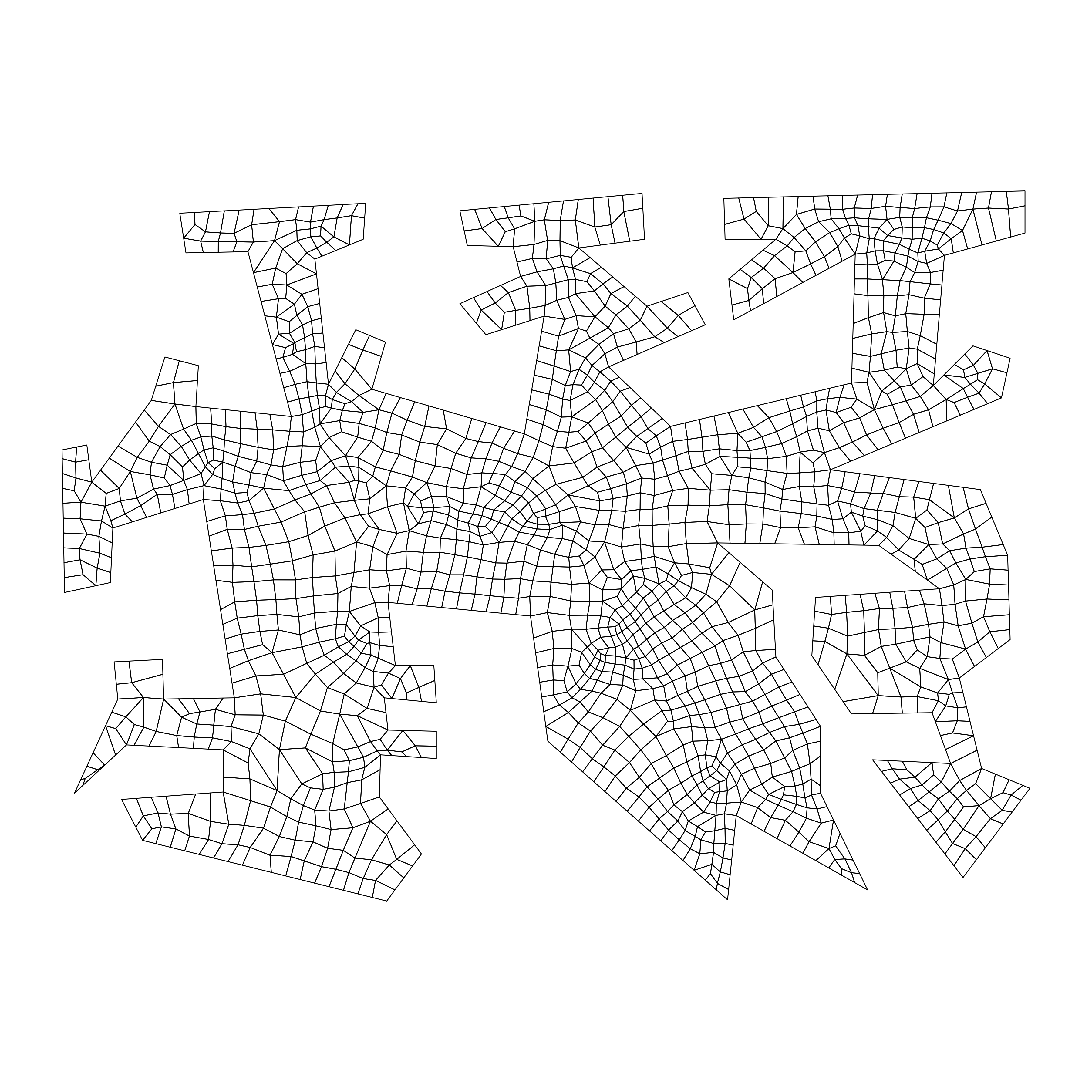}
     \end{subfigure}
     \hfill
     \begin{subfigure}[b]{0.31\textwidth}
         \centering
         \includegraphics[width=\textwidth]{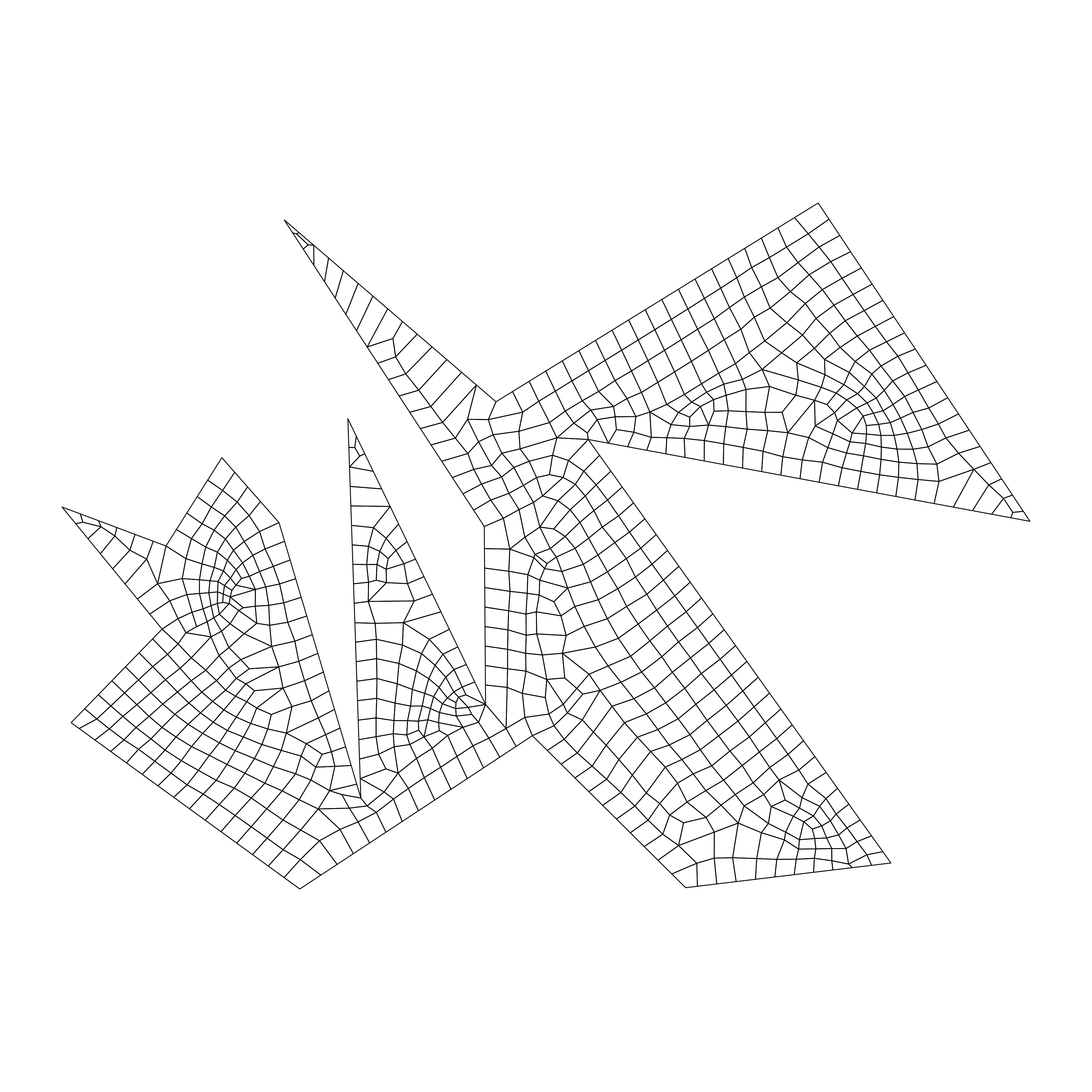}
     \end{subfigure}
     \hfill
     \begin{subfigure}[b]{0.31\textwidth}
         \centering
         \includegraphics[width=\textwidth]{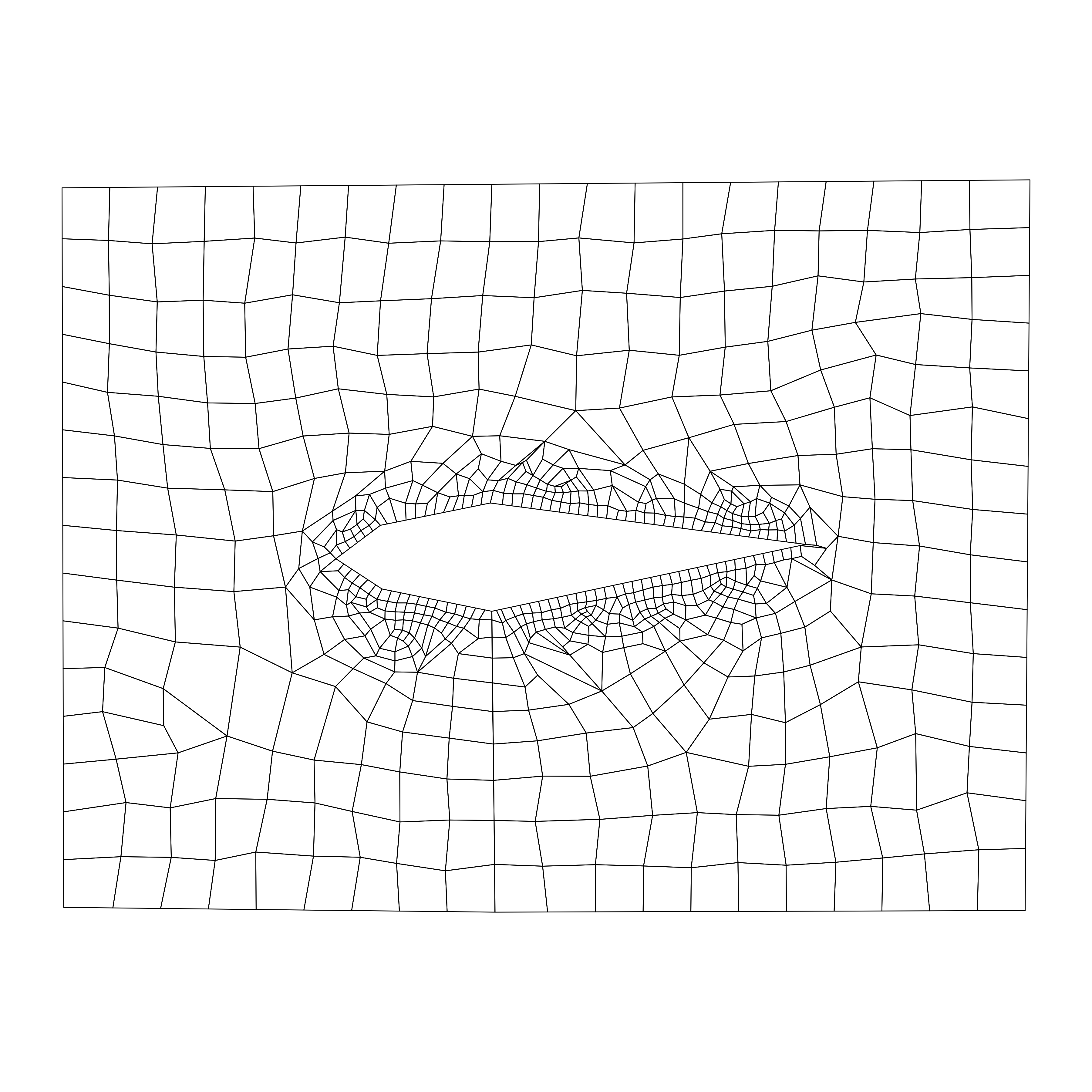}
     \end{subfigure}
     \hfill

     \vspace{-20pt} 

     \begin{subfigure}[b]{0.31\textwidth}
         \centering
         \includegraphics[width=\textwidth]{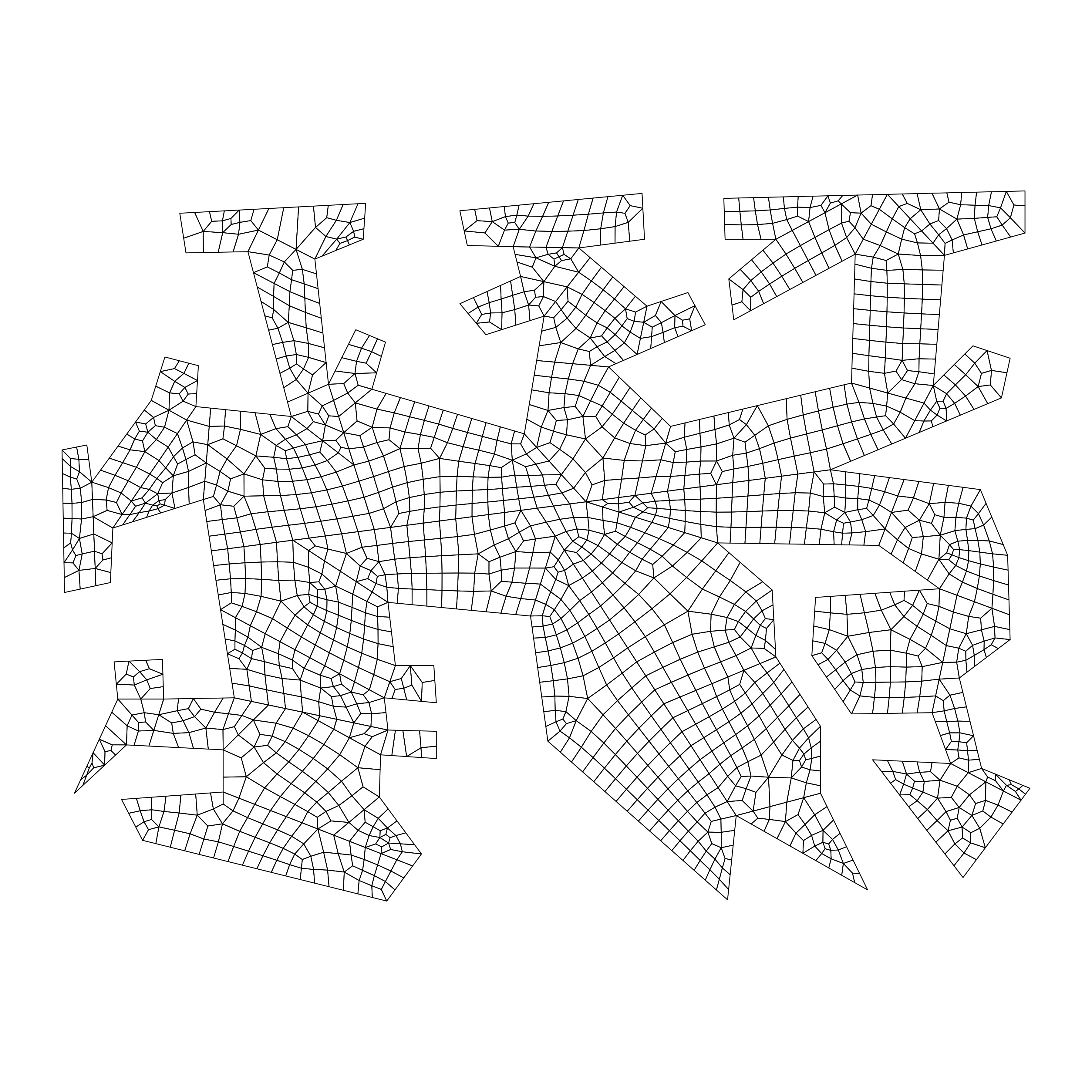}
         \caption{Domain 1}
     \end{subfigure}
     \hfill
     \begin{subfigure}[b]{0.31\textwidth}
         \centering
         \includegraphics[width=\textwidth]{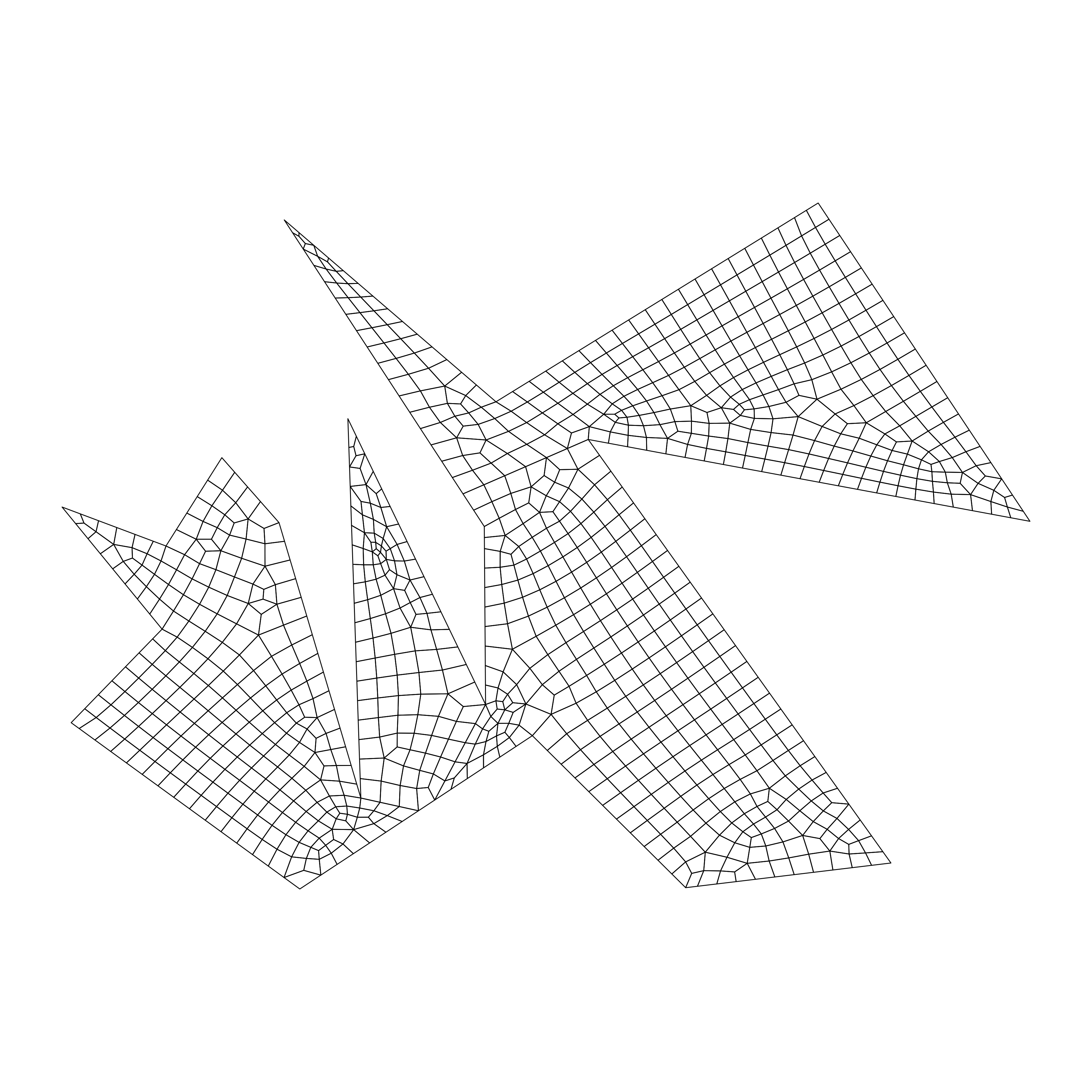}
         \caption{Domain 2}
     \end{subfigure}
     \hfill
     \begin{subfigure}[b]{0.31\textwidth}
         \centering
         \includegraphics[width=\textwidth]{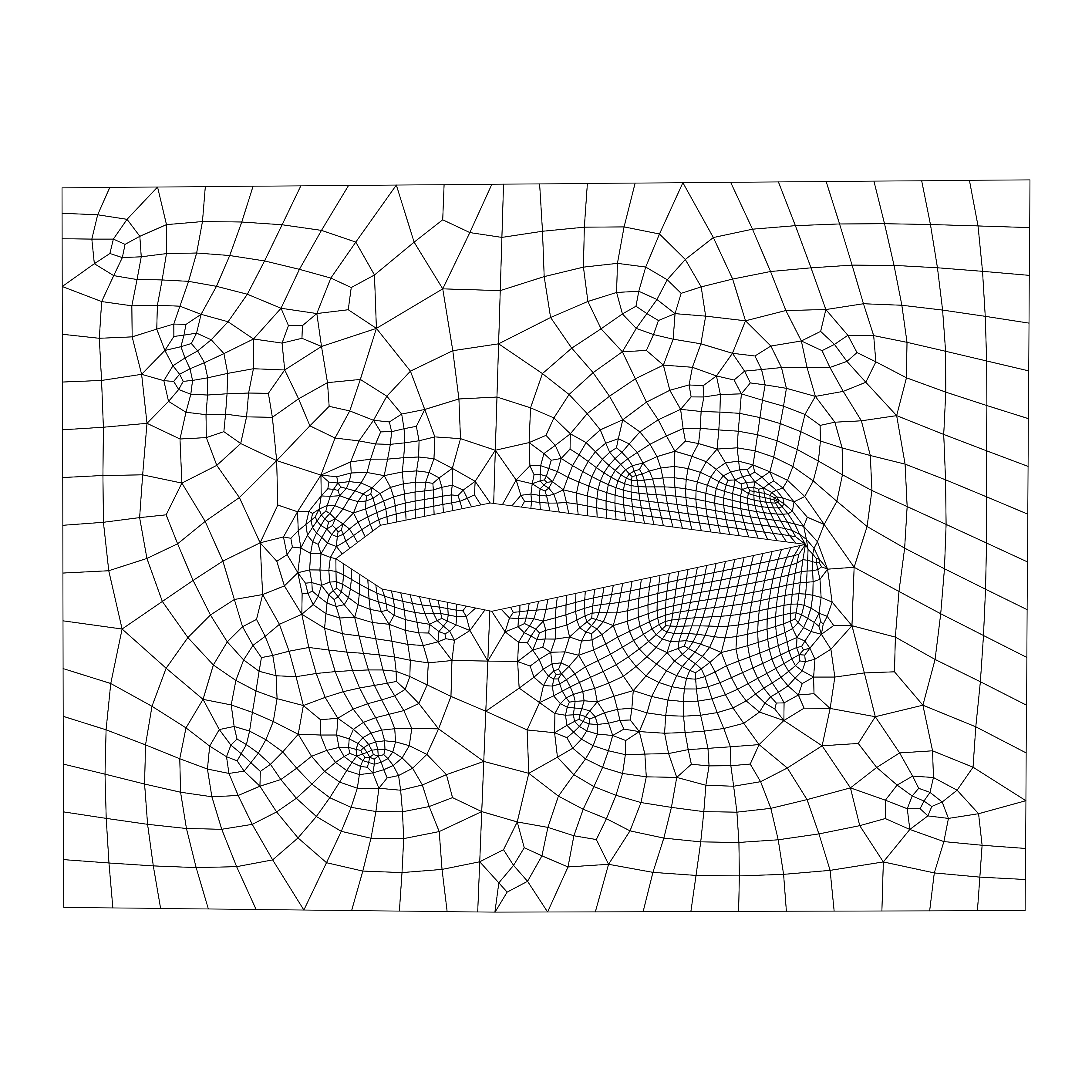}
         \caption{Domain 3}
     \end{subfigure}
     \end{minipage}
     \caption{Qualitative comparison between FreeMesh-RL (top row) and the proposed $Dmsh$ framework (bottom row) across three benchmark geometries. The proposed framework produces more uniform element distributions and improved handling of concave regions and hole geometries, particularly in domains requiring decomposition before advancing-front meshing.}
     \label{fig:comparo}
\end{figure}

\begin{table}[htpb]
    \centering
    \caption{Quantitative mesh quality comparison between FreeMesh-RL (current state-of-the-art reinforcement learning meshing method) and the proposed $Dmsh$ framework across three benchmark 2D geometries of increasing topological complexity: Domain 1 (highly non-convex star-like geometry), Domain 2 (multi-spike planar domain with sharp concavities), and Domain 3 (planar geometry containing an internal hole, modelled as a single convex block for FreeMesh-RL and as separate exterior/interior regions for $Dmsh$). Quality is measured using the Scaled Jacobian metric; higher average quality and global minimum, and lower variance, denote better mesh quality. Bold entries indicate the superior result for each metric and domain. Values are reported as FreeMesh-RL | $Dmsh$.}
    \label{tab:comparison}
    \begin{tabular}{lccc} 
        \toprule
        \textbf{Quality (Scaled Jacobian)} & \textbf{Domain 1 } & \textbf{Domain 2 }& \textbf{Domain 3 } \\ 
        \midrule
        Average Element Quality        & 0.8938 | \textbf{0.9054}    & 0.9187 | \textbf{0.9257}  & 0.8731 | \textbf{0.8887}\\ 
        Element Quality Variance       & 0.0157 | \textbf{0.0138}    & 0.0176 | \textbf{0.0127}  & 0.0331 | \textbf{0.0131}       \\ 
        Global Minimum                 & \textbf{0.0387} | 0.0001    & 0.0699 | \textbf{0.2792}  & \textbf{0.0606} | 0.0388     \\ 
        \bottomrule
    \end{tabular}
\end{table}

As we can see in Table~\ref{tab:comparison} and Figure~\ref{fig:comparo}, even in geometries where both FreeMesh-RL and $Dmsh$ are applied, $Dmsh$ manages to create better quality elements with lesser variance. In Domain 2, the advantage of the block decomposer present in the $Dmsh$ pipeline is especially seen with the global minimum being 4 times higher than FreeMesh-RL. 

It is worth mentioning that for Domain 3, the geometry is modeled as a single convex block wrapping around the cavity to be successfully meshed with FreeMesh-RL, whereas the geometries were independently passed as exterior and interior elements to $Dmsh$. This distinction allows $Dmsh$ to mesh geometries with multiple holes, as seen in Figure \ref{fig:holeDecompositionResult}. 

Moreover, $Dmsh$ is parallel by execution, allowing multiple agents to mesh the geometry at the same time, enabling exponentially faster inference as the domain size increases, thereby enabling scaling to tens of thousands of elements. This, combined with the contextual understanding of the hole and block decomposers, enables $Dmsh$ to mesh 3D surfaces. 

\subsection{Surface Meshes}

\begin{figure}[htpb]

     \centering
     \begin{minipage}{0.99\textwidth}
     \begin{subfigure}[b]{0.32\textwidth}
         \centering
         \includegraphics[width=\textwidth]{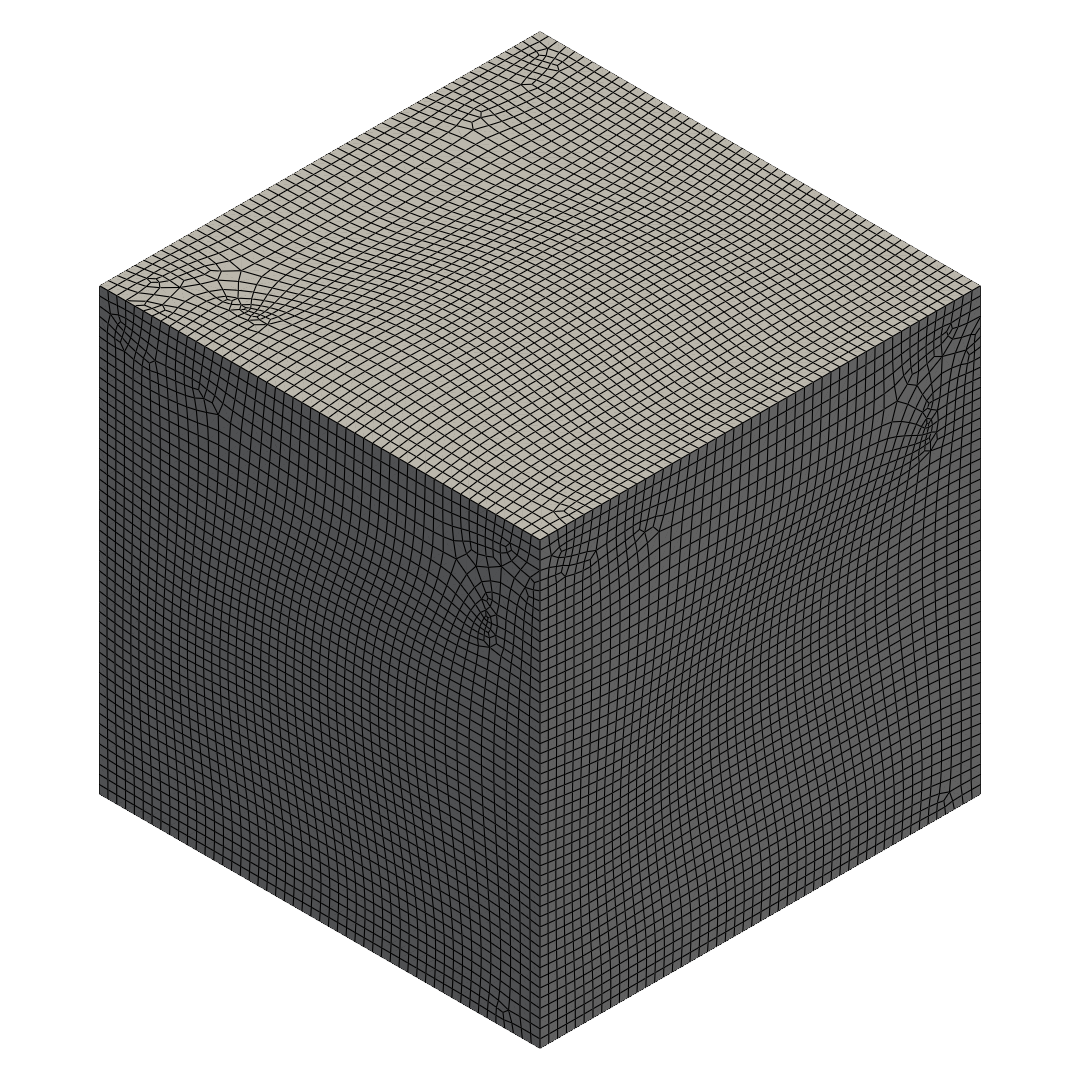}
     \end{subfigure}
     \hfill
     \begin{subfigure}[b]{0.32\textwidth}
         \centering
         \includegraphics[width=\textwidth]{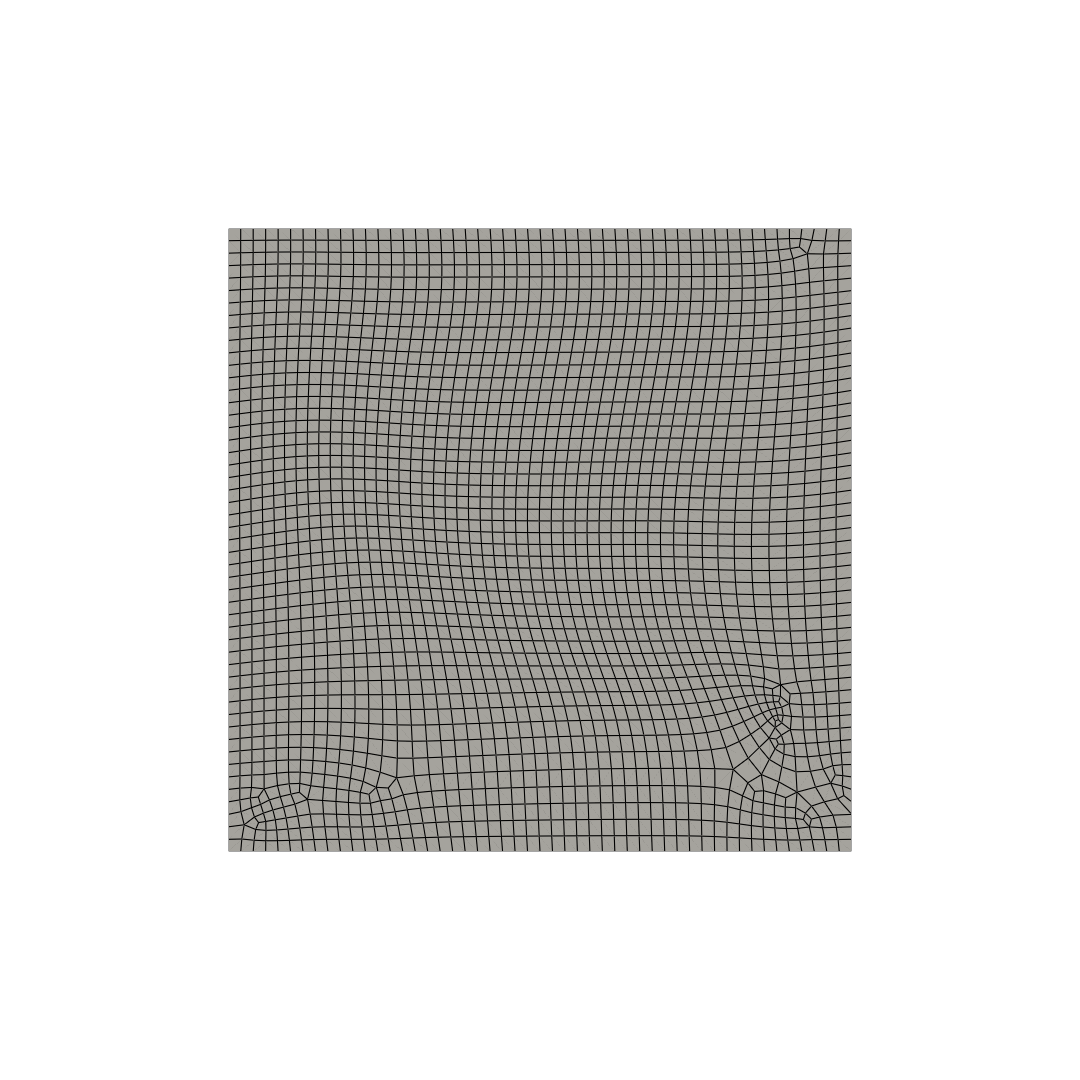}
     \end{subfigure}
     \hfill
     \begin{subfigure}[b]{0.32\textwidth}
         \centering
         \includegraphics[width=\textwidth]{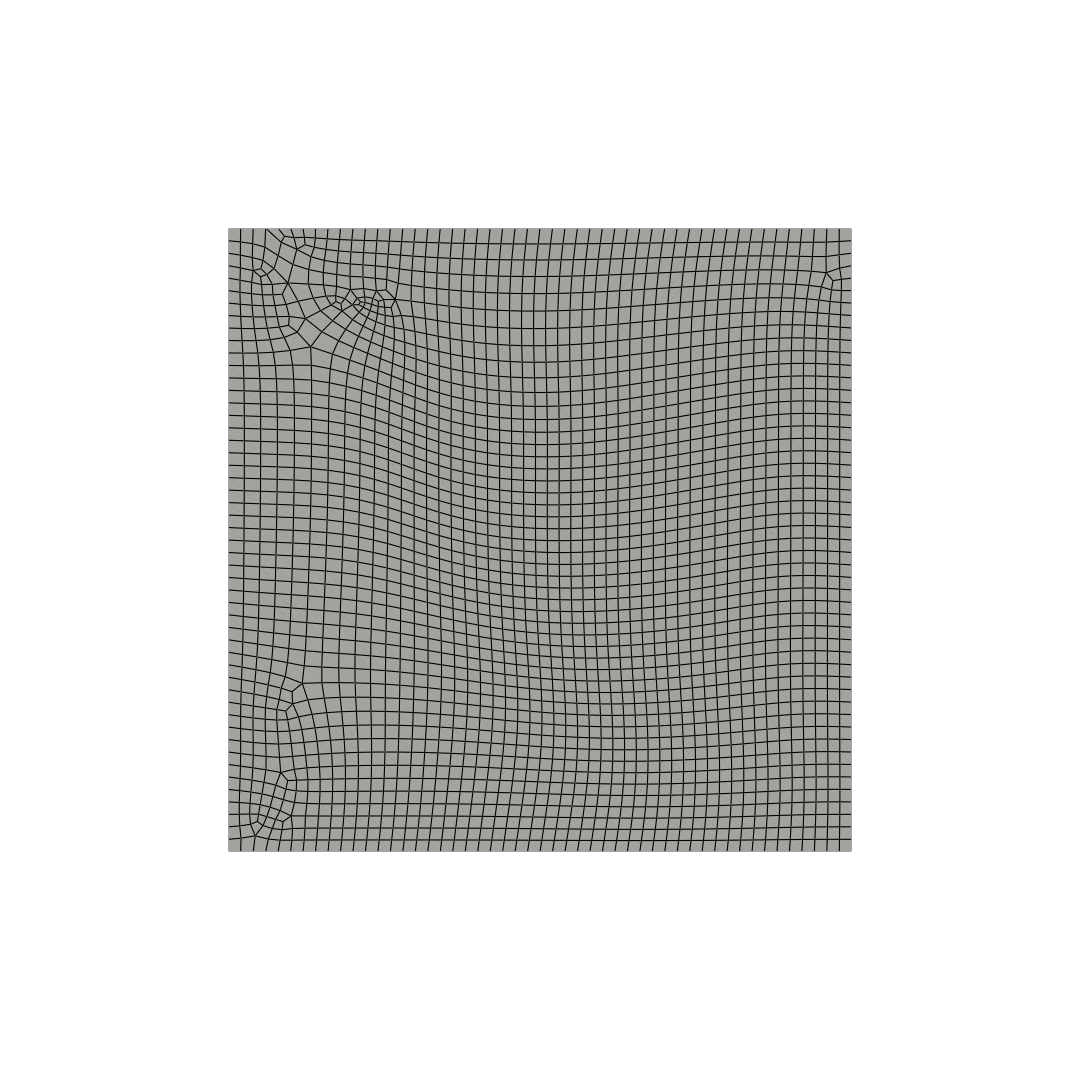}
     \end{subfigure}
     \hfill


     \begin{subfigure}[b]{0.32\textwidth}
     \centering
     \includegraphics[width=\textwidth]{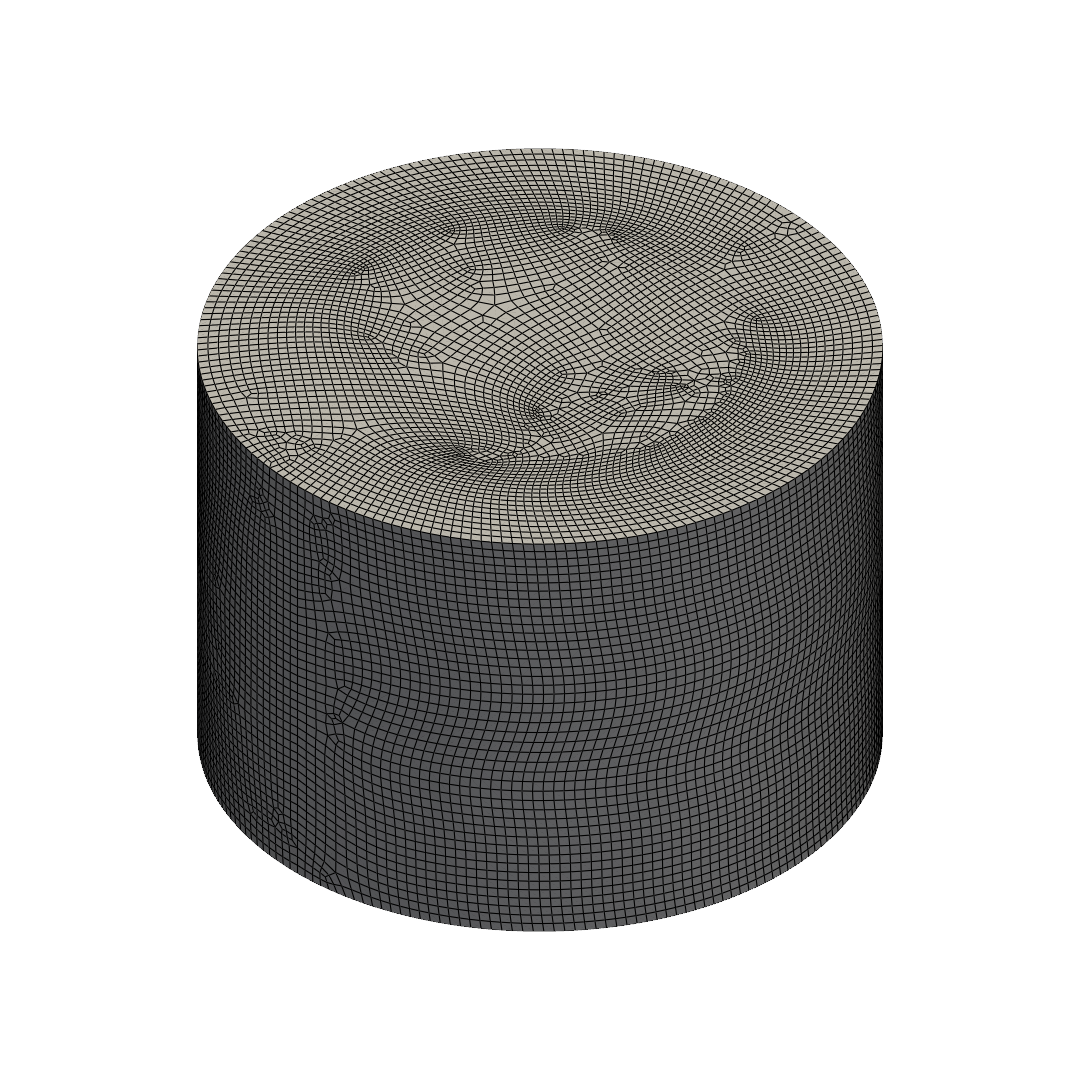}
     \end{subfigure}
     \hfill
     \begin{subfigure}[b]{0.32\textwidth}
         \centering
         \includegraphics[width=\textwidth]{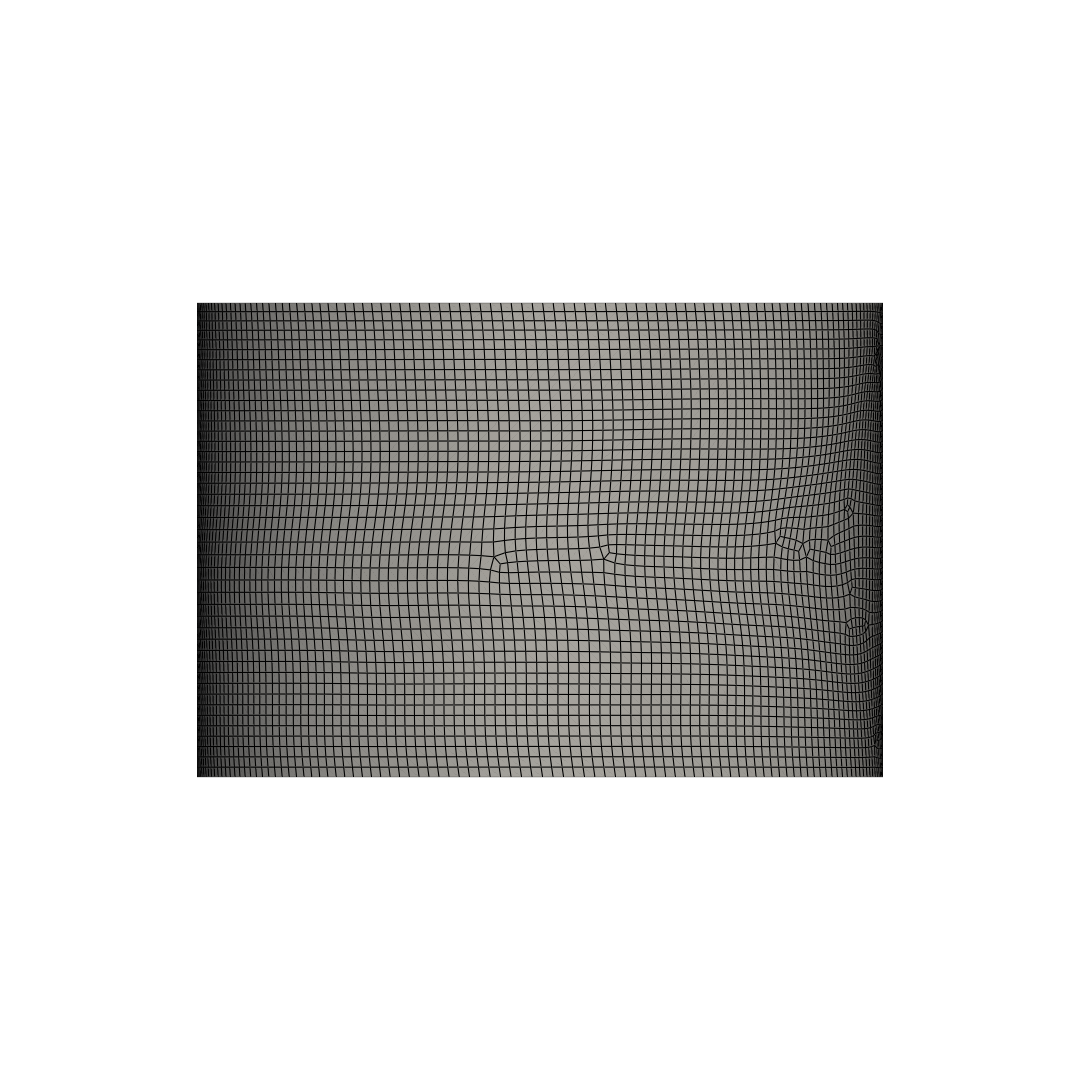}
     \end{subfigure}
     \hfill
     \begin{subfigure}[b]{0.32\textwidth}
         \centering
         \includegraphics[width=\textwidth]{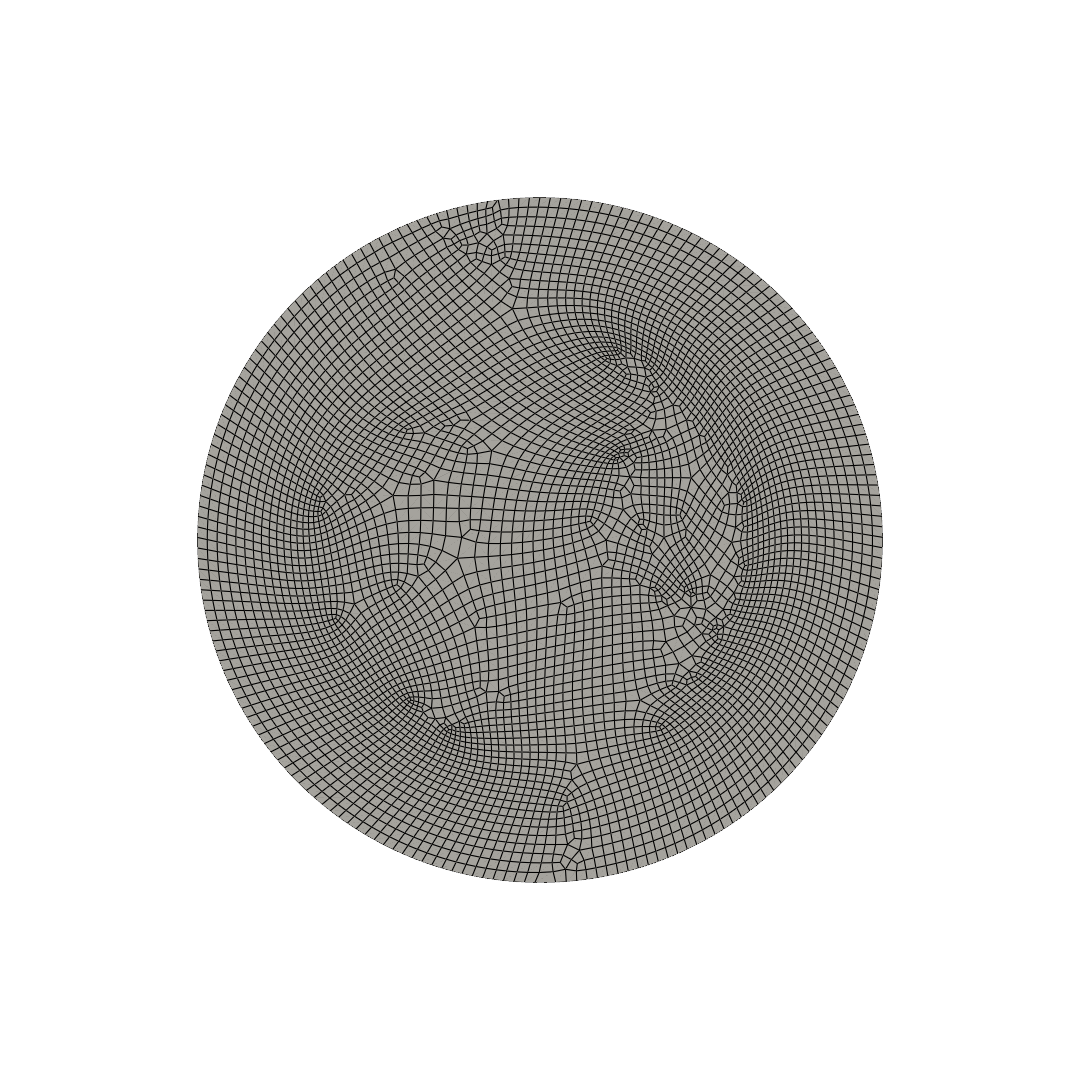}
     \end{subfigure}
     \hfill


     \begin{subfigure}[b]{0.32\textwidth}
         \centering
         \includegraphics[width=\textwidth]{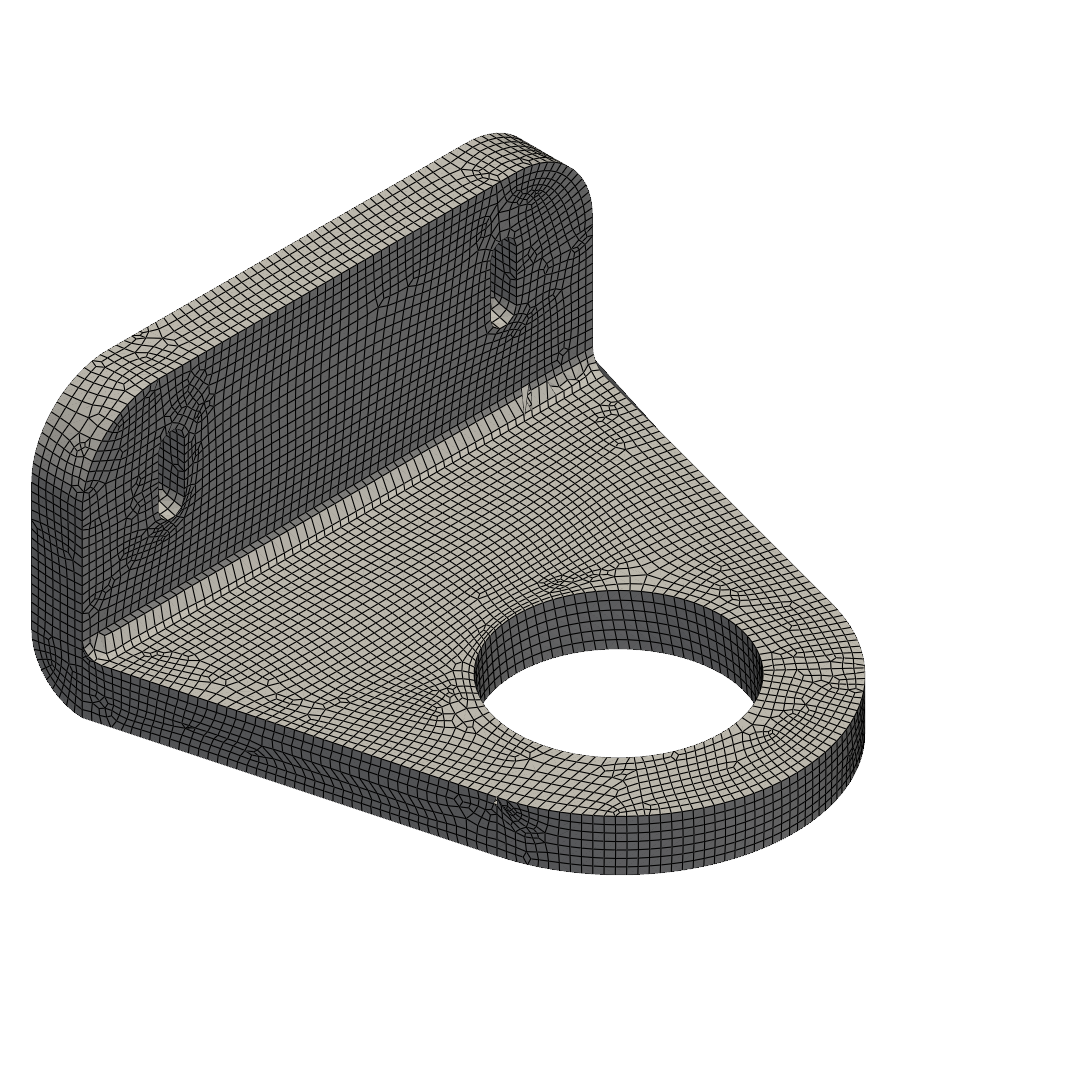}
     \end{subfigure}
     \hfill
     \begin{subfigure}[b]{0.32\textwidth}
         \centering
         \includegraphics[width=\textwidth]{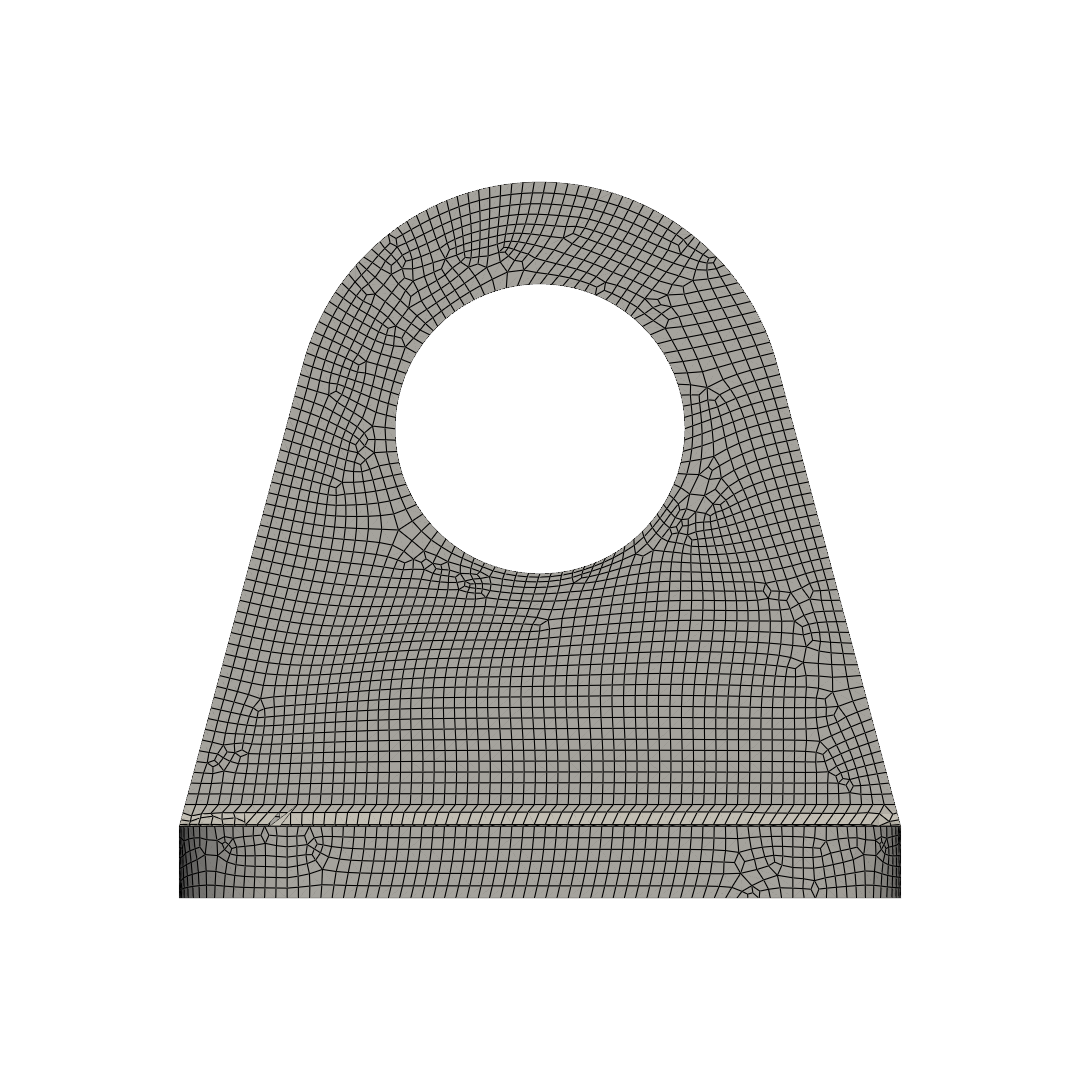}
     \end{subfigure}
     \hfill
     \begin{subfigure}[b]{0.32\textwidth}
         \centering
         \includegraphics[width=\textwidth]{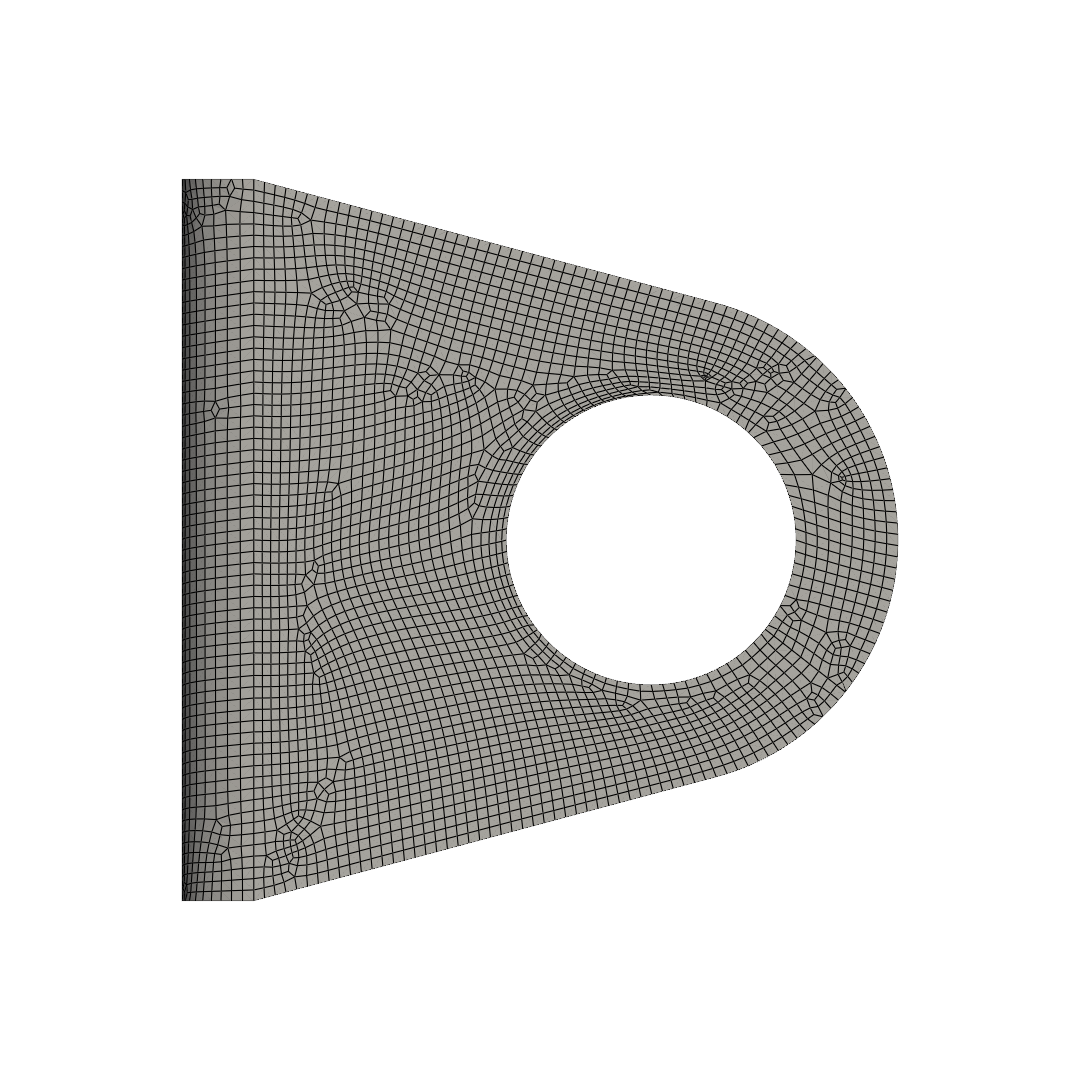}
     \end{subfigure}
     \end{minipage}
     \caption{Isometric, Top and Bottom views of a cube, a cylinder, and a bracket meshed by $Dmsh$}
     \label{fig:SurfaceMeshView}
\end{figure}

We extend these agents to a 3D surface by building a surface meshing pipeline that individually meshes each face of the 3D surface. To achieve this, we use the Python plug-in of the OpenCASCADE library, enabling us to mesh the parametric form of the boundary representation and convert the new nodes and node mappings back into spatial coordinates.  We demonstrate this by meshing a cube (14,586 elements), a cylinder (39,366 elements), and a bracket (12,244 elements) (see Figure~\ref{fig:SurfaceMeshView}). 

The robustness of $$Dmsh$$ allows it to mesh almost any shape, with or without holes, in the parametric space, thereby creating an all-quad surface mesh. We tabulate the salient additions our framework provides on top of existing research.
\begin{table}[htpb]
    \centering
    \caption{Feature-level capability comparison between FreeMesh-RL (current state-of-the-art RL-based all-quadrilateral meshing method) and the proposed $$Dmsh$$ framework. }
    \vspace{0.2cm}
    \begin{tabular}{lcc}
        \hline
        \textbf{Capability} & \textbf{FreeMesh-RL} & \textbf{$Dmsh$} \\
        \hline
        All Quad & $\checkmark$ & $\checkmark$ \\
        Single Hole & $\checkmark^*$ & $\checkmark$ \\
        Block Decomposition & $\times$ & $\checkmark$ \\
        Multiple Holes & $\times$ & $\checkmark$ \\
        Parallel Meshing & $\times$ & $\checkmark$ \\
        Surface Meshing & $\times$ & $\checkmark$ \\
        \hline
        *boundary modification required
    \end{tabular}
    \label{tab:freemeshrl_Dmsh_comparison}
\end{table}

\section{Conclusion}
In this work, we presented a fully autonomous framework for the generation of all-quadrilateral surface meshes, formulated as a multi-agent reinforcement learning problem involving three interacting agents. The proposed approach integrates computer vision-based representations to guide block and hole decomposition, effectively approximating heuristic strategies traditionally employed by human experts. Quantitative evaluation on established benchmark geometries demonstrates that the framework achieves performance comparable to state-of-the-art techniques, while offering improved generality and adaptability across diverse geometric configurations, including convex regions. Future research directions include scaling the learning framework through expanded training datasets and introducing a Mixture of Experts (MoE) paradigm, wherein specialized agents could be trained to handle distinct classes of geometric features and topologies. Additionally, extending the formulation to volumetric mesh generation and incorporating parallelization strategies, such as domain decomposition, is expected to significantly enhance computational efficiency and scalability, enabling the generation of large-scale meshes with millions of elements.

\section*{Acknowledgment}
Anirudh Kalyan would like to gratefully acknowledge the support provided by the DAAD KOSPIE scholarship for the academic year 2025–2026.

\bibliographystyle{elsarticle-num} 
\bibliography{myRefs}
\end{document}